\documentclass[11pt]{article}
\usepackage{amsfonts}

\usepackage{CJK}
\usepackage{amsfonts,amssymb,amsmath,mathrsfs}

\usepackage{color,latexsym,amsfonts}
 \newcommand{\qed}{\hfill\rule{2mm}{3mm}\vspace{4mm}}

 \newtheorem{theorem}{Theorem}[section]
 \newtheorem{lemma}[theorem]{Lemma}
 \newtheorem{corollary}[theorem]{Corollary}
 \newtheorem{proposition}[theorem]{Proposition}
 
 \newtheorem{Definition}[theorem]{Definition}
 \newtheorem{remark}[theorem]{Remark}
 \newtheorem{condition}[theorem]{Assumption}

 \def\blemma{\begin{lemma}\sl{}\def\elemma{\end{lemma}}}
 \def\btheorem{\begin{theorem}\sl{}\def\etheorem{\end{theorem}}}
 
 \def\bdefinition{\begin{Definition}\sl{}\def\edefinition{\end{Definition}}}
 \def\bproposition{\begin{proposition}\sl{}\def\eproposition{\end{proposition}}}
 \def\bremark{\begin{remark}\sl{}\def\eremark{\end{remark}}}
 \def\bcondition{\begin{condition}\sl{}\def\econdition{\end{condition}}}

 \def\beqlb{\begin{eqnarray}}\def\eeqlb{\end{eqnarray}}
 \def\beqnn{\begin{eqnarray*}}\def\eeqnn{\end{eqnarray*}}

 \def\mcr{\mathscr}\def\mbb{\mathbb}\def\mbf{\mathbf}

 \def\<{\langle}\def\>{\rangle}

 \def\ar{&\!\!}\def\qqquad{\qquad\qquad}

 \def\eqref#1{{\rm(\ref{#1})}}

 \def\proof{\noindent{\it Proof.~}}\def\qed{\hfill$\Box$\medskip}

\begin{document}
\noindent{(Draft: 2016/09/19)}

\bigskip\bigskip

\centerline{\Large\bf Pathwise uniqueness for an SPDE}

\smallskip

\centerline{\Large\bf with H\"older continuous coefficient}

\smallskip

\centerline{\Large\bf driven by $\alpha$-stable noise}

\bigskip

\centerline{Xu Yang\footnote{Supported by NSFC (Nos.~11401012 and 11602003) and NSF of Ningxia
(No.~NZ15095)}
and Xiaowen Zhou
\footnote{Supported by  NSERC grant 249554-2011}}

\bigskip

\centerline{\small School of Mathematics and Information Science,
Beifang University of Nationalities,}

\centerline{\small Yinchuan 750021, People's Republic of China.}

 \smallskip

\centerline{\small Department of Mathematics and Statistics, Concordia University,}
\centerline{\small 1455 de Maisonneuve Blvd. West, Montreal, Quebec, H3G 1M8, Canada}

\bigskip

\centerline{{\tt xuyang@mail.bnu.edu.cn} and {\tt
xiaowen.zhou@concordia.ca}}

\bigskip\bigskip

{\narrower{\narrower

\noindent{\bf Abstract.} In this paper we study the pathwise
uniqueness of nonnegative solution for the following stochastic
partial differential equation with H\"older continuous coefficient:
 \beqnn
\frac{\partial X_t(x)}{\partial t}=\frac{1}{2} \Delta X_t(x)
+G(X_t(x))+H(X_{t-}(x))
\dot{L}_t(x),\quad t>0,
~x\in\mbb{R},
 \eeqnn
where $\dot{L}$ denotes an $\alpha$-stable white noise on
$\mbb{R}_+\times \mbb{R}$ without negative jumps, $G$ satisfies a
condition weaker than Lipschitz, and $H$ is nondecreasing and
$\beta$-H\"older continuous for $1<\alpha<2$ and $0<\beta<1$.

For  $G\equiv0$ and $H(x)=x^\beta$, in Mytnik (2002) a weak solution
to the above stochastic heat equation was constructed and the
pathwise uniqueness of the nonnegative solution was left as an open
problem. In this paper we give an affirmative answer to this problem
for certain values of $\alpha$ and $\beta$. In particular, for
$\alpha\beta=1$ the solution to the above equation is the
density of a super-Brownian motion with $\alpha$-stable branching
(see also Mytnik (2002)) and  our result leads to its pathwise
uniqueness for $1<\alpha<\sqrt{5}-1$.

The local H\"older continuity
of the solution is also obtained in this paper for fixed time $t>0$.

\bigskip


\bigskip


\par}\par}

\section{Introduction}\label{sec_1}

\setcounter{equation}{0}

\subsection{Background and motivation}

It was proved by Konno and Shiga (1988) \cite{KoS88} and by Reimers (1989) \cite{Rei89}
that for an arbitrary initial measure the one-dimensional binary branching super-Brownian motion has a
jointly continuous density that is a random field
$\{X_t(x):t>0,x\in\mbb{R}\}$ satisfying the following
continuous-type stochastic partial differential equation (SPDE):
 \beqlb\label{1.02}
\frac{\partial}{\partial t} X_t(x)
 =
\frac{1}{2}\Delta X_t(x) + \sqrt{X_t(x)}\dot{W}_t(x), \quad t>0,
~x\in\mbb{R},
 \eeqlb
where $\Delta$ denotes the one-dimensional Laplacian operator and
$\{\dot{W}_t(x): t>0, x\in \mbb{R}\}$ denotes the derivative of a
space-time Gaussian white noise.

The weak uniqueness of solution to the above stochastic heat equation
follows from that of a martingale problem for super-Brownian motion.
The pathwise uniqueness of nonnegative solution to SPDE \eqref{1.02}
remained open even though it had been studied by many authors. The
main difficulty comes from the non-Lipschitz diffusion coefficient.
Progresses have been made in considering modified forms of the SPDE.
When the random field $\{W_t(x): t>0, x\in \mbb{R}\}$ is colored in
space and white in time, the strong uniqueness of  nonnegative
solution to the SPDE was obtained by Mytnik \textit{et al.} (2006)
\cite{MPS06} and further work can be found in Rippl and Sturm (2013)
\cite{RS13} and in Neuman (2014) \cite{N14}. Xiong (2013)
\cite{Xio13} proved the pathwise uniqueness of a  SPDE satisfied by
the ``distribution function" of the  super-Brownian on $\mathbb{R}
$. When $\{W_t(x): t>0, x\in \mbb{R}\}$ is a space-time Gaussian
white noise, the solutions are allowed to take both positive and
negative values and $\sqrt{X_t(x)}$ is replaced by
$\sigma(t,x,X_t(x))$ in SPDE \eqref{1.02},  the pathwise uniqueness
of the solution was proved by Mytnik and Perkins (2011) \cite{MyP11}
for $\sigma(\cdot,\cdot,u)$ with H\"older continuity in $u$ of index
$\beta_0>3/4$. Further work can be found in Mytnik and Neuman (2015)
\cite{MN15}. Recently, some negative results were obtained. When
$\sqrt{X_t(x)}$ is replaced by $|X_t(x)|^{\beta_1}$ in the SPDE
\eqref{1.02}, Burdzy \textit{et al.} (2010) \cite{BMP10} showed a
non-uniqueness result for $0<\beta_1<1/2$ and Mueller \textit{et
al.} (2014) \cite{MMyP2014} proved a non-uniqueness result for
$1/2\le\beta_1<3/4$.

Mytnik (2002) \cite{M02} considered the following jump-type SPDE and
constructed a weak solution:
  \beqlb\label{1.001}
\frac{\partial X_t(x)}{\partial t}=\frac{1}{2} \Delta X_t(x)+X_{t-}(x)^\beta
\dot{L}_t(x),\quad t>0,
~x\in\mbb{R},
 \eeqlb
where $0<\beta<1$ and for $1<\alpha<2$, $\dot{L}$ is a one sided
$\alpha$-stable white noise on $\mbb{R}_+\times \mbb{R}$ without
negative jumps. Put $p:=\alpha\beta<2$. The solution to
\eqref{1.001} with $p=1$ is the density of a super-Brownian motion
with $\alpha$-stable branching and the weak uniqueness of the
solution holds; see \cite[Theorem 1.6]{M02}. But for the other
values of $p$ the uniqueness  for \eqref{1.001}
was left as an open problem; see \cite[Remark 1.7]{M02}. During the
past ten years there have been a number of very interesting results
on the solution of SPDE \eqref{1.001} for  $p=1$. In particular,
Mytnik and Perkins (2003) \cite{MP03} showed that the solution has a
continuous modification at any fixed time. Fleischmann \textit{et
al. }(2010) \cite{FMW10} showed that this continuous modification is
locally H\"older continuous with index $\eta_c:=2/\alpha-1$, and
Fleischmann \textit{et al. }(2011) \cite{FMW11} further showed that
it is H\"older continuous with index
$\bar{\eta}_c:=(3/\alpha-1)\wedge1$ at any given spatial point. A
more precise analysis on the regularity of the solution was given in
Mytnik and Wachtel (2015) \cite{MW15}. He \textit{et al.} (2014)
\cite{HLY12} showed that another jump-type and \eqref{1.001} related
SPDE on the distribution-function-valued process
 is pathwise unique. For $p\neq1$, the uniqueness of
solution (including the weak uniqueness) to SPDE \eqref{1.001} and
the regularities of the solution $X_t(\cdot)$ at a fixed time $t$
are also left as open problems; see \cite[Remark 5.9]{M02}.

In this paper we want  to establish the pathwise uniqueness of
nonnegative solution for \eqref{1.001}.
For this purpose we consider a  SPDE more general than
\eqref{1.001}:
  \beqlb\label{1.01}
\frac{\partial X_t(x)}{\partial t}=\frac{1}{2} \Delta X_t(x)
+G(X_t(x))+H(X_{t-}(x))
\dot{L}_t(x),\quad t>0,
~x\in\mbb{R},
 \eeqlb
where $G$ and $H$ are non-negative functions and satisfy the following conditions:\\
(C1) (Linear growth condition) There is a constant $C$ so that
 \beqnn
0\le G(x)\le C(x+1),\qquad x\ge0.
 \eeqnn
(C2) Function $G$ is continuous and there is a non-decreasing
and concave function $r_0$ on $[0,\infty)$ so that $r_0(0)=0$,
$\int_{0+}r_0(z)^{-1}dz=\infty$ and
 \beqnn
\mathop{\rm sgn}(x-y)(G(x)-G(y))\le r_0(|x-y|),\qquad x,y\ge0,
 \eeqnn
where $\mathop{\rm sgn}(x):=1_{(0,\infty)}(x)-1_{(-\infty,0)}(x)$.

{\noindent}(C3) ($\beta$-H\"older continuity) There is a constant
$C$ so that
 \beqnn
|H(x)-H(y)|\le C|x-y|^\beta,\qquad x,y\ge0.
 \eeqnn
(C4) $H(x)$ is a nondecreasing function.

There have been many results on SPDEs driven by stable noises; see
e.g. \cite{AWZ98,saint,Mueller98,AW00,CZ11}. In \cite{AWZ98}, the
existence and uniqueness were established for solutions of
stochastic reaction equations driven by Poisson random measures. The
existence of weak solutions and pathwise uniqueness for stochastic
evolution equations driven by L\'evy processes can be found  in
\cite{CZ11}. It was also shown in \cite{CZ11} that the pathwise
uniqueness holds if the coefficient of the L\'evy noise satisfies a
condition weaker than Lipschitz continuity but stronger than
H\"older continuity.  The main results of \cite{AW00,saint,AWZ98} are the
strong existence and uniqueness of solution to \eqref{1.01} with
general L\'evy noise $\dot{L}$ and Lipschitz continuous coefficient
$H$. In this paper we use a
Yamada-Watanabe argument that is different from
\cite{CZ11}, and we consider a special L\'evy noise of stable noise
without negative jumps. The stable noise had not been treated in the
above mentioned papers although technically it is not hard to extend
their results in that direction under the Lipschitz condition on
$H$. One contribution of this paper is that we are able to get rid of
the Lipschitz condition on $H$ since we only need it to be H\"older
continuous.

The SPDE \eqref{1.001} was studied in Mueller (1998) \cite{Mueller98}
 for $\alpha$-stable noise $\dot{L}$ with  $0<\alpha<1$.
We also refer to Peszat and Zabczyk (2007) \cite{PZ07} for early work
on SPDEs driven by L\'evy noises.

\textit{Throughout this paper, we always assume that $1<\alpha<2$,
$0<\beta<1$ and the solutions to \eqref{1.001} and \eqref{1.01} are
nonnegative. Our goal} is to establish the pathwise uniqueness of
solution to \eqref{1.01} under conditions (C1)--(C4) and further
restrictions on $\alpha$ and $\beta$. In particular, for $p=1$ we
show that the pathwise uniqueness holds for $1<\alpha<\sqrt{5}-1$.
To prove the pathwise uniqueness we need to show a local H\"older
continuity of the solution at fixed time $t>0$, which  also extends
the regularity results for super-Brownian motion with
$\alpha$-stable branching obtained in Fleischmann \textit{et al.
}(2010) \cite{FMW10}.

To continue with the introduction we present some notation. Let
$\mathscr{B}(\mbb{R})$ be the set of Borel functions on $\mbb{R}$.
Let $B(\mbb{R})$ denote the Banach space of bounded Borel functions
on $\mbb{R}$ furnished with the supremum norm $\|\cdot\|$. We use
$C(\mbb{R})$ to denote the subset of $B(\mbb{R})$ of bounded
continuous functions. For any integer $n\ge 1$ let $C^n(\mbb{R})$ be
the subset of $C(\mbb{R})$ of functions with bounded continuous
derivatives up to the $n$th order. Let $C_c^n(\mbb{R})$ be the
subset of $C^n(\mbb{R})$ of functions  with compact supports. We use the
superscript ``+'' to denote the subsets of positive elements of the
function spaces, e.g., $B(\mbb{R})^+$. For $f,g\in
\mathscr{B}(\mbb{R})$ write $\<f,g\>= \int_{\mbb{R}}f(x)g(x)dx$ whenever
it exists. Let $M(\mbb{R})$ be the space of finite Borel measures on
$\mbb{R}$ endowed with the weak convergence topology. For $\mu\in
M(\mbb{R})$ and $f\in B(\mbb{R})$ we also write $\mu(f) = \int f
d\mu$.

Equation \eqref{1.01} is a formal SPDE that is understood in the
following sense: For any $f\in\mathscr{S}(\mbb{R})$,  the (Schwartz)
space of rapidly decreasing and infinitely differentiable functions on
$\mbb{R}$,
 \beqlb\label{1.1}
\<X_t,f\>
 \ar=\ar
X_0(f)+ \frac12\int_0^t\<X_s,f''\>ds
+\int_0^tds\int_{\mbb{R}}G(X_s(x))f(x)dx \cr
 \ar\ar
+ \int_0^t\int_{\mathbb{R}}H(X_{s-}(x)) f(x)L(ds,dx),\qquad t\ge0,
 \eeqlb
where $X_0\in M(\mbb{R})$ and $L(ds,dx)$ is a one-sided $\alpha$-stable white noise on
$\mbb{R}_+\times \mbb{R}$ without negative jumps.


\bdefinition\label{t1.10} SPDE \eqref{1.1} has a weak
solution $(X,L)$ with initial value $X_0\in M(\mbb{R})$ if there is a
pair $(X,L)$ defined on the same filtered probability space
$(\Omega,\mathscr{F},\mathscr{F}_t,\mbf{P})$ satisfying  the
following conditions.

\noindent(i) $L$ is an $\alpha$-stable white noise on
$\mbb{R}_+\times \mbb{R}$ without negative jumps.

\noindent(ii) The two-parameter nonnegative process
$X=\{X_t(x):t>0,x\in\mbb{R}\}$ is progressively measurable on
$\mbb{R}_+\times \mbb{R}\times\Omega$, and
$\{1_{\{t=0\}}X_0(dx)+1_{\{t>0\}}X_t(x)dx:t\ge0\}$ is a
$M(\mbb{R})$-valued c\'adl\'ag process.

\noindent(iii) For each $f\in \mathscr{S}(\mbb{R})$, $(X,L)$
satisfies \eqref{1.1}.
\edefinition
The definition of this kind of $\alpha$-stable white noise
$L(ds,dx)$ and Definition \ref{t1.10} can be found in \cite{M02}.

\subsection{The main results and approaches}

Given $t>0$, we say $\tilde{X}_t$ is a {\it continuous modification}
of $X_t$ if $\tilde{X}_t(x)$ is continuous in $x$  and
$\mbf{P}\{\tilde{X}_t(x)=X_t(x)$ for all $x\}=1$.
The following first theorem gives the local H\"older continuity (in
the spatial variable) for the
continuous modification of the solution to \eqref{1.1}.

\btheorem\label{t1.1} (Local H\"older
continuity) For any fixed $t>0$, $X_t$ has a continuous
modification $\tilde{X}_t$. Moreover, for each
$\eta<\eta_c:=\frac{2}{\alpha}-1$, with probability one the
continuous modification $\tilde{X}_t$ is locally H\"older continuous
of exponent $\eta$, i.e. for any compact set
$\mathbb{K}\subset\mbb{R}$,
 \beqlb\label{2.19}
\sup_{x,z\in \mbb{K},x\neq
z}\frac{|\tilde{X}_t(x)-\tilde{X}_t(z)|}{|x-z|^\eta}<\infty,  \,\,\,
\mbf{P}\mbox{-a.s.}
 \eeqlb
In addition, if $\beta<1/\alpha+(\alpha-1)/2$, then for each $T>0$ and
subsequence $\{n':n'\ge1\}$ of $\{n:n\ge1\}$, we have
 \beqlb\label{2.17}
\liminf_{n'\to\infty} \frac{1}{2^{n'}}\sum_{k=1}^{2^{n'}}
\sup_{x,z\in \mbb{K},x\neq z}
\frac{|\tilde{X}_{{n'}_kT}(x)-\tilde{X}_{{n'}_kT}(z)|}
{|x-z|^\eta}<\infty, \,\,\, \mbf{P}\mbox{-a.s.},
 \eeqlb
where ${n'}_k:=\frac{k}{2^{n'}}$ for $1\le k\le2^{n'}$.
 \etheorem

\bremark Theorem \ref{t1.1} gives an answer to
\cite[Conjecture 1.5]{FMW10} when the fractional Laplacian $\Delta_\alpha$ is the
Laplacian operator $\Delta$ and the function $g$ there is replaced by $H$.
It also gives an answer to the open problem of
\cite[Remark 5.9]{M02}.
 \eremark

\bcondition\label{c1.1} For $p:=\alpha\beta>1$,   there is a constant
 $q>\frac{3p}{3-\alpha}$ so that for any weak
solution $(X,L)$ to \eqref{1.1}  it holds that
 \beqnn
\mbf{P}\Big\{\int_0^tds\int_{\mbb{R}}X_s(x)^qdx<\infty \mbox{ for
all } t>0\Big\}=1.
 \eeqnn
\econdition

\btheorem\label{t4.1} (Pathwise uniqueness)
Suppose that conditions (C1)--(C4) hold, and
that
 \beqlb\label{4.17}
2(\alpha-1)/(2-\alpha)^2<\beta<1/\alpha+(\alpha-1)/2.
 \eeqlb
We also assume that Assumption \ref{c1.1} holds for $p>1$. If
$(X,L)$ and $(Y,L)$, with $X_0=Y_0\in M(\mbb{R})$, are  two weak
solutions to equation \eqref{1.1} defined on the same filtered
probability space $(\Omega,\mathscr{F},\mathscr{F}_t,\mbf{P})$,
 then with probability one, for
each $t>0$ we have
 \beqlb\label{4.16}
X_t(x)=Y_t(x),\qquad\lambda_0\mbox{-a.e. }x,
 \eeqlb
where $\lambda_0$ denotes the Lebesgue measure on $\mbb{R}$.
\etheorem

\bremark
(i) Since we assume that $\beta\in(0,1)$, it follows from
the first inequality of \eqref{4.17} that the theorem makes sense for $\alpha\in(1,3-\sqrt{3})$.

(ii) Theorem \ref{t4.1} gives an affirmative answers to the open problem of
\cite[Remark 1.7]{M02} for $\alpha$ and $\beta$ satisfying \eqref{4.17}.

(iii) If $p=1$, inequality \eqref{4.17} is equivalent to
$1<\alpha<\sqrt{5}-1$. So, for super-Brownian motion, i.e. $G\equiv0$, $H(x)=x^\beta$ and
$p=1$, Theorem \ref{t4.1} also leads to  the pathwise uniqueness of
\eqref{1.001} for $1<\alpha<\sqrt{5}-1$, which is a key  result of
this paper.

(iv) We stress here that in Theorem \ref{t4.1}, Assumption \ref{c1.1} is not needed if $0<p=\alpha\beta\le1$.

(v) If $G\equiv0$ and $H(x)=x^\beta$, then SPDE \eqref{1.1} has a
weak solution satisfying Assumption \ref{c1.1} by \cite[Proposition
5.1]{M02} and the proof of \cite[Theorem 1.5]{M02}.

(vi) The non-negativity assumption on $H$ and $G$, which makes the proof a bit
simple and may be needed on the existence of the solution to the SPDE,
is in fact not necessary.
\eremark

 To prove the uniqueness we need a local H\"older
continuity of the solution at fixed time $t>0$ (Theorem \ref{t1.1}). For super-Brownian
motion, the proof for the local H\"older continuity of $X_t(x)$ is
based on the following equation from Fleischmann \textit{et al.
}(2010) \cite{FMW10}:
 \beqlb\label{1.03}
\<X_t,f\>=X_0(f)+\frac12\int_0^t\<X_s,f''\>ds
+\int_0^t\int_0^\infty\int_{\mbb{R}}f(x)zM(ds,dz,dx),
 \eeqlb
where $M(ds,dz,dx)$ denotes a compensated Poisson random measure  on $(0,\infty)^2\times
\mbb{R}$ with compensator $\hat{M}(ds,dz,dx)=dsm_0(dz)X_s(x)dx$ for
measure $m_0(dz):=c_0z^{-1-\alpha}1_{\{z>0\}}dz$ with
$c_0:={\alpha(\alpha-1)}/{\Gamma(2-\alpha)}$ and
Gamma function $\Gamma$. Equation \eqref{1.03} is established for
super-Brownian motion. But for the other cases, the solution to
\eqref{1.1} may not be a density of super-Brownian motion and we can
not obtain the equivalent of equation \eqref{1.03}. So, inspired by
Dawson and Li (2006, 2012)\cite{DLi06,DLi12}, we reformulate
\eqref{1.1} as the following SPDE in Proposition \ref{t1.9}:
 \beqlb\label{1.4}
\<X_t,f\>
 \ar=\ar
X_0(f)+
\frac12\int_0^t\<X_s,f''\>ds+\int_0^tds\int_{\mbb{R}}G(X_s(x))f(x)dx
\cr
 \ar\ar
+ \int_0^t\int_0^\infty\int_{\mbb{R}}\int_0^{H(X_{s-}(u))^\alpha }
zf(u) \tilde{N}_0(ds,dz,du,dv),
 \eeqlb
where $f\in \mathscr{S}(\mbb{R})$ and
$\tilde{N}_0(ds,dz,du,dv)$ is a compensated Poisson random measure
on $(0,\infty)^2\times\mbb{R}\times(0,\infty)$ with intensity
$dsm_0(dz)dudv$. By modifying the proof of \cite[Theroem 1.2(a)]{FMW10} and
using \eqref{1.4}, we can obtain Theorem \ref{t1.1}. Notice that
$\eta_c\uparrow 1$ as $\alpha\downarrow 1$, which  is quite
different from  that of a continuous-type SPDE whose local H\"older
index is typically smaller than $\frac12$. This observation is key to
proving the pathwise uniqueness.

We now outline our approach. By  an infinite-dimensional version of
the Yamada-Watanabe argument for ordinary stochastic differential
equations (see Mytnik \textit{et al. }(2006)), showing the pathwise
uniqueness is reduced to showing that the analogue of the local time
term is zero; see the proofs of Theorem \ref{t4.1} and Lemma
\ref{t4.4}. That is to show that
 \beqlb\label{1.06}
\mbf{E}\{I_5^{m,n}(t\wedge\tau_k)\} \to0
 \eeqlb
as $m,n\to\infty$, where
 \beqnn
I_5^{m,n}(t\wedge\tau_k) := \int_0^{t\wedge\tau_k}ds\int_0^\infty
m_0(dz)\int_{\mbb{R}}
\<D_n(\<U_s,\Phi_{\cdot}^m\>,zV_s(y)\Phi_{\cdot}^m(y)),\Psi_s\>dy
 \eeqnn
 for $\tau_k:=\gamma_k\wedge\sigma_k$, and
$\gamma_k$ and $\sigma_k$ are two stopping times to be defined later
in (\ref{def_gamma}) and (\ref{def_sigma}), respectively. Here $U_s$
is the difference of two weak solutions to \eqref{1.1}, $V_s$ denotes
the difference of compositions of these two solutions into function
$H$, respectively, $\Psi$ is a test function, $\Phi_x^m$ is a
mollifier, $D_n(y,z):=\phi_n(y+z)-\phi_n(y)-z\phi_n'(y)$, $\phi_n$
($supp(\phi_n'')\subset (a_n,a_{n-1})$ and $a_n\downarrow0$) is the
function satisfying $\phi_n(x)\to |x|$  from  Yamada-Watanabe's
proof. To prove \eqref{1.06}, we divide
$\mbf{E}\{I_5^{m,n}(t\wedge\tau_k)\}$  into two terms
 \beqnn
I_{5,1}^{m,n,k,i}(t)
 \ar:=\ar
\mbf{E}\Big\{\int_0^{t\wedge\gamma_k}ds\int_0^{1/i}
m_0(dz)\int_{\mbb{R}}
\<D_n(\<U_s,\Phi_{\cdot}^m\>,zV_s(y)\Phi_{\cdot}^m(y)),\Psi_s\>dy\Big\},\cr
 \eeqnn
 and
 \beqnn
I_{5,2}^{m,n,k,i}(t)
 \ar:=\ar
\mbf{E}\Big\{\int_0^{t\wedge\sigma_k}ds\int_{1/i}^\infty
m_0(dz)\int_{\mbb{R}}
\<D_n(\<U_s,\Phi_{\cdot}^m\>,zV_s(y)\Phi_{\cdot}^m(y)),\Psi_s\>dy\Big\},
 \eeqnn
so that $\mbf{E}\{I_5^{m,n}(t\wedge\sigma_k)\} \le
I_{5,1}^{m,n,k,i}(t)+I_{5,2}^{m,n,k,i}(t)$ for all $i\ge1$.

Using the fact $\phi''_n\le 2(na_n)^{-1}$, we can show that
$I_{5,1}^{m,n,k,i}(t)$ goes to zero as $m,n,i\to\infty$ in a
dependent way (see Lemma \ref{t4.9}). So, the difficult part  is to
show that $I_{5,2}^{m,n,k,i}(t)$ goes to zero as $m,n,i\to\infty$.
To this end we use the local H\"older continuity of the solutions
and the monotonicity of $H$ to estimate $I_{5,2}^{m,n,k,i}(t)$,
which is elaborated in the following. The proof is inspired by an
argument of Mytnik and Perkins (2011) \cite{MyP11},
for fixed $s,m$ and $x$, denote by $x_{s,m}\in[-1,1]$  a value
satisfying
 \beqnn
|\tilde{V}_s(x-\frac{x_{s,m}}{m})|
=\inf_{y\in[-1,1]}|\tilde{V}_s(x-\frac{y}{m})|,
 \eeqnn
where $\tilde{V}_s$ and $\tilde{U}_s$ are the continuous
modifications of $V_s$ and $U_s$, respectively. The key to proving
that $I_{5,2}^{m,n,k,i}(t)$ goes to zero is to split it into two
terms again, where one term is bounded from above by
 \beqnn
I_{5,2,1}^{m,n,k,i}(t)
  \ar:=\ar
\mbf{E}\Big\{\int_0^{t\wedge\sigma_{k}}ds
\int_{-K}^K\Psi_s(x)dx\int_{1/i}^\infty zm_0(dz)
\int_{-1}^1\Phi(y)dy \cr \ar\ar\quad
\times\int_0^1\Big|\tilde{D}_n\big(\<\tilde{U}_s,\Phi_x^m\>,mzh\tilde{V}_s(x-\frac{y}{m})\big)
[V_s(x-\frac{y}{m})-\tilde{V}_s(x-\frac{x_{s,m}}{m})]\Big|dh \Big\},
 \eeqnn
and the other term is bounded from above by
 \beqnn
I_{5,2,2}^{m,n,i}(t)
 \ar:=\ar
\mbf{E}\Big\{\int_0^tds \int_{-K}^K\Psi_s(x)dx\int_{1/i}^\infty
zm_0(dz) \int_{-1}^1\Phi(y)dy\cr
 \ar\ar\quad
\times\int_0^1\tilde{D}_n\big(\<\tilde{U}_s,\Phi_x^m\>,mzh\tilde{V}_s(x-\frac{y}{m})\big)
|\tilde{V}_s(x-\frac{x_{s,m}}{m})|
1_{\{\tilde{V}_s(x-\frac{x_{s,m}}{m})\neq0\}}dh \Big\},
 \eeqnn
where $\tilde{D}_n(y,z)=\phi_n'(y+z)-\phi_n'(y)$.

The local H\"older continuity of the solutions is used to estimate
$I_{5,2,1}^{m,n,k,i}(t)$ and the nondecreasingness of $H$ is used to
estimate $I_{5,2,2}^{m,n,i}(t)$. Observe that for fixed $s$, the
continuous modification $\tilde{X}_s$ of the weak solution to
\eqref{1.1} satisfies
 \beqlb\label{1.04}
 \ar\ar
\sup_{|x|\le K,|y|\vee|v|\le1}|\tilde{X}_s(x-\frac{y}{m})
-\tilde{X}_s(x-\frac{v}{m})|^{\beta}  \cr
 \ar\ar\qqquad\le
(2/m)^{\eta\beta} \sup_{|x|\le K,|y|\vee|v|\le1,y\neq v}
\frac{|\tilde{X}_s(x-\frac{y}{m})
-\tilde{X}_s(x-\frac{v}{m})|^{\beta}} {|y/m-v/m|^{\eta\beta}},
 \eeqlb
where $K>0$ and $0<\eta<\eta_c=2/\alpha-1$. So, it is natural to
apply the H\"older continuity of $x\mapsto\tilde{X}_s(x)$ to find a
collection of suitable stopping times $(\sigma_{k})_{k\ge1}$ so that
$\lim_{k\to\infty}\sigma_{k}=\infty$ almost surely, and using the
$\beta$-H\"older continuity of $H$ (condition (C3)), the term
$I_{5,2,1}^{m,n,k,i}(t)$ can be bounded by
$m^{-\eta\beta}i^{\alpha-1}$ which tends to zero as $m,n,i$ jointly
go to infinity in a certain way. It is hard to show that the
supremum or integral with respect to $s\in(0,T]$ on the right hand side
of \eqref{1.04} is finite. To this end, the time $\sigma_{k}$ is
chosen so that a Riemann type ``integral'' of the right hand side of
\eqref{1.04} over $s\in[0,\sigma_{k}]$ is finite. One can find the
details in the Step 1 of the proof of Lemma \ref{t4.8}.

Concerning the second term $I_{5,2,2}^{m,n,i}(t)$, if
$\tilde{V}_s(x-\frac{x_{s,m}}{m})\neq0$, then the function
$[-1,1]\ni y\mapsto\tilde{V}_s(x-\frac{y}{m})$ is bounded away from
zero. The nondecreasingness of $H$ (that is condition (C4) and this
condition will only be used here) ensures that
$\tilde{V}_s(x-\frac{y}{m})$ and $\tilde{U}_s(x-\frac{y}{m})$ always
have the same sign, which means
$\tilde{D}_n\big(\<\tilde{U}_s,\Phi_x^m\>,mzh\tilde{V}_s(x-\frac{y}{m})\big)=0$
for $|\<\tilde{U}_s,\Phi_x^m\>|\ge a_{n-1}$ for all $z,h\ge0$ (here
we use the fact $supp(\phi_n'')\subset (a_n,a_{n-1})$ and
$a_n\downarrow0$). Thus by the $\beta$-H\"older continuity of $H$
(condition (C3)),
 \beqnn
|\tilde{V}_s(x-\frac{x_{s,m}}{m})|\le
C|\tilde{U}_s(x-\frac{x_{s,m}}{m})|^\beta \le
C|\<\tilde{U}_s,\Phi_x^m\>|\leq Ca_{n-1}^\beta,
 \eeqnn
 which implies $I_{5,2,2}^{m,n,i}(t)$
also converges to zero as $m,n,i\to\infty$ under certain conditions
of $\alpha$ and $\beta$ (see the details in the Step 2 of the proof
of Lemma \ref{t4.8}).

\subsection{Comments on the main results with general $G$ and $H$}
The main results, Theorems
\ref{t1.1} and \ref{t4.1}, also hold if functions $G(x)$ and $H(x)$
are replaced by $G(t,x,y)$ and $H(t,x,y)$, respectively, as in
\cite{MyP11,MPS06}. More specifically, we can consider an SPDE more
general than \eqref{1.01}:
 \beqlb\label{1.05}
\frac{\partial X_t(x)}{\partial t}=\frac{1}{2} \Delta X_t(x)
+G(t,x,X_t(x))+H(t,x,X_{t-}(x))
\dot{L}_t(x),\quad t>0,
~x\in\mbb{R},
 \eeqlb
where $G$ and $H$ satisfy the following growth and continuity
conditions:\\
(1) The mapping
$(t,x,y)\mapsto(G(t,x,y),H(t,x,y))$ is continuous and there is a
constant $C$ so that
 \beqnn
|G(t,x,y)|+|H(t,x,y)|\le C(1+y),\qquad t,y\ge0,~x\in\mbb{R}.
 \eeqnn
(2) Let $r_0$ be the concave function defined in condition (C2).
Then
 \beqnn
\mathop{\rm sgn}(y_1-y_2)(G(t,x,y_1)-G(t,x,y_2))\le
r_0(|y_1-y_2|),\quad t,y_1,y_2\ge0,~x\in\mbb{R}.
 \eeqnn
(3) ($\beta$-H\"older continuity) There is a constant $C$ so that
 \beqnn
|H(t,x,y_1)-H(t,x,y_2)|\le C|y_1-y_2|^\beta,\qquad
t,y_1,y_2\ge0,~x\in\mbb{R}.
 \eeqnn
(4) For fixed $t\ge0$ and $x\in\mbb{R}$, $H(t,x,y)$ is
nondecreasing in $y$.\\
Under the above conditions, by the same arguments in this paper, we
can show that the results of Theorems \ref{t1.1} and \ref{t4.1} also
hold for SPDE \eqref{1.05}. For simplicity we only study the SPDE
\eqref{1.01} in this paper.

The paper is organized as follows. In Section~2 we first present
some properties of the weak solution to equation \eqref{1.1}.
The proofs of Theorems \ref{t1.1} and \ref{t4.1} are established in
Sections 3 and 4, respectively. In
Section 5, the proofs of Proposition \ref{t1.5} and Lemma \ref{t1.4}
are presented.

{\bf Notation:} Throughout this paper, we adopt the conventions
 \beqnn
\int_x^y=\int_{(x,y]}\mbox{ and
}\int_x^\infty=\int_{(x,\infty)}
 \eeqnn
for any $y\ge x\ge 0$. Let $C$  denote a positive constant whose
value might change from line to line. We write $C_\varepsilon$ or
$C'_\varepsilon$  if the constant depends on another value
$\varepsilon\ge 0$. Let $\mbb{Q}$ be the notation for the set of
rational numbers. We sometimes write $\mbb{R}_+$ for $[0,\infty)$.
Let $(P_t)_{t\ge0}$ denote the transition semigroup of a one-dimensional
Brownian motion. For $t>0$ and $x\in\mbb{R}$ write $p_t(x):=(2\pi
t)^{-\frac12}\exp\{-x^2/(2t)\}$. We always use $N_0(ds,dz,du,dv)$ to
denote the Poisson random measure corresponding to the compensated
Poisson measure $\tilde{N}_0(ds,dz,du,dv)$.

\section{Properties of the weak solution}

\setcounter{equation}{0}

In this section we establish some properties of the weak solution to
\eqref{1.1}, which will be used in the next two sections. Recall the
measure $m_0(dz)=c_0z^{-1-\alpha}1_{\{z>0\}}dz$ for
$c_0={\alpha(\alpha-1)}/{\Gamma(2-\alpha)}$ where $\Gamma$ is the
Gamma function. By the proof of Theorem 1.1(a) of Mytnik and Perkins
(2003), there is a Poisson random measure $N(ds,dz,dx)$ on
$(0,\infty)^2\times \mbb{R}$ with intensity $dsm_0(dz)dx$ so that
 \beqlb\label{2.6}
L(ds,dx)=\int_0^\infty z\tilde{N}(ds,dz,dx),
 \eeqlb
where $\tilde{N}(ds,dz,du)$ is the compensated measure for
$N(ds,dz,dx)$. Thus, if $\{X_t:t\ge0\}$ is a weak solution of
\eqref{1.1}, then for each $f\in\mathscr{S}(\mbb{R})$ we have
 \beqlb\label{1.1b}
\<X_t,f\>
 \ar=\ar
X_0(f)+
\frac12\int_0^t\<X_s,f''\>ds+\int_0^tds\int_{\mbb{R}}G(X_s(x))f(x)dx
\cr
 \ar\ar
+ \int_0^t\int_0^\infty\int_{\mathbb{R}}H(X_{s-}(x))
f(x)z\tilde{N}(ds,dz,dx), \quad t>0,
 \eeqlb
which will be used to obtain \eqref{1.4}. For this we need
Assumption \ref{c1.1} on the weak solution of \eqref{1.1} for the
case $p>1$. For $0<p\le1$, by Definition \ref{t1.10} and the
H\"older inequality it is easy to check that the It\^o integrals in
\eqref{1.1} and \eqref{1.1b} are well defined. For $1<p<2$,
under Assumption \ref{c1.1} and by a similar argument it is easy to check
that the It\^o integrals in \eqref{1.1} and \eqref{1.1b} are also
well defined; see the details in Lemma \ref{r1.1}. By
\cite[Proposition 5.1]{M02} and the proof of \cite[Theorem
1.5]{M02}, for $G\equiv0$ and $H(x)=x^\beta$, the solution to
\eqref{1.1} exists and satisfies  Assumption \ref{c1.1}. In the
following proposition we always assume that conditions (C1) and (C3)
are satisfied and  Assumption \ref{c1.1} holds for the weak solution
$(X,L)$  to \eqref{1.1}.

\bproposition\label{t1.9} (i) If $(X,L)$ is a weak solution of
\eqref{1.1}, then there is, on an enlarged probability space, a
compensated Poisson random measure $\tilde{N}_0(ds,dz,du,dv)$ on
$(0,\infty)^2\times\mbb{R}\times(0,\infty)$ with intensity
$dsm_0(dz)dudv$ so that \eqref{1.4} holds. (ii) Conversely, if $X$
satisfies \eqref{1.4}, then there is an $\alpha$-stable white noise
$L(ds,dx)$ on $\mbb{R}_+\times \mbb{R}$ without negative jumps so
that \eqref{1.1} holds. \eproposition \proof (i) Suppose that
$(X,L)$ is a weak solution of \eqref{1.1}. Then by the argument at
the beginning of this section, \eqref{1.1b} holds. Define a
predictable $(0,\infty)\times(\mbb{R}\cup\{\infty\})$-valued process
$\theta(s,z,u,v)$ by
$\theta(s,z,u,v)=(\theta_1(s,z,u),\theta_2(s,u,v))$ with
 \beqnn
\theta_1(s,z,u):=\frac{z}{H(X_{s-}(u))} 1_{\{H(X_{s-}(u))\neq0\}}
+z1_{\{H(X_{s-}(u))=0\}}
 \eeqnn
and
 \beqnn
\theta_2(s,u,v):=\tilde{\theta}(s,u,v) 1_{\{H(X_{s-}(u))\neq0\}}
+\bar{\theta}(u,v)1_{\{H(X_{s-}(u))=0\}},
 \eeqnn
where \beqnn \tilde{\theta}(s,u,v):=\bigg\{
\begin{array}{l}
u,\quad v\le H(X_{s-}(u))^\alpha \\
\infty, \quad v> H(X_{s-}(u))^\alpha  \end{array}
,\qquad \bar{\theta}(u,v):=\bigg\{
\begin{array}{l}
u, \quad v\in(0,1)
\\
\infty,\quad v\in(0,1)^c
 \end{array}
 \eeqnn
and we use the convention that $0\cdot\infty=0$. Then for all
$B\in\mathscr{B}(0,\infty)$ and $a\le b\in\mbb{R}$,
 \beqnn
 \ar\ar
1_{B\times(a,b]}(\theta(s,z,u,v))
=1_{B\times(a,b]}(\theta_1(s,z,u),\theta_2(s,u,v)) \cr
 \ar=\ar
1_{\{H(X_{s-}(u))\neq0,u\in(a,b],v\le H(X_{s-}(u))^\alpha\}}1_B\Big(\frac{z}{H(X_{s-}(u))}\Big)
+ 1_{\{H(X_{s-}(u))=0,z\in B,u\in(a,b],v\in (0,1)\}}.
 \eeqnn
Moreover, recalling $m_0(dz)=c_0z^{-1-\alpha}1_{\{z>0\}}dz$, by a
change of variable it is easy to see that
 \beqnn
 \ar\ar
\int_0^\infty\int_{\mbb{R}}\int_0^\infty 1_{B\times
(a,b]}(\theta(s,z,u,v)) m_0(dz)du dv  \cr
 \ar\ar\quad=
\int_0^\infty m_0(dz)\int_a^b du \int_0^{H(X_{s-}(u))^\alpha}
1_{\{H(X_{s-}(u))\neq0\}} 1_B\Big(\frac{z}{H(X_{s-}(u))}\Big)dv
\cr
 \ar\ar\qquad
+ \int_0^\infty\int_a^b 1_{\{H(X_{s-}(u))=0\}}1_B(z) m_0(dz) du \cr
 \ar\ar\quad=
\int_0^\infty\int_{\mbb{R}}1_{B\times(a,b]}(z,u)m_0(dz)du.
 \eeqnn
Then by \cite[p.93]{IkW89}, on an extension of the probability
space, there exists a Poisson random measure $N_0(ds,dz,du,dv)$ on
$(0,\infty)^2\times\mbb{R}\times(0,\infty)$ with intensity
$dsm_0(dz)dudv$ so that
 \beqnn
N((0,t]\times B\times (a,b])
= \int_0^t\int_0^\infty\int_{\mbb{R}}\int_0^\infty
1_{B\times(a,b]}(\theta(s,z,u,v))N_0(ds,dz,du,dv).
 \eeqnn
Let $\tilde{N}_0(ds,dz,du,dv)=N_0(ds,dz,du,dv)-dsm_0(dz)dudv$. Then
by \eqref{2.6} it is easy to see that for each $f\in
\mathscr{S}(\mbb{R})$,
 \beqnn
 \ar\ar
\int_0^t\int_{\mathbb{R}}H(X_{s-}(u))f(u)L(ds,du)=
\int_0^t\int_0^\infty\int_{\mbb{R}}\int_0^{H(X_{s-}(u))^\alpha}
zf(u)\tilde{N}_0(ds,dz,du,dv).
 \eeqnn

(ii) The proof is essentially the same as that of \cite[Theorem 9.32]{Li11}.
Suppose that $\{X_t:t>0,x\in\mbb{R}\}$ satisfies \eqref{1.4}.
Define the random measure $N(ds,dz,du)$ on $(0,\infty)^3$ by
 \beqnn
 \ar\ar
N((0,t]\times B\times (a,b]) \cr
 \ar\ar\quad:=
\int_0^t\int_0^\infty\int_a^b\int_0^{H(X_{s-}(u))^\alpha}
1_{\{H(X_{s-}(u))\neq0\}} 1_B\Big( \frac{z}{H(X_{s-}(u))} \Big)N_0(ds,dz,du,dv) \cr
 \ar\ar\qquad \,
+ \int_0^t\int_0^\infty\int_a^b\int_0^1 1_{\{H(X_{s-}(u))=0\}}1_B(z) N_0(ds,dz,du,dv).
 \eeqnn
It is easy to see that $N(ds,dz,du)$ has a predictable compensator
 \beqnn
\hat{N}((0,t]\times B\times (a,b])
 \ar=\ar
\int_0^t\int_0^\infty\int_a^b\int_0^{H(X_{s-}(u))^\alpha}
1_{\{H(X_{s-}(u))\neq0\}} 1_B\Big( \frac{z}{H(X_{s-}(u))}
\Big)dsm_0(dz)dudv \cr
 \ar\ar\quad
+ \int_0^t\int_0^\infty\int_a^b\int_0^1 1_{\{H(X_{s-}(u))=0\}}
1_B(z) dsm_0(dz)dudv \cr
 \ar=\ar
\int_0^t\int_0^\infty \int_a^b 1_B(z)dsm_0(dz)du.
 \eeqnn
Then $N(ds,dz,du)$ is a Poisson random measure with intensity
$dsm_0(dz)du$; see Theorems II.1.8 and II.4.8 in \cite{JSh87}.
Define the $\alpha$-stable white noise $L$ by
 \beqnn
L_t(a,b]=\int_0^t\int_0^\infty\int_a^bz\tilde{N}(ds,dz,du).
 \eeqnn
We then have
 \beqnn
 \ar\ar
\int_0^t\int_{\mathbb{R}}H(X_{s-}(u))f(u)L(ds,du)
=
\int_0^t\int_0^\infty\int_{\mathbb{R}}H(X_{s-}(u)) f(u)z\tilde{N}(ds,dz,du) \cr
 \ar=\ar
\int_0^t\int_0^\infty\int_{\mbb{R}}\int_0^{H(X_{s-}(u))^\alpha }
zf(u) \tilde{N}_0(ds,dz,du,dv)
 \eeqnn
for each $f\in \mathscr{S}(\mbb{R})$.  $(X,L)$ is thus a weak
solution to \eqref{1.1}. \qed

In the rest of this section, we always assume that conditions (C1)
and (C3) are satisfied and $(X,L)$ is a weak solution to \eqref{1.1}
with deterministic initial value $X_0\in M(\mbb{R})$ and Assumption
\ref{c1.1} satisfied. Then it follows from Proposition \ref{t1.9},
$\{X_t(x):t>0,x\in\mbb{R}\}$ satisfies \eqref{1.4}. Recall that
$(P_t)_{t\ge0}$ is the transition semigroup of a one-dimensional
Brownian motion and $p_t(x)=(2\pi t)^{-\frac12}\exp\{-x^2/(2t)\}$
for $t>0$ and $x\in\mbb{R}$.

\bproposition\label{t1.5} For any $t>0$ and $f\in B(\mbb{R})$
satisfying $\lambda_0(|f|)<\infty$ we have
 \beqlb\label{1.2}
\<X_t,f\> \ar=\ar X_0(P_tf) +
\int_0^tds\int_{\mbb{R}}G(X_s(x))P_{t-s}f(x)dx \cr
 \ar\ar
+\int_0^t\int_0^\infty\int_{\mbb{R}}\int_0^{H(X_{s-}(u))^\alpha}
z P_{t-s}f(u)\tilde{N}_0(ds,dz,du,dv),~~ \mbf{P}\mbox{-a.s}.
 \eeqlb
 Moreover,
 \beqlb\label{1.3}
X_t(x)
 \ar=\ar
\int_{\mbb{R}}p_t(x-z)X_0(dz)+
\int_0^tds\int_{\mbb{R}}p_{t-s}(x-z)G(X_s(z))dz \cr
 \ar\ar
+ \int_0^t\int_0^\infty\int_{\mbb{R}}\int_0^{H(X_{s-}(u))^\alpha }
z p_{t-s}(x-u)\tilde{N}_0(ds,dz,du,dv),
\mbf{P}\mbox{-a.s., }\lambda_0\mbox{-a.e. }x.
 \eeqlb
\eproposition The proof is given in the Appendix.

We refer to \cite[Theorem 8.6]{PZ07} for
the stochastic integration with respect to a Poisson random measure.
For $k>0$ let $\tilde{\tau}_k$ be a stopping time  defined by
 \beqlb\label{1.21}
\tilde{\tau}_k:=\inf\{t: F(t)>k\}
 \eeqlb
with the convention $\inf \emptyset=\infty$, where
$F(t):=(\int_0^tds\int_{\mbb{R}}X_s(x)^qdx)\vee \<X_t,1\>$ for the
case $p>1$ and $F(t):=\<X_t,1\>$ for the case $p\le1$. Then it follows from
Definition \ref{t1.10} and Assumption \ref{t1.9} that
 \beqlb\label{1.21b}
\lim_{k\to\infty}\tilde{\tau}_k=\infty,\qquad \mbf{P}\mbox{-a.s}.
 \eeqlb
The following lemma says that Assumption \ref{c1.1}
also assures that the It\^o integrals in \eqref{1.2} and \eqref{1.3}
are well defined.
 \blemma\label{r1.1}
If Assumption \ref{c1.1} holds,
the It\^o integrals in \eqref{1.2} and \eqref{1.3}
are well defined.
 \elemma
 \proof
Since the reasoning is similar, we only explain
that of \eqref{1.3} in the following. We first
consider the case $0<p\le1$. Since $u^p\le u+1$ for $u\ge0$, then for each
$1\le\bar{\alpha}<2$ and any $x\in \mbb{R}$,
 \beqlb\label{1.32}
 \ar\ar
\mbf{E}\Big\{\int_0^{t\wedge\tilde{\tau}_k}ds\int_{\mbb{R}}X_s(u)^p
p_{t-s}(x-u)^{\bar{\alpha}} du\Big\} \cr
 \ar\ar\le
\mbf{E}\Big\{\int_0^t[2\pi(t-s)]^{-\frac{\bar{\alpha}-1}{2}}ds
\int_{\mbb{R}}[X_s(u)+1] p_{t-s}(x-u)1_{\{s\le \tilde{\tau}_k\}}
du\Big\} \cr
 \ar\ar\le
\int_0^t[2\pi(t-s)]^{-\frac{\bar{\alpha}-1}{2}}
\mbf{E}\Big\{1+[2\pi(t-s)]^{-1/2}\<X_s,
1\>1_{\{s\le \tilde{\tau}_k\}}\Big\}ds  \cr
 \ar\ar\le
\int_0^t[2\pi(t-s)]^{-\frac{\bar{\alpha}-1}{2}}
\Big\{1+[2\pi(t-s)]^{-1/2}k\Big\}ds <\infty.
 \eeqlb
Therefore, by condition (C3), for any $k\ge1$,
 \beqnn
 \ar\ar
\mbf{E}\Big\{\int_0^{t\wedge\tilde{\tau}_k}ds\int_1^\infty
zm_0(dz)\int_{\mbb{R}}du\int_0^{H(X_{s-}(u))^\alpha} | p_{t-s}(x-u)|dv\Big\} \cr
 \ar\ar\qquad\le
C\,\mbf{E}\Big\{\int_0^{t\wedge\tilde{\tau}_k}ds\int_1^\infty
zm_0(dz)\int_{\mbb{R}} [1+X_s(u)^p]p_{t-s}(x-u)du\Big\}<\infty
 \eeqnn
and
 \beqnn
 \ar\ar
\mbf{E}\Big\{\int_0^{t\wedge\tilde{\tau}_k}ds\int_0^1
z^{\bar{\bar{\alpha}}}m_0(dz)
\int_{\mbb{R}}du\int_0^{H(X_{s-}(u))^\alpha} |  p_{t-s}(x-u)|^{\bar{\bar{\alpha}}}dv\Big\} \cr
 \ar\ar\qquad\le
C\,\mbf{E}\Big\{\int_0^{t\wedge\tilde{\tau}_k}ds\int_0^1
z^{\bar{\bar{\alpha}}}m_0(dz) \int_{\mbb{R}}
[1+X_s(u)^p]p_{t-s}(x-u)^{\bar{\bar{\alpha}}}du\Big\}<\infty
 \eeqnn
for $\alpha<\bar{\bar{\alpha}}<2$,
which ensures that the stochastic integral in
\eqref{1.3} is well defined.

In the following we consider the case $p>1$.
Observe that $\frac{q}{p}>\frac{3}{3-\alpha}>\frac32$,
 which implies $\frac12(\frac{q}{q-p}-1)<1$.
Thus by the H\"older inequality,
 \beqlb\label{1.17}
 \ar\ar
\mbf{E}\Big\{\int_0^{t\wedge\tilde{\tau}_k}ds\int_{\mbb{R}}X_s(u)^p
p_{t-s}(x-u)du\Big\} \cr \ar\ar\quad\le
\Big\{\mbf{E}\Big[\int_0^{t\wedge\tilde{\tau}_k}ds
\int_{\mbb{R}}X_s(u)^qdu\Big]\Big\}^{\frac{p}{q}}
\Big\{\int_0^tds\int_{\mbb{R}}
p_{t-s}(x-u)^{\frac{q}{q-p}}du\Big\}^{1-\frac pq}  \cr
 \ar\ar\quad\le
k^{\frac pq}\Big\{\int_0^t[2\pi(t-s)]^{-\frac12(\frac{q}{q-p}-1)}ds\int_{\mbb{R}}
p_{t-s}(x-u)du\Big\}^{1-\frac pq} \cr
 \ar\ar\quad=
k^{\frac pq}\Big\{\int_0^t[2\pi(t-s)]^{-\frac12(\frac{q}{q-p}-1)}ds\Big\}^{1-\frac pq}
<\infty.
 \eeqlb
Observe that $q>\frac{3p}{3-\alpha}$ implies
$q\alpha/(q-p)<3$.
Then similar to \eqref{1.17}, there is a constant $\alpha<\bar{\alpha}'<2$ so that
$q\bar{\alpha}'/(q-p)<3$ and
 \beqlb\label{1.10}
\mbf{E}\Big\{\int_0^{t\wedge\tilde{\tau}_k}ds\int_{\mbb{R}}X_s(u)^p
p_{t-s}(x-u)^{\bar{\alpha}'}du\Big\}
 \le
k^{\frac pq}
\Big\{\int_0^t[2\pi(t-s)]^{-\frac12(\frac{q\bar{\alpha}'}{q-p}-1)}ds\Big\}^{1-\frac pq}
<\infty.
 \eeqlb
Therefore, by condition (C3) again, for any $k\ge1$,
 \beqnn
 \ar\ar
\mbf{E}\Big\{\int_0^{t\wedge\tilde{\tau}_k}ds\int_1^\infty
zm_0(dz)\int_{\mbb{R}}du\int_0^{H(X_{s-}(u))^\alpha} | p_{t-s}(x-u)|dv\Big\} \cr
 \ar\ar\qquad\le
C\,\mbf{E}\Big\{\int_0^{t\wedge\tilde{\tau}_k}ds \int_1^\infty
zm_0(dz)\int_{\mbb{R}} [1+X_s(u)^p]p_{t-s}(x-u)du\Big\}<\infty
 \eeqnn
and
 \beqnn
 \ar\ar
\mbf{E}\Big\{\int_0^{t\wedge\tilde{\tau}_k}ds\int_0^1
z^{\bar{\alpha}'}m_0(dz)
\int_{\mbb{R}}du\int_0^{H(X_{s-}(u))^\alpha} | p_{t-s}(x-u)|^{\bar{\alpha}'}dv\Big\} \cr
 \ar\ar\qquad\le
C\,\mbf{E}\Big\{\int_0^{t\wedge\tilde{\tau}_k}ds\int_0^1
z^{\bar{\alpha}'}m_0(dz)\int_{\mbb{R}}
[1+X_s(u)^p]p_{t-s}(x-u)^{\bar{\alpha}'}du\Big\}<\infty.
 \eeqnn
So, the stochastic integral in \eqref{1.3} is also well defined.
 \qed

\blemma\label{t1.4} Let $0< \bar{p}< \alpha$   be fixed.
Then for  any $T>0$ and any $0<t\le T$, there is a set $K_t\subset\mbb{R}$ of
Lebesgue measure zero so that
 \beqlb\label{1.22}
\mbf{E}\{X_t(x)^{ \bar{p}}\}\le C_Tt^{-\frac{\bar{p}}{2}},
\qquad x\in\mbb{R} \backslash K_t.
 \eeqlb
\elemma The proof is also given in the Appendix.

\blemma\label{t1.6} Suppose that $T>0$, $\delta\in(1,\alpha)$,
$\delta_1\in(\alpha,2)$ and $0<
r<\min\{1,\frac{3-\delta_1}{\delta_1}\}$. Then for each $0<t\le T$
and the set $K_t\subset\mbb{R}$ from Lemma \ref{t1.4} we have
 \beqlb\label{1.5}
\mbf{E}\{|X_t(x_1)-X_t(x_2)|^\delta\}\le
C_Tt^{-\frac{(r+1)\delta}{2}}|x_1-x_2|^{r\delta},\qquad x_1,x_2\in\mbb{R}\backslash K_t.
 \eeqlb
\elemma
\proof
For $t>0$ and $x\in\mbb{R}$ let
 \beqnn
Z_1(t,x):=\int_0^t\int_0^1\int_{\mbb{R}}\int_0^{H(X_{s-}(u))^\alpha }
z p_{t-s}(x-u)\tilde{N}_0(ds,dz,du,dv)
 \eeqnn
and
 \beqnn
Z_2(t,x):=\int_0^t\int_1^\infty\int_{\mbb{R}}\int_0^{H(X_{s-}(u))^\alpha }
z p_{t-s}(x-u)\tilde{N}_0(ds,dz,du,dv).
 \eeqnn
By (2.4e) in \cite{Rosen}, for all $t>0$, $\theta\in[0,1]$ and
$u\in\mbb{R}$ we have
 \beqlb\label{1.35}
|p_t(x_1-u)-p_t(x_2-u)|\le C|x_1-x_2|^\theta t^{-\theta/2}[p_t(x_1-u)+p_t(x_2-u)].
 \eeqlb
Then by (1.6) in \cite{saint}, condition (C3) and Lemma \ref{t1.4},
 \beqlb\label{1.14}
 \ar\ar
\mbf{E}\{|Z_1(t,x_1)-Z_1(t,x_2)|^{\delta_1}\} \cr
 \ar\ar\qquad\le
C\int_0^1z^{\delta_1}m_0(dz)\int_0^tds\int_{\mbb{R}}
\mbf{E}\{H(X_{s-}(u))^\alpha\}|p_{t-s}(x_1-u)-p_{t-s}(x_2-u)|^{\delta_1}du
\cr
 \ar\ar\qquad\le
C\int_0^1z^{\delta_1}m_0(dz)\int_0^tds\int_{\mbb{R}}
\mbf{E}\{1+X_s(u)^p\}|p_{t-s}(x_1-u)-p_{t-s}(x_2-u)|^{\delta_1}du
\cr
 \ar\ar\qquad\le
C_T|x_1-x_2|^{r\delta_1}\int_0^t[1+s^{-p/2}](t-s)^{-\frac{r\delta_1+\delta_1-1}{2}}ds
\int_{\mbb{R}}[p_{t-s}(x_1-u)+p_{t-s}(x_2-u)]du \cr
 \ar\ar\qquad\le
C_T|x_1-x_2|^{r\delta_1}\int_0^t[1+s^{-p/2}](t-s)^{-\frac{r\delta_1+\delta_1-1}{2}}ds.
 \eeqlb
It follows from the  H\"older inequality that
 \beqlb\label{1.11}
\mbf{E}\{|Z_1(t,x_1)-Z_1(t,x_2)|^\delta\}
\le
\Big\{\mbf{E}\big[|Z_1(t,x_1)-Z_1(t,x_2)|^{\delta_1}\big]\Big\}^{\frac{\delta}{\delta_1}}.
 \eeqlb

Similar to \eqref{1.14} we have
 \beqlb\label{1.12}
\mbf{E}\{|Z_2(t,x_1)-Z_2(t,x_2)|^\delta\}
\le
C_T|x_1-x_2|^{r\delta}\int_0^t[1+s^{-p/2}](t-s)^{-\frac{r\delta+\delta-1}{2}}ds.
 \eeqlb
Combining \eqref{1.14}--\eqref{1.12} one has
 \beqlb\label{1.13}
\mbf{E}\Big\{|Z_1(t,x_1)-Z_1(t,x_2)|^\delta
+|Z_2(t,x_1)-Z_2(t,x_2)|^\delta\Big\}
\le C_Tt^{-\frac{(r+1)\delta}{2}}|x_1-x_2|^{r\delta}.
 \eeqlb
By the H\"older inequality and condition (C1), \beqlb\label{1.24}
 \ar\ar
\mbf{E}\Big\{\Big|\int_0^tds\int_{\mbb{R}}p_{t-s}(x_1-z)G(X_s(z))dz
-\int_0^tds\int_{\mbb{R}}p_{t-s}(x_2-z)G(X_s(z))dz\Big|^\delta\Big\} \cr
 \ar\ar\qquad\le
2\,\mbf{E}\Big\{\int_0^tds\int_{\mbb{R}}
|p_{t-s}(x_1-z)-p_{t-s}(x_2-z)|G(X_s(z))^\delta dz\Big\}  \cr
 \ar\ar\qquad\le
C\,\mbf{E}\Big\{\int_0^tds\int_{\mbb{R}}
|p_{t-s}(x_1-z)-p_{t-s}(x_2-z)|[1+X_s(z)^\delta]dz\Big\} \cr
 \ar\ar\qquad\le
C|x_1-x_2|^{r\delta}\int_0^t(t-s)^{-r\delta/2}[1+s^{-\delta/2}]ds
\int_{\mbb{R}}[p_{t-s}(x_1-u)+p_{t-s}(x_2-u)]du \cr
 \ar\ar\qquad\le
C_T|x_1-x_2|t^{-\frac{(r+1)\delta}{2}},
 \eeqlb
where Lemma \ref{t1.4} and \eqref{1.35} was used in the third
inequality. By \eqref{1.35} again we have
 \beqnn
 \ar\ar
\Big|\int_{\mbb{R}}p_t(x_1-y)X_0(dy)-\int_{\mbb{R}}p_t(x_2-y)X_0(dy)\Big| \cr
 \ar\ar\qquad\le
\int_{\mbb{R}}|p_t(x_1-y)-p_t(x_2-y)|X_0(dy) \le
C|x_1-x_2|^rt^{-\frac{r+1}{2}}X_0(1),
 \eeqnn
which together with \eqref{1.3} and \eqref{1.13}--\eqref{1.24} implies
\eqref{1.5}.
\qed

\blemma\label{t1.8}
For each $t>0$ and $t_n>0$ satisfying $t_n\to t$ as $n\to\infty$,
there is a set $K_t\subset\mbb{R}$ of Lebesgue measure zero so that
 \beqnn
\lim_{n\to\infty}\mbf{E}\{|X_{t_n}(x)-X_t(x)|\}=0,\qquad
x\in\mbb{R} \backslash K_t.
 \eeqnn
\elemma
\proof
For $t_0,t>0$, by \eqref{1.3},
 \beqnn
 \ar\ar
|X_{t_0+t}(x)-X_{t_0}(x)|  \cr
 \ar\ar\quad\le
\int_{\mbb{R}}|p_{t_0+t}(x-y)-p_{t_0}(x-y)|X_0(dy)
+\int_{t_0}^{t_0+t}ds\int_{\mbb{R}}p_{t_0+t-s}(x-z)G(X_s(z))dz \cr
 \ar\ar\quad\quad
+\int_0^{t_0}ds\int_{\mbb{R}}\big|[p_{t_0+t-s}(x-z)-p_{t_0-s}(x-z)]G(X_s(z))\big|dz  \cr
 \ar\ar\quad\quad
+\Big|\int_{t_0}^{t_0+t}\int_0^\infty \int_{\mbb{R}}
\int_0^{H(X_{s-}(u))^\alpha }
z
p_{t_0+t-s}(x-u)\tilde{N}_0(ds,dz,du,dv)\Big| \cr
 \ar\ar\quad\quad
+\Big|\int_0^{t_0}\int_0^\infty \int_{\mbb{R}}
\int_0^{H(X_{s-}(u))^\alpha }
z [p_{t_0+t-s}(x-u)-p_{t_0-s}(x-u)]\tilde{N}_0(ds,dz,du,dv)\Big|
\cr
 \ar\ar\quad=:
I_1(t_0,t)+I_2(t_0,t)+I_3(t_0,t) +|I_4(t_0,t)|+|I_5(t_0,t)|.
 \eeqnn
By the dominated convergence, $I_1(t_0,t)$ tends to zero as $t\to0$.

By condition (C1), Lemma \ref{t1.4} and the dominated convergence,
we have that both
 \beqnn
\mbf{E}\{I_2(t_0,t)\}
\le
C\,\mbf{E}\Big\{\int_{t_0}^{t_0+t}ds\int_{\mbb{R}}
p_{t_0+t-s}(x-z)[1+X_s(z)]dz\Big\}
\le
C\int_{t_0}^{t_0+t}(1+s^{-\frac12})ds
 \eeqnn
and
 \beqnn
\mbf{E}\{I_3(t_0,t)\}
 \ar\le\ar
C\,\mbf{E}\Big\{\int_0^{t_0}ds\int_{\mbb{R}}
|p_{t_0+t-s}(x-z)-p_{t_0-s}(x-z)|[1+X_s(z)]dz\Big\} \cr
 \ar\le\ar
C_T\int_0^{t_0}(1+s^{-\frac12})ds\int_{\mbb{R}}|p_{t_0+t-s}(x-z)-p_{t_0-s}(x-z)|dz
 \eeqnn
go to $0$ as $t\to0$.

Let
 \beqnn
I_{4,1}(t_0,t):=\int_{t_0}^{t_0+t}\int_0^1
\int_{\mbb{R}}\int_0^{H(X_{s-}(u))^\alpha }
z
p_{t_0+t-s}(x-u)\tilde{N}_0(ds,dz,du,dv)
 \eeqnn
and
 \beqnn
I_{5,1}(t_0,t)
 :=
\int_0^{t_0}\int_0^1
\int_{\mbb{R}}\int_0^{H(X_{s-}(u))^\alpha }
z
[p_{t_0+t-s}(x-u)
-p_{t_0-s}(x-u)]\tilde{N}_0(ds,dz,du,dv).
 \eeqnn
Let $I_{4,2}(t_0,t):=I_4(t_0,t)-I_{4,1}(t_0,t)$ and
$I_{5,2}(t_0,t):=I_5(t_0,t)-I_{5,1}(t_0,t)$. Then by the H\"older
continuity of $H$, Lemma \ref{t1.4} and the dominated convergence
again, both
 \beqnn
\mbf{E}\{I_{4,1}(t_0,t)^2\}
 \ar=\ar
\int_0^1 z^2m_0(dz)\int_{t_0}^{t_0+t}ds\int_{\mbb{R}}
\mbf{E}\{H(X_s(u))^\alpha\}p_{t_0+t-s}(x-u)^2du \cr
 \ar\le\ar
C\int_0^1 z^2m_0(dz)\int_{t_0}^{t_0+t}ds\int_{\mbb{R}}
\mbf{E}\{1+X_s(u)^p\}p_{t_0+t-s}(x-u)^2du \cr
 \ar\le\ar
C_T\int_{t_0}^{t_0+t}[1+s^{-\frac{p}{2}}](t_0+t-s)^{-\frac12}ds
 \eeqnn
and
 \beqnn
\mbf{E}\{I_{5,1}(t_0,t)^2\}
 \ar=\ar
\int_0^1 z^2m_0(dz)\int_0^{t_0}ds\int_{\mbb{R}}
\mbf{E}\{H(X_s(u))^\alpha\}[p_{t_0+t-s}(x-u)-p_{t_0-s}(x-u)]^2du
\cr
 \ar\le\ar
C_T\int_0^{t_0}[1+s^{-\frac{p}{2}}]ds\int_{\mbb{R}}[p_{t_0+t-s}(x-u)-p_{t_0-s}(x-u)]^2du
 \eeqnn
go to $0$ as $t\to0$.

 Similarly, both
 \beqnn
 \ar\ar
\mbf{E}\{|I_{4,2}(t_0,t)|\} \cr
 \ar\ar\quad\le
\mbf{E}\Big\{\Big|\int_{t_0}^{t_0+t}\int_1^\infty \int_{\mbb{R}}
\int_0^{H(X_{s-}(u))^\alpha }
z
p_{t_0+t-s}(x-u)N_0(ds,dz,du,dv)\Big|\Big\} \cr
 \ar\ar\qquad
+\mbf{E}\Big\{\int_1^\infty
zm_0(dz)\int_{t_0}^{t_0+t}ds\int_{\mbb{R}} H(X_{s-}(u))^\alpha p_{t_0+t-s}(x-u)du\Big\} \cr
 \ar\ar\quad\le
2\int_1^\infty zm_0(dz)\int_{t_0}^{t_0+t}ds\int_{\mbb{R}}
\mbf{E}\{H(X_s(u))^\alpha\}p_{t_0+t-s}(x-u)du
 \eeqnn
and
 \beqnn
\mbf{E}\{|I_{5,2}(t_0,t)|\} \le 2\int_1^\infty
zm_0(dz)\int_0^{t_0}ds\int_{\mbb{R}}
\mbf{E}\{H(X_s(u))^\alpha\}|p_{t_0+t-s}(x-u)-p_{t_0-s}(x-u)|du
 \eeqnn
go to $0$ as $t\to0$. The proof is thus completed. \qed

\section{Proof of Theorem \ref{t1.1}}

\setcounter{equation}{0}

In this section we establish the proof of Theorem \ref{t1.1}. Throughout this
section we always assume that conditions (C1) and (C3) hold and that
$(X,L)$ is a weak solution to \eqref{1.1} satisfying Assumption
\ref{c1.1} with deterministic initial value $X_0\in M(\mbb{R})$. For
$n\ge1$ and $0\le k\le 2^n$, put $n_k:=k/2^n$. Define ${n'}_k$ similarly
for $n'\ge1$.

For any $x\in\mbb{R}$ and $s>0$ let
 \beqlb\label{2.18}
Z_s(x)=\int_0^s\int_0^\infty\int_{\mbb{R}}
\int_0^{H(X_{s_1-}(u))^\alpha} zp_{s-s_1}(x-u)\tilde{N}_0(ds_1,dz,du,dv).
 \eeqlb

Before presenting the proof of Theorem \ref{t1.1}, we first
establish a weaker version of the result which will be used to
conclude the proof of Theorem \ref{t1.1}.

 \blemma\label{t1.7} The results of Theorem \ref{t1.1} hold
with $\eta_c$ replaced by $\eta_c'=\eta_c1_{\{\alpha\ge\frac32\}}
+\frac{\alpha-1}{\alpha}1_{\{\alpha<\frac32\}}$. \elemma \proof Let
$r,\delta$ and $\delta_1$ satisfy the conditions in Lemma \ref{t1.6}
and $r\delta>1$. By \eqref{1.13} and the proof of Corollary 1.2(ii)
of \cite{Walsh}, for each $0<\varepsilon<r-1/\delta$ and $T>0$,
there is a constant $C_T$ independent of $t\in(0,T]$ so that
 \beqnn
\mbf{E}\Big\{ \sup_{x,y\in \mathbb{K},x\neq y}
\frac{|\tilde{Z}_t(x)-\tilde{Z}_t(y)|}{|x-y|^\varepsilon}\Big\} \le
C_Tt^{-\frac{r+1}{2}},
 \eeqnn
where $\tilde{Z}_t(x)$ denotes a continuous modification of $Z_t(x)$
for each $t>0$. Then by Fatou's lemma, for each subsequence
$\{n':n'\ge1\}$ of $\{n:n\ge1\}$,
 \beqnn
 \ar\ar
\mbf{E}\Big\{\liminf_{n\to\infty}
\frac{1}{2^n}\sum_{k=1}^{2^n}\sup_{x,y\in \mathbb{K},x\neq y}
\frac{|\tilde{Z}_{n_kT}(x)-\tilde{Z}_{n_kT}(y)|}{|x-y|^\varepsilon}\Big\}
\cr
 \ar\ar\qquad\le
\liminf_{n\to\infty}\mbf{E}\Big\{
\frac{1}{2^n}\sum_{k=1}^{2^n}\sup_{x,y\in \mathbb{K},x\neq y}
\frac{|\tilde{Z}_{n_kT}(x)-\tilde{Z}_{n_kT}(y)|}{|x-y|^\varepsilon}\Big\}
\cr
 \ar\ar\qquad =
\liminf_{n\to\infty}\frac{1}{2^n}\sum_{k=1}^{2^n}\mbf{E}\Big\{
\sup_{x,y\in \mathbb{K},x\neq y}
\frac{|\tilde{Z}_{n_kT}(x)-\tilde{Z}_{n_kT}(y)|}{|x-y|^\varepsilon}\Big\}
\cr
 \ar\ar\qquad\le
C_T\liminf_{n\to\infty}\frac{1}{2^n}\sum_{k=1}^{2^n}(n_kT)^{-\frac{r+1}{2}}
=C_T\int_0^Ts^{-\frac{r+1}{2}}ds <\infty,
 \eeqnn
which implies
 \beqnn
\liminf_{n'\to\infty} \frac{1}{2^{n'}}\sum_{k=1}^{2^{n'}}
\sup_{x,z\in \mathbb{K},x\neq z}
\frac{|\tilde{Z}_{{n'}_kT}(x)-\tilde{Z}_{{n'}_kT}(z)|}
{|x-z|^\varepsilon} <\infty, \,\,\,\,\, \mbf{P}\mbox{-a.s.}
 \eeqnn
Let
 \beqnn
\tilde{X}_t(x)=\int_{\mbb{R}}p_t(x-u)X_0(du)
+\int_0^tds\int_{\mbb{R}}p_{t-s}(x-z) G(X_s(z))dz +\tilde{Z}_t(x).
 \eeqnn
Then it follows from \eqref{1.3} that $\tilde{X}_t(x)$ is a
continuous modification of $X_t(x)$ for each fixed $t>0$. We first
show that \eqref{2.17} holds for $\eta=\varepsilon\in(0,r-1/\delta)$
in the following. By \eqref{1.35} and condition (C1) for each
$\varepsilon'\in(0,1)$,
 \beqnn
\int_{\mbb{R}}|p_t(x-u)-p_t(y-u)|X_0(du) \le C|x-y|^{\varepsilon'}
t^{-\frac{{\varepsilon'}+1}{2}}X_0(1)
 \eeqnn
and
 \beqnn
 \ar\ar
\int_0^tds\int_{\mbb{R}}|p_{t-s}(x-z)-p_{t-s}(y-z)|G(X_s(z))dz \cr
 \ar\ar\qquad\le
C\int_0^tds\int_{\mbb{R}}|p_{t-s}(x-z)-p_{t-s}(y-z)|[X_s(z)+1]dz \cr
 \ar\ar\qquad\le
C_T|x-y|^{\varepsilon'}\int_0^t
(t-s)^{-\frac{\varepsilon'+1}{2}}[\<X_s,1\>+1]ds
 \eeqnn
for each $t>0$.
Then
 \beqlb\label{1.25}
 \ar\ar
\limsup_{n'\to\infty} \frac{1}{2^{n'}}\sum_{k=1}^{2^{n'}}
\sup_{x,y\in \mathbb{K},x\neq y}
\frac{\int_{\mbb{R}}|p_{{n'}_kT}(x-u)-p_{{n'}_kT}(y-u)|X_0(du)}{|x-y|^{\varepsilon'}}
\cr
 \ar\ar\qquad\le
\limsup_{n'\to\infty} \frac{C}{2^{n'}}\sum_{k=1}^{2^{n'}}
({n'}_kT)^{-\frac{\varepsilon'+1}{2}}X_0(1)
=CX_0(1)\int_0^Tt^{-\frac{\varepsilon'+1}{2}}ds<\infty
 \eeqlb
and
 \beqlb\label{1.26}
 \ar\ar
\limsup_{n'\to\infty} \frac{1}{2^{n'}}\sum_{k=1}^{2^{n'}}
\sup_{x,y\in \mathbb{K},x\neq y} \frac{1}{|x-y|^{\varepsilon'}}
\int_0^{{n'}_kT}ds\int_{\mbb{R}}|p_{{{n'}_kT}-s}(x-z) \cr
 \ar\ar\qqquad\qqquad\qqquad\quad
-p_{{{n'}_kT}-s}(y-z)|G(X_s(z))dz \cr
 \ar\ar\quad\le
\limsup_{n'\to\infty} \frac{C}{2^{n'}}\sum_{k=1}^{2^{n'}}
\int_0^{{n'}_kT}
({{n'}_kT}-s)^{-\frac{\varepsilon'+1}{2}}[\<X_s,1\>+1]ds \cr
 \ar\ar\quad=
C_T\int_0^Tdt\int_0^t
(t-s)^{-\frac{\varepsilon'+1}{2}}[\<X_s,1\>+1]ds \cr
 \ar\ar\quad=
C_T\big[\sup_{s\in(0,T]}\<X_s,1\>+1\big]\int_0^Tdt\int_0^t
(t-s)^{-\frac{\varepsilon'+1}{2}}ds<\infty,\quad\mbf{P}\mbox{-a.s.},
 \eeqlb
where the fact $\sup_{s\in(0,T]}\<X_s,1\><\infty$ $\mbf{P}$-a.s. was
used in the last inequality. Therefore, \eqref{2.17} holds for
$\eta<r-1/\delta$. Let $\delta=\alpha-\sigma$ and
$\delta_1=\alpha+\sigma$ for small enough $\sigma>0$. This means
that \eqref{2.17} holds for
 \beqnn
\eta<\min\Big\{1,\frac{3}{\alpha+\sigma}-1\Big\}-\frac{1}{\alpha-\sigma}
=\min\Big\{1-\frac{1}{\alpha-\sigma},
\frac{3}{\alpha+\sigma}-\frac{1}{\alpha-\sigma}-1\Big\}.
 \eeqnn
Letting $\sigma\to0$ one can finish the proof. \qed

\blemma\label{t2.1} For any fixed $t>0$, let $\tilde{X}_t$ be a
continuous modification of $X_t$. Then for any compact subset
$\mathbb{K}$ of $\mbb{R}$ and $\delta\in(1,\alpha)$,
 \beqlb\label{1.6}
\sup_{t\in(0,T]}t^{\frac{\delta}{2}}\, \mbf{E}\Big\{\sup_{x\in
\mathbb{K}}\tilde{X}_t(x)^\delta\Big\} <\infty.
 \eeqlb
\elemma \proof By Lemma \ref{t1.4}, for each $t\in(0,T]$, there is a
sequence $\{y_t(n):n\ge1\}\subset \mathbb{K}\cup[-1,1]$ so that
$y_t(n)\to0$ as $n\to\infty$ and
 \beqnn
\mbf{E}\{\tilde{X}_t(0)^\delta\}
 \ar=\ar
\mbf{E}\Big\{\lim_{n\to\infty}\tilde{X}_t(y_t(n))^\delta\Big\}
 \le\liminf_{n\to\infty}
\mbf{E}\Big\{\tilde{X}_t(y_t(n))^\delta\Big\} \cr
 \ar\le\ar
\liminf_{n\to\infty} \mbf{E}\Big\{X_t(y_t(n))^\delta\Big\} \le C_T
t^{-\frac\delta2},
 \eeqnn
which implies
 \beqlb\label{1.29}
\sup_{t\in(0,T]}t^{\frac{\delta}{2}}\,
\mbf{E}\big\{\tilde{X}_t(0)^\delta\big\} <\infty.
 \eeqlb
Then the desired result follows from \eqref{1.29}, Lemma \ref{t1.6}
and \cite[Corollary 1.2(iii)]{Walsh}. \qed

By Proposition \ref{t1.9},
$\{X_t(x):t>0,x\in\mbb{R}\}$ satisfies \eqref{1.4} with $X_0\in
M(\mbb{R})$. Similar to \cite[Lemma 2.12]{FMW10} we can prove the
following lemma.
 \blemma\label{t3.2} Fix $\delta,\delta'\in[1,3)$,
$r,r'\in[0,1]$ with $r<\frac{3-\delta}{\delta}$ and
$\frac{\alpha(r'\delta'+\delta')}{2(\alpha+1-p\vee1)}<1$, and a nonempty compact
$\mathbb{K}\subset\mbb{R}$. Define
 \beqnn
\bar{\mathcal{V}}_n:=\int_0^1\frac{1}{2^n}\sum_{k=1}^{2^n}1_{(n_{k-1},n_k]}(s)
\sum_{i=k}^{2^n}(n_i-s)^{-\frac{r\delta+\delta-1}{2}}\bar{\mathcal{V}}_{n,i}(s)ds
 \eeqnn
and
 \beqnn
\bar{\mathcal{U}}_n:=\sup_{1\le i\le
2^n}\int_0^{n_i}(n_i-s)^{-\frac{r'\delta'+\delta'-1}{2}}
\bar{\mathcal{V}}_{n,i}(s) ds
 \eeqnn
with
 \beqnn
\bar{\mathcal{V}}_{n,i}(s):=\sup_{x_1,x_2\in
\mbb{K}}\int_{\mbb{R}}H(X_s(u))^\alpha
|p_{n_i-s}(x_1-u)+p_{n_i-s}(x_2-u)| du.
 \eeqnn
Then we have
 \beqlb\label{1.8}
\sup_{n\ge1} \mbf{P}[\bar{\mathcal{V}}_n\ge C_\varepsilon]\le
\varepsilon,\quad \sup_{n\ge1}\mbf{P}[\bar{\mathcal{U}}_n\ge
C_\varepsilon]\le \varepsilon, \qqquad \varepsilon>0.
  \eeqlb
Moreover, for each $x_1,x_2\in \mbb{K}$ and $n\ge1$ we have
 \beqlb\label{1.7}
 \ar\ar
\int_0^1ds\int_{\mbb{R}}H(X_s(u))^\alpha
\frac{1}{2^n}\sum_{k=1}^{2^n}1_{(n_{k-1},n_k]}(s)\sum_{i=k}^{2^n}
|p_{n_i-s}(x_1-u) \cr
  \ar\ar\qqquad\qqquad\quad\qquad
-p_{n_i-s}(x_2-u)|^\delta du \le
C\bar{\mathcal{V}}_n|x_1-x_2|^{r\delta}
 \eeqlb
and
 \beqlb\label{3.17}
 \ar\ar
\int_0^{n_i}ds\int_{\mbb{R}}H(X_s(u))^\alpha
|p_{n_i-s}(x_1-u) \cr
  \ar\ar\qquad\qquad
-p_{n_i-s}(x_2-u)|^{\delta'} du \le
C\bar{\mathcal{U}}_n|x_1-x_2|^{r'\delta'}, \quad 1\le i\le 2^n.
 \eeqlb
\elemma

\proof We assume that $\mathbb{K}\subset[0,1]$ for simplicity.
Observe that $\mbf{P}$\mbox{-a.s.},
 \beqlb\label{3.16}
\bar{\mathcal{V}}_{n,i}(s)
 \ar\le\ar
C\sup_{x_1,x_2\in[0,1]}\int_{\mbb{R}}[1+X_s(y)^p]
\Big[p_{n_i-s}(x_1-y) +p_{n_i-s}(x_2-y)\Big]dy \cr
 \ar\le\ar
C+C\sup_{x_1,x_2\in[0,1]}\int_{|y|\ge2}X_s(y)^p
\Big[p_{n_i-s}(x_1-y) +p_{n_i-s}(x_2-y)\Big]dy
+C\sup_{|y|\le2}\tilde{X}_s(y)^p \cr
 \ar\le\ar
C+C\int_{\mbb{R}}X_s(y)^p \Big[p_{n_i-s}(y+1) +p_{n_i-s}(y-1)\Big]dy
+C\sup_{|y|\le2}\tilde{X}_s(y)^p.
 \eeqlb
 Then by Lemmas \ref{t1.4} and \ref{t2.1},
 \beqnn
\mbf{E}\{\bar{\mathcal{V}}_{n,i}(s)\}\le C[1+s^{-\frac{p}{2}}].
 \eeqnn
It is elementary to check that
 \beqnn
\sup_{n\ge1}\mbf{E}\{\bar{\mathcal{V}}_n\}
 \ar\le\ar
\sup_{n\ge1}C\int_0^1\frac{1}{2^n}\sum_{k=1}^{2^n}1_{(n_{k-1},n_k]}(s)
\sum_{i=k}^{2^n}(n_i-s)^{-\frac{r\delta+\delta-1}{2}}[1+s^{-\frac{p}{2}}]ds
\cr
 \ar\le\ar
\sup_{n\ge1}\frac{C}{2^n}\sum_{k=1}^{2^n}\sum_{i=k}^{2^n}\int_{n_{k-1}}^{n_k}
(n_i-s)^{-\frac{r\delta+\delta-1}{2}}[1+s^{-\frac{p}{2}}]ds \cr
 \ar=\ar
\sup_{n\ge1}\frac{C}{2^n}\sum_{i=1}^{2^n}\int_0^{n_i}
(n_i-s)^{-\frac{r\delta+\delta-1}{2}}[1+s^{-\frac{p}{2}}]ds<\infty.
 \eeqnn
Then the first assertion of \eqref{1.8} follows from  the Markov inequality.
Using
\eqref{1.35} one gets \eqref{1.7} and \eqref{3.17}.

We now prove the second assertion of \eqref{1.8}. Observe that for
each constant $\theta>0$,
 \beqnn
(n_i-s)^{-\theta}\int_{y\ge2}X_s(y)^pp_{n_i-s}(y-1)dy
\le
C\int_{\mbb{R}}X_s(y)^pp_2(y-1)dy
 \eeqnn
and
 \beqnn
(n_i-s)^{-\theta}\int_{y\le-2}X_s(y)^pp_{n_i-s}(y+1)dy
\le
C\int_{\mbb{R}}X_s(y)^pp_2(y+1)dy.
 \eeqnn
Then it is easy to check that
 \beqlb\label{2.14}
 \ar\ar
(n_i-s)^{-\frac{r'\delta'+\delta'-1}{2}}\bar{\mathcal{V}}_{n,i}(s)
\cr
 \ar\ar\quad \le
C(n_i-s)^{-\frac{r'\delta'+\delta'-1}{2}}\Big\{
1+\int_{y\ge2}X_s(y)^pp_{n_i-s}(y-1)dy
+\int_{y\le-2}X_s(y)^pp_{n_i-s}(y+1)dy \cr
 \ar\ar\qquad
+\sup_{x_1,x_2\in[0,1]}\int_{|y|\le2}X_s(y)^p
\big[p_{n_i-s}(x_1+y) +p_{n_i-s}(x_2+y)\big]dy\Big\} \cr
 \ar\ar\quad \le
C(n_i-s)^{-\frac{r'\delta'+\delta'-1}{2}}
+C\int_{\mbb{R}}X_s(y)^p[p_2(y+1)+p_2(y-1)]dy \cr
 \ar\ar\qquad
+(n_i-s)^{-\frac{r'\delta'+\delta'-1}{2}}
\sup_{x_1,x_2\in[0,1]}\int_{|y|\le2}X_s(y)^p \big[p_{n_i-s}(x_1+y)
+p_{n_i-s}(x_2+y)\big]dy\Big\}
 \eeqlb
 for each $1\le i\le 2^n$. Observe that for each $x\in[0,1]$,
 \beqnn
\int_{|y|\le2}X_s(y)^pp_{n_i-s}(x+y)dy
\le
\Big|\int_{|y|\le2}X_s(y)p_{n_i-s}(x+y)dy\Big|^p
\le[2\pi(n_i-s)^{-\frac12}\<X_s,1\>]^p
 \eeqnn
for  $0<p\le 1$ and
 \beqnn
\int_{|y|\le2}X_s(y)^pp_{n_i-s}(x+y)dy
\le
[2\pi(n_i-s)]^{-\frac{1}{2}}\<X_s,1\> \sup_{|y|\le2}\tilde{X}_s(y)^{p-1},
\,\,\,\mbf{P}\mbox{-a.s.}
 \eeqnn
for  $1< p<2$. Combining with \eqref{2.14} we have
 \beqlb\label{2.15}
\bar{\mathcal{U}}_n
 \ar\le\ar
C\sup_{1\le i\le 2^n}
\int_0^{n_i}(n_i-s)^{-\frac{r'\delta'+\delta'-1}{2}}ds
+C\int_0^1ds\int_{\mbb{R}}X_s(y)^p[p_2(y+1)+p_2(y-1)]dy \cr
 \ar\ar
+C\sup_{s\in(0,1]}\<X_s,1\>^p \sup_{1\le i\le
2^n}\int_0^{n_i}(n_i-s)^{-\frac{r'\delta'+\delta'+p-1}{2}}ds \cr
 \ar\le\ar
C+C\int_0^1ds\int_{\mbb{R}}X_s(y)^p[p_2(y+1)+p_2(y-1)]dy
+C\sup_{s\in(0,1]}\<X_s,1\>^p
 \eeqlb
for  $0<p\le 1$ and
 \beqnn
\bar{\mathcal{U}}_n
 \ar\le\ar
C\sup_{1\le i\le 2^n}
\int_0^{n_i}(n_i-s)^{-\frac{r'\delta'+\delta'-1}{2}}ds
+C\int_0^1ds\int_{\mbb{R}}X_s(y)^p[p_2(y+1)+p_2(y-1)]dy \cr
 \ar\ar
+C\sup_{s\in(0,1]}\<X_s,1\> \sup_{1\le i\le
2^n}\int_0^{n_i}(n_i-s)^{-\frac{r'\delta'+\delta'}{2}}
\sup_{|y|\le2}\tilde{X}_s(y)^{p-1}ds \cr
 \ar\le\ar
C+C\int_0^1ds\int_{\mbb{R}}X_s(y)^p[p_2(y+1)+p_2(y-1)]dy \cr
 \ar\ar
+C\sup_{s\in(0,1]}\<X_s,1\>\sup_{1\le i\le 2^n}
\int_0^{n_i}(n_i-s)^{-\frac{r'\delta'+\delta'}{2}}
\sup_{|y|\le2}\tilde{X}_s(y)^{p-1}ds,\,\,\,\mbf{P}\mbox{-a.s.}
 \eeqnn
for  $1<p<2$. Taking $\alpha_-\in(1,\alpha)$ with
$\frac{\alpha_-(r'\delta'+\delta')}{2(\alpha_-+1-p)}<1$, we have
 \beqnn
 \ar\ar
\int_0^{n_i}(n_i-s)^{-\frac{r'\delta'+\delta'}{2}}
\sup_{|y|\le2}\tilde{X}_s(y)^{p-1}ds \cr
 \ar\ar\quad\le
\Big|\int_0^{n_i}(n_i-s)^{-\frac{\alpha_-(r'\delta'+\delta')}{2(\alpha_-+1-p)}}ds
\Big|^{\frac{\alpha_-+1-p}{\alpha_-}}
\Big|\int_0^{n_i}\sup_{|y|\le2}\tilde{X}_s(y)^{\alpha_-}ds\Big|^{\frac{p-1}{\alpha_-}} \cr
 \ar\ar\quad\le
C\Big|\int_0^1\sup_{|y|\le2}\tilde{X}_s(y)^{\alpha_-}ds\Big|^{\frac{\alpha_-}{p-1}}
 \eeqnn
for the case $1<p<2$. It then follows that
 \beqlb\label{2.16}
\bar{\mathcal{U}}_n
 \ar\le\ar
C+C\int_0^1ds\int_{\mbb{R}}X_s(y)^p[p_2(y+1)+p_2(y-1)]dy  \cr
 \ar\ar
+C\sup_{s\in(0,1]}\<X_s,1\>
\Big|\int_0^1\sup_{|y|\le2}X_s(y)^{\alpha_-}ds\Big|^{\frac{\alpha_-}{p-1}}
 \eeqlb
for the case $1<p<2$.
By Lemmas \ref{t1.4} and \ref{t2.1},
  \beqnn
\mbf{E}\Big\{\int_0^1ds\int_{\mbb{R}}X_s(y)^p[p_2(y+1)+p_2(y-1)]dy\Big\}
\le
C\int_0^1s^{-\frac p2}ds<\infty
 \eeqnn
and
 \beqnn
\mbf{E}\Big\{\int_0^1\sup_{|y|\le2}X_s(y)^{\alpha_-}ds\Big\}
\le\int_0^1\mbf{E}\Big\{\sup_{|y|\le2}X_s(y)^{\alpha_-}\Big\}ds
\le C\int_0^1s^{-\alpha_-/2}ds<\infty,
 \eeqnn
which imply
 \beqnn
\int_0^1ds\int_{\mbb{R}}X_s(y)^p[p_2(y+1)+p_2(y-1)]dy
+\int_0^1\sup_{|y|\le2}X_s(y)^{\alpha_-}ds<\infty, \,\,\,\,
\mbf{P}\mbox{-a.s.}
 \eeqnn
 Combining with \eqref{2.15}--\eqref{2.16} and the fact
 \beqnn
\sup_{s\in(0,1]}\<X_s,1\><\infty, \,\,\, \mbf{P}\mbox{-a.s.},
 \eeqnn
 we have
 \beqnn
\sup_{n\ge1}\bar{\mathcal{U}}_n<\infty, \,\,\, \mbf{P}\mbox{-a.s.},
 \eeqnn
 which implies the second assertion of \eqref{1.8}. \qed

For $t\ge0$ and $\psi\in B(\mbb{R})$ define  discontinuous
martingales
 \beqnn
t\mapsto
M_t^1(\psi):=\int_0^t\int_0^\infty\int_{\mbb{R}}\int_0^{H(X_{s-}(u))^\alpha}
z\psi(u)1_{\{|u|\le \mbb{K}_1\}} \tilde{N}_0(ds,dz,du,dv)
 \eeqnn
and
 \beqnn
t\mapsto M_t^2(\psi):=\int_0^t\int_0^\infty\int_{\mbb{R}}
\int_0^{H(X_{s-}(u))^\alpha} \frac{z\psi(u)}{|u|}1_{\{|u|>
\mbb{K}_1\}} \tilde{N}_0(ds,dz,du,dv),
 \eeqnn
where $\mbb{K}_1:=\mbb{K}_0+1$ with $\mbb{K}_0:=\sup_{x\in
\mbb{K}}|x|$. For $i=1,2$ let $\Delta M_s^i(y)$ denote the jumps of
$M^i(ds,dy)$. Similar to \cite[Lemma 2.14]{FMW10} one can show the
following result.

 \blemma\label{t3.3} Let
$\gamma\in(0,\alpha^{-1})$ and $\lambda:=\alpha^{-1}-\gamma$. Then
for each $\varepsilon>0$ there exists a constant $C_\varepsilon>0$
independent of $n$ so that
 \beqnn
\mbf{P}\Big(\bigcup_{k=1}^{2^n} \Big\{\Delta M_s^1(y)> 2^{\lambda
n}C_\varepsilon(n_k-s)^\lambda \mbox{~for some
}s\in[n_{k-1},n_k)\mbox{~and~}|y|\le \mathbb{K}_1 \Big\}\Big) \le
\varepsilon
 \eeqnn
and
 \beqnn
\mbf{P}\Big(\bigcup_{k=1}^{2^n} \Big\{\Delta M_s^2(y)> 2^{\lambda
n}C_\varepsilon(n_k-s)^\lambda \mbox{~for some
}s\in[n_{k-1},n_k)\mbox{~and~} |y|> \mathbb{K}_1\Big)\le
\varepsilon.
 \eeqnn
\elemma \proof Since the proofs are similar, we only present the
first one. Let $c>0$. For $n\ge1$ and $1\le k\le 2^n$ put
 \beqnn
Y_{n,k}^1 := N_0\Big((s,z,u,v):s\in[n_{k-1},n_k),z\ge 2^{\lambda
n}c(n_k-s)^\lambda, |u|\le \mathbb{K}_1,v\le
H(X_{s-}(u))^\alpha\Big).
 \eeqnn
Then by the Markov inequality,
 \beqnn
\mbf{P}\Big\{\frac{\Delta M_s^1(y)}{2^{\lambda n}(n_k-s)^\lambda}> c
\mbox{~for some }s\in[n_{k-1},n_k)\mbox{~and~}|y|\le
\mathbb{K}_1\Big\} =\mbf{P}\{Y_{n,k}^1\ge1\} \le
\mbf{E}\{Y_{n,k}^1\}.
 \eeqnn
By Lemma \ref{t1.4},
 \beqnn
\mbf{E}\{Y_{n,k}^1\}
 \ar=\ar
\mbf{E}\Big\{\int_{n_{k-1}}^{n_k}ds\int_0^\infty
m_0(dz)\int_{-\mbb{K}_1}^{\mbb{K}_1}du \int_0^{H(X_s(u))^\alpha}
1_{\{z\ge c2^{\lambda n}(n_k-s)^\lambda\}} dv \Big\}\cr
 \ar=\ar
Cc^{-\alpha}\int_{n_{k-1}}^{n_k}ds\int_{-\mbb{K}_1}^{\mbb{K}_1}
\mbf{E}\{H(X_s(u))^\alpha\} 2^{-\alpha\lambda
n}(n_k-s)^{-\alpha\lambda} du \cr
 \ar\le\ar
Cc^{-\alpha}\int_{n_{k-1}}^{n_k}ds\int_{-\mbb{K}_1}^{\mbb{K}_1}\mbf{E}\{1+X_s(u)^p\}
2^{-\alpha\lambda n}(n_k-s)^{-\alpha\lambda} du \cr
 \ar\le\ar
Cc^{-\alpha}2^{-\alpha\lambda n}\int_{n_{k-1}}^{n_k}
(n_k-s)^{-\alpha\lambda}[1+s^{-\frac{p}{2}}]ds
\le Cc^{-\alpha}n_{k-1}^{-p/2}2^{- n}
 \eeqnn
for $2\le k\le 2^n$.  Similarly
 \beqnn
\mbf{E}\{Y_{n,1}^1\} \le Cc^{-\alpha}2^{-\alpha\lambda
n}\int_0^{n_1} (n_1-s)^{-\alpha\lambda}[1+s^{-\frac{p}{2}}]ds \le
Cc^{-\alpha}n_1^{-p/2}2^{- n},
 \eeqnn
Thus
 \beqnn
\sum_{k=1}^{2^n}\mbf{E}\{Y_{n,k}^1\}\le Cc^{-\alpha}.
 \eeqnn
The desired result then follows. \qed

Let
$\mathscr{L}$
be the space of measurable functions $\psi$ on $\mbb{R}_+\times\mbb{R}$
so that
 \beqnn
\int_0^ts^{-\frac{p}{2}}ds\int_{\mbb{R}}
\Big[|\psi(s,x)|+|\psi(s,x)|^2\Big]dx<\infty,\qquad t>0. \eeqnn
Similar to \cite[Lemma 2.15]{FMW10}, we have the next result.
\blemma\label{t3.1} Given $\psi\in \mathscr{L}$ with $\psi\ge0$,
there exist a spectrally positive $\alpha$-stable process
$\{L_t:t\ge0\}$  so that for $t\ge0$,
 \beqnn
Z_t(\psi):=
\int_0^t\int_0^\infty\int_{\mbb{R}}\int_0^\infty
1_{\{v\le H(X_s(u))^\alpha\}}z\psi(s,u)\tilde{N}_0(ds,dz,du,dv)
=L_{T'(t)},
 \eeqnn
where
 \beqnn
T'(t):=\int_0^tds\int_{\mbb{R}}
[H(X_s(u))\psi(s,u)]^\alpha du.
 \eeqnn
\elemma
The proof is similar to that of \cite[Lemma 2.15]{FMW10}.

\noindent{\it Proof of Theorem \ref{t1.1}.} By Lemma \ref{t1.7} we
only consider the case $\alpha<3/2$. Since the proof of \eqref{2.19}
is essentially same to that of \cite[Theroem 1.2(a)]{FMW10} based on
equation \eqref{1.4}, then we only state the proof of \eqref{2.17},
which is a modification of that of \cite[Theroem 1.2(a)]{FMW10} and
proceeds as follows. Let $X_0\in M(\mbb{R})$ be fixed. We assume
that $T=1$ in this proof. Recall \eqref{2.18}  and
$\lambda=\frac1\alpha-\gamma$ with $\gamma\in(0,\alpha^{-1})$. Also
recall that $n_k=\frac{k}{2^n}$ for $n\ge1$ and $0\le k\le 2^n$. Let
$f(x)^+$ and $f(x)^-$ be, respectively, the positive part and the
negative part of  $f(x)$. For $x_1, x_2\in \mbb{R}$ and $s\in[0,1]$
define
 \beqnn
\mathcal{U}_s(x_1,x_2):=\int_0^s\int_0^\infty\int_{\mbb{R}}
\int_0^{H(X_{s-}(u))^\alpha} z
[p_{s-s_1}(x_1-u)-p_{s-s_1}(x_2-u)]^+\tilde{N}_0(ds_1,dz,du,dv).
 \eeqnn
Let $\mathcal{V}_s(x_1,x_2)$ be defined as
$\mathcal{U}_s(x_1,x_2)$ with
$[p_{s-s_1}(x_1-u)-p_{s-s_1}(x_2-u)]^+$ replaced by
$[p_{s-s_1}(x_1-u)-p_{s-s_1}(x_2-u)]^-$. Then
 \beqlb\label{3.26}
Z_s(x_1)-Z_s(x_2)
=[\mathcal{U}_s^+(x_1,x_2)-\mathcal{U}_s^-(x_1,x_2)]
-[\mathcal{V}_s^+(x_1,x_2)-\mathcal{V}_s^-(x_1,x_2)].
 \eeqlb

Observe that
 \beqnn
 \ar\ar
\frac{1}{2^n}\sum_{k=1}^{2^n}\mathcal{U}_{n_k}(x_1,x_2) \cr
 \ar=\ar
\frac{1}{2^n}\sum_{k=1}^{2^n}\int_0^{n_k}
\int_0^\infty\int_{\mbb{R}}\int_0^{H(X_{s-}(u))^\alpha}
z[p_{n_k-s}(x_1-u)
-p_{n_k-s}(x_2-u)]^+\tilde{N}_0(ds,dz,du,dv) \cr
 \ar=\ar
\int_0^1
\int_0^\infty\int_{\mbb{R}}\int_0^{H(X_{s-}(u))^\alpha}
\frac{1}{2^n}\sum_{k=1}^{2^n}1_{(0,n_k]}(s)
z[p_{n_k-s}(x_1-u) \cr
 \ar\ar\qqquad\qqquad\qqquad\qqquad
-p_{n_k-s}(x_2-u)]^+\tilde{N}_0(ds,dz,du,dv)   \cr
 \ar=\ar
\int_0^1\int_0^\infty\int_{\mbb{R}}\int_0^{H(X_{s-}(u))^\alpha}
\frac{1}{2^n}\sum_{k=1}^{2^n}1_{(n_{k-1},n_k]}(s)\sum_{i=k}^{2^n}
z[p_{n_i-s}(x_1-u) \cr
 \ar\ar\qqquad\qqquad\qqquad\qqquad
-p_{n_i-s}(x_2-u)]^+\tilde{N}_0(ds,dz,du,dv) \cr
 \ar=\ar
\int_0^1\int_0^\infty\int_{\mbb{R}}\int_0^\infty
1_{\{v\le H(X_{s-}(u))^\alpha\}}
z\psi_n(s,u)\tilde{N}_0(ds,dz,du,dv),
 \eeqnn
where
 \beqnn
\psi_n(s,u):=\frac{1}{2^n}\sum_{k=1}^{2^n}1_{(n_{k-1},n_k]}(s)\sum_{i=k}^{2^n}
[p_{n_i-s}(x_1-u)-p_{n_i-s}(x_2-u)]^+.
 \eeqnn
One can see that $\psi_n(s,u)$ satisfies
the assumptions of Lemma \ref{t3.1}, and there is a stable process
$\{L_t:t\ge0\}$ so that
 \beqlb\label{1.9}
\frac{1}{2^n}\sum_{k=1}^{2^n}\mathcal{U}_{n_k}(x_1,x_2) =L_{T_n},
 \eeqlb
where
\beqnn
T_n:=\int_0^1ds\int_{\mbb{R}}[H(X_s(u))\psi_n(s,u)]^\alpha du.
 \eeqnn
Let $\varepsilon\in(0,1)$ be fixed. Let $\bar{\mathcal{V}}_n$ and
$\bar{\mathcal{U}}_n$ be defined in Lemma \ref{t3.2}. Then
 \beqlb\label{3.11}
\sup_{n\ge1}\Big\{\mbf{P}(\bar{\mathcal{V}}_n>
C_\varepsilon)+\mbf{P}(\bar{\mathcal{U}}_n> C_\varepsilon)\Big\}
\le2\varepsilon.
 \eeqlb

Set
 \beqnn
A_n^\varepsilon
 \ar:=\ar
\bigcap_{k=1}^{2^n} \Bigg(\Big\{\frac{\Delta M_s^1(y)}{2^{\lambda n}
(n_k-s)^\lambda}\le C_\varepsilon \mbox{~for all
}s\in[n_{k-1},n_k)\mbox{~and~}|y|\le \mathbb{K}_1\Big\} \cr
 \ar\ar\qquad
\bigcap \Big\{\frac{\Delta M_s^2(y)}{2^{\lambda n}
(n_k-s)^\lambda}\le C_\varepsilon \mbox{~for all
}s\in[n_{k-1},n_k)\mbox{~and~}|y|>\mathbb{K}_1\Big\}\Bigg)  \cr
 \ar\ar
\bigcap\Big\{\bar{\mathcal{V}}_n\le
C_\varepsilon,\bar{\mathcal{U}}_n\le C_\varepsilon\Big\}.
 \eeqnn
By Lemmas \ref{t3.2} and \ref{t3.3},
 \beqnn
\sup_{n\ge1}\mbf{P}(A_n^{\varepsilon,c})\le4\varepsilon,
 \eeqnn
where $A_n^{\varepsilon,c}$ denotes the complement of
$A_n^\varepsilon$. Define
$\mathcal{U}_s^{n,\varepsilon}(x_1,x_2):=\mathcal{U}_s(x_1,x_2)1_{A_n^\varepsilon}$.

In order to complete our proof we need to establish the following two lemmas.
 \blemma\label{t3.5}
For each $n\ge1$ and $r>0$,
 \beqnn
\mbf{P}\Big\{ \frac{1}{2^n}\sum_{k=1}^{2^n}
\mathcal{U}_{n_k}^{n,\varepsilon}(x_1,x_2) \ge r|x_1-x_2|^\eta\Big\}
\le (C_\varepsilon r^{-1}|x_1-x_2|)^{C'_\varepsilon
r|x_1-x_2|^{(\eta-\eta_c)/2}}.
 \eeqnn
 \elemma
\proof
By \eqref{1.9}, we have
 \beqlb\label{1.27}
\mbf{P}\Big\{ \frac{1}{2^n}\sum_{k=1}^{2^n}
\mathcal{U}_{n_k}^{n,\varepsilon}(x_1,x_2) \ge r|x_1-x_2|^\eta\Big\}
= \mbf{P}\Big\{L_{T_n}>r|x_1-x_2|^\eta,A_n^\varepsilon\Big\}.
 \eeqlb
Note that on event $A_n^\varepsilon$ the jumps of $M_s^1(x)$ do not
exceed
 \beqnn
C_\varepsilon2^{\lambda n}
(n_k-s)^\lambda,\qquad s\in[n_{k-1},n_k).
 \eeqnn
Then the jumps of
 \beqnn
(0,1)\ni l\mapsto\int_0^l\int_0^\infty\int_{\mbb{R}}\int_0^\infty
1_{\{v\le H(X_{s-}(u))^\alpha,|u|\le \mbb{K}_1\}}
z\psi_n(s,u)\tilde{N}_0(ds,dz,du,dv)
 \eeqnn
are bounded by
 \beqlb\label{3.10}
I_n
 \ar:=\ar
C_\varepsilon2^{\lambda n}
\sup_{1\le k\le 2^n}\sup_{(s,y)\in [n_{k-1},n_k)\times\mbb{R}}
(n_k-s)^\lambda\psi_n(s,y) \cr
 \ar\le\ar
C_\varepsilon2^{\lambda n}
\sup_{1\le k\le 2^n}\sup_{(s,y)\in [n_{k-1},n_k)\times\mbb{R}}
\frac{1}{2^n}\sum_{i=k}^{2^n}
(n_k-s)^\lambda|p_{n_i-s}(x_1-y)-p_{n_i-s}(x_2-y)|.
 \eeqlb
Applying \eqref{1.35} with $\theta=\eta_c-2\gamma$ gives
 \beqnn
 \ar\ar
\sup_{y\in \mbb{R}}\frac{1}{2^n}\sum_{i=k}^{2^n}
(n_k-s)^\lambda|p_{n_i-s}(x_1-y)-p_{n_i-s}(x_2-y)| \cr
 \ar\ar\qquad\le
C|x_1-x_2|^{\eta_c-2\gamma}
\frac{1}{2^n}\sum_{i=k}^{2^n}(n_i-s)^{-\eta_c/2+\gamma}(n_k-s)^\lambda
\sup_{y\in \mbb{R}}p_{n_i-s}(y) \cr
 \ar\ar\qquad\le
C|x_1-x_2|^{\eta_c-2\gamma}
\frac{1}{2^n}\sum_{i=k}^{2^n}
\Big(\frac{n_k-s}{n_i-s}\Big)^\lambda \cr
 \ar\ar\qquad\le
C|x_1-x_2|^{\eta_c-2\gamma}
\frac{1}{2^n}\sum_{i=k}^{2^n}
\Big(\frac{1}{i-k+1}\Big)^\lambda
\le C|x_1-x_2|^{\eta_c-2\gamma}2^{-\lambda n}
 \eeqnn
for $s\in[n_{k-1},n_k)$.
This implies
 \beqlb\label{3.12}
I_n\le
C_\varepsilon |x_1-x_2|^{\eta_c-2\gamma}.
 \eeqlb

Similarly, one can see that the jumps of
 \beqnn
(0,1)\ni l\mapsto\int_0^l\int_0^\infty\int_{\mbb{R}}\int_0^\infty
1_{\{v\le H(X_{s-}(u))^\alpha,|u|>\mbb{K}_1\}}
z\psi_n(s,u)\tilde{N}_0(ds,dz,du,dv)
 \eeqnn
are bounded by
 \beqnn
C_\varepsilon |x_1-x_2|^{\eta_c-2\gamma}.
 \eeqnn
Combining with \eqref{3.12} we conclude that the jumps of
 \beqnn
(0,1)\ni l\mapsto\int_0^l\int_0^\infty\int_{\mbb{R}}\int_0^\infty
1_{\{v\le H(X_{s-}(u))^\alpha\}}
z\psi_n(s,u)\tilde{N}_0(ds,dz,du,dv)
 \eeqnn
on $A_n^\varepsilon$ are bounded by
 \beqnn
C_\varepsilon |x_1-x_2|^{\eta_c-2\gamma}.
 \eeqnn

Observe that
 \beqlb\label{3.4}
 \ar\ar
\mbf{P}\Big\{L_{T_n}\ge r|x_1-x_2|^\eta,A_n^\varepsilon\Big\} \cr
 \ar\ar\qquad=
\mbf{P}\Big\{L_{T_n}\ge r|x_1-x_2|^\eta,
\sup_{u<T_n}\Delta L_u\le C_\varepsilon|x_1-x_2|^{\eta_c-2\gamma},
A_n^\varepsilon\Big\} \cr
 \ar\ar\qquad\le
\mbf{P}\Big\{\sup_{v\le T_n} L_v
1_{\{\sup_{u<v}\Delta L_u\le C_\varepsilon|x_1-x_2|^{\eta_c-2\gamma}\}}
\ge r|x_1-x_2|^\eta,A^\varepsilon_n\Big\}.
 \eeqlb
Moreover,
 \beqnn
T_n
 \ar\le\ar
\int_0^1ds\int_{\mbb{R}}H(X_s(u))^\alpha
\sum_{k=1}^{2^n}1_{(n_{k-1},n_k]}(s)\Big|\frac{1}{2^n}\sum_{i=k}^{2^n}
|p_{n_i-s}(x_1-u)-p_{n_i-s}(x_2-u)|\Big|^\alpha du \cr
 \ar\le\ar
\int_0^1ds\int_{\mbb{R}}H(X_s(u))^\alpha
\sum_{k=1}^{2^n}1_{(n_{k-1},n_k]}(s)\frac{1}{2^n}\sum_{i=k}^{2^n}
|p_{n_i-s}(x_1-u)-p_{n_i-s}(x_2-u)|^\alpha du.
 \eeqnn
Applying Lemma \ref{t3.2} with $\delta=\alpha$
and $r=1$ one gets
 \beqnn
T_n \le C_\varepsilon |x_1-x_2|^\alpha \qqquad
\mbox{on~}\{\bar{\mathcal{V}}_n\le C_\varepsilon\}.
 \eeqnn
Combining this with \eqref{3.4} we have
 \beqnn
 \ar\ar
\mbf{P}\Big\{L_{T_n}\ge r|x_1-x_2|^\eta,
A_n^\varepsilon\Big\} \cr
 \ar\ar\quad\le
\mbf{P}\Big\{\sup_{v\le C_\varepsilon|x_1-x_2|^\alpha} L_v
1_{\{\sup_{u<v}\Delta L_u\le C_\varepsilon|x_1-x_2|^{\eta_c-2\gamma}\}}
\ge r|x_1-x_2|^\eta\Big\}.
 \eeqnn
Using (3.14) in \cite{FMW10}, and \cite[Lemma 2.3]{FMW10} with
$\kappa=\alpha$, $t=C_\varepsilon |x_1-x_2|^\alpha$,
$x=r|x_1-x_2|^\eta$, and $y=C_\varepsilon
|x_1-x_2|^{\eta_c-2\gamma}$, one obtains
 \beqlb\label{1.28}
\mbf{P}\Big\{L_{T_n}\ge r|x_1-x_2|^\eta, A_n^\varepsilon\Big\} \le
\Big( C_\varepsilon r^{-1}|x_1-x_2|^{2\alpha-2}
\Big)^{C'_\varepsilon r|x_1-x_2|^{\eta-\eta_c+2\gamma}}.
 \eeqlb
Taking $\gamma:=\frac{\eta_c-\eta}{4}$, we have
 \beqnn
\mbf{P}\Big\{L_{T_n}\ge r|x_1-x_2|^\eta, A_n^\varepsilon\Big\}\le
(C_\varepsilon r^{-1}|x_1-x_2|)^{C'_\varepsilon
r|x_1-x_2|^{(\eta-\eta_c)/2}},
 \eeqnn
which together with \eqref{1.27} proves the lemma.
\qed

 \blemma\label{t3.4}
For each $n\ge1$ and $1\le i\le 2^n$ let $\{L_{n_i}(t):t\ge0\}$
be a spectrally positive $\alpha$-stable process. Let
$L_{n_i}^-(t)$ be the negative part of $L_{n_i}(t)$, and $T(t)$
be defined as in Lemma \ref{t3.1} with $\psi(s,u)$ replaced by
$[p_{t-s}(x_1-u)-p_{t-s}(x_2-u)]^+$. Then for each $x>0$ and
$n\ge1$,
 \beqnn
\mbf{P}\Big\{\frac{1}{2^n}\sum_{i=1}^{2^n}L_{n_i}^-(T(n_i))>x,
\bar{\mathcal{U}}_n\le C_\varepsilon\Big\} \le
C_\varepsilon\exp\Big\{ -\frac{C'_\varepsilon
x^{\alpha/(\alpha-1)}}{|x_1-x_2|^{r'\delta'/(\alpha-1)}}\Big\},
 \eeqnn
where $r'$ and $\delta'$ are defined in Lemma \ref{t3.2}.
 \elemma
\proof
It is easy to see for all $h>0$,
 \beqlb\label{3.23}
 \ar\ar
\mbf{P}\Big\{\frac{1}{2^n}\sum_{i=1}^{2^n}L_{n_i}^-(T(n_i))>x,
\bar{\mathcal{U}}_n\le C_\varepsilon\Big\} \cr
 \ar\ar\qquad=
\mbf{P}\Big\{\exp\Big[\frac{h}{2^n}\sum_{i=1}^{2^n}L_{n_i}^-(T(n_i))\Big]>e^{hx},
\bar{\mathcal{U}}_n\le C_\varepsilon\Big\} \cr
  \ar\ar\qquad\le
e^{-hx}\mbf{E}\Big\{\exp\Big[\frac{h}{2^n}\sum_{i=1}^{2^n}L_{n_i}^-(T(n_i))\Big]
1_{\{\bar{\mathcal{U}}_n\le C_\varepsilon\}}\Big\} \cr
  \ar\ar\qquad\le
e^{-hx}\prod_{i=1}^{2^n}\Big|\mbf{E}\Big[\exp\Big(h
L_{n_i}^-(T(n_i))\Big)1_{\{\bar{\mathcal{U}}_n\le
C_\varepsilon\}}\Big]\Big|^{\frac{1}{2^n}}.
 \eeqlb
Observe that
 \beqlb\label{3.18}
 \ar\ar
\mbf{E}\Big[\exp\Big(h
L_{n_i}^-(T(n_i))\Big)1_{\{\bar{\mathcal{U}}_n\le
C_\varepsilon\}}\Big] \cr
 \ar=\ar
\mbf{P}\Big[L_{n_i}^-(T(n_i))=0,\bar{\mathcal{U}}_n\le
C_\varepsilon\Big] +\mbf{E}\Big[\exp\Big(h
L_{n_i}^-(T(n_i))\Big)1_{\{\bar{\mathcal{U}}_n\le
C_\varepsilon,L_{n_i}^-(T(n_i))>0\}}\Big] \cr
 \ar\le\ar
1+ \int_0^\infty
e^{hy}\mbf{P}\big[L_{n_i}^-(T(n_i))>y,\bar{\mathcal{U}}_n\le
C_\varepsilon\big]dy.
 \eeqlb
By Lemma \ref{t3.2},
  \beqnn
T(n_i)\le C\bar{\mathcal{U}}_n|x_1-x_2|^{r'\delta'}\le
C_\varepsilon |x_1-x_2|^{r'\delta'}
 \eeqnn
on $\{\bar{\mathcal{U}}_n\le C_\varepsilon\}$. Then using Lemma 2.4
of \cite{FMW10}, for each $y>0$,
 \beqlb\label{3.19}
 \ar\ar
\mbf{P}\big\{L_{n_i}^-(T(n_i))>y,\bar{\mathcal{U}}_n\le
C_\varepsilon\big\}
 =
\mbf{P}\Big\{L_{n_i}(T(n_i))<-y,\bar{\mathcal{U}}_n\le
C_\varepsilon\Big\} \cr
 \ar\ar\quad\le
\mbf{P}\Big\{\inf_{u\le C_\varepsilon
|x_1-x_2|^{r'\delta'}}L_{n_i}(u)<-y\Big\} \le
\exp\Big\{-C_\varepsilon
y^{\alpha/(\alpha-1)}|x_1-x_2|^{-r'\delta'/(\alpha-1)}\Big\}.
 \eeqlb
Since for all $a,b\ge0$,
 \beqnn
ab\le (1- \alpha^{-1})a^{\alpha/(\alpha-1)} + \alpha^{-1} b^\alpha,
 \eeqnn
then combining \eqref{3.18} and \eqref{3.19} we have
 \beqnn
 \ar\ar
\mbf{E}\Big[\exp\Big(h
L_{n_i}^-(T(n_i))\Big)1_{\{\bar{\mathcal{U}}_n\le
C_\varepsilon\}}\Big] \cr
 \ar\ar\quad\le
1+\int_0^\infty e^{hy}\exp\Big\{-C_\varepsilon
y^{\alpha/(\alpha-1)}|x_1-x_2|^{-r'\delta'/(\alpha-1)}\Big\}dy \cr
 \ar\ar\quad\le
1+\exp\big\{C_\varepsilon'|x_1-x_2|^{r'\delta'}h^\alpha\big\}
\int_0^\infty \exp\Big\{-C_\varepsilon
y^{\alpha/(\alpha-1)}|x_1-x_2|^{-r'\delta'/(\alpha-1)}\Big\}dy \cr
 \ar\ar\quad\le
1+C_\varepsilon\exp\big\{C'_\varepsilon|x_1-x_2|^{r'\delta'}h^\alpha\big\}
\le
C_\varepsilon\exp\big\{C'_\varepsilon|x_1-x_2|^{r'\delta'}h^\alpha\big\}.
 \eeqnn
Then it follows from \eqref{3.23} that
 \beqnn
\mbf{P}\Big\{\frac{1}{2^n}\sum_{i=1}^{2^n}L_{n_i}^-(T(n_i))>x,
\bar{\mathcal{U}}_n\le C_\varepsilon\Big\} \le
C_\varepsilon\exp\big\{C'_\varepsilon|x_1-x_2|^{r'\delta'}h^\alpha-hx\big\}.
 \eeqnn
Minimizing the function
$h\mapsto C_\varepsilon'|x_1-x_2|^{r'\delta'}h^\alpha-hx$,
the desired result follows.
\qed

We now return to the proof of Theorem \ref{t1.1}. By Lemma
\ref{t3.1}, for each $0<s\le1$, there exists a spectrally positive
$\alpha$-stable processes $\{L_s(t):t\ge0\}$ so that
 \beqnn
\mathcal{U}_s(x_1,x_2)=L_s(T(s)),
 \eeqnn
where
 \beqnn T(s)=\int_0^sds_1\int_{\mbb{R}}
\Big[H(X_{s_1}(u))[p_{s-s_1}(x_1-u)-p_{s-s_1}(x_2-u)]^+\Big]^\alpha
du.
 \eeqnn
Since $p<1+\alpha(\alpha-1)/2$,
then applying Lemma \ref{t3.4} with $r'=2-\alpha$ and $\delta'=1$, we get
 \beqlb\label{3.25}
\mbf{P}\Big\{\frac{1}{2^n}\sum_{i=1}^{2^n}\mathcal{U}_{n_i}(x_1,x_2)^->
r|x_1-x_2|^\eta,A_n^\varepsilon\Big\} \le C_\varepsilon\exp\Big\{
-\frac{C'_\varepsilon
r^{\alpha/(\alpha-1)}}{|x_1-x_2|^{(2-\alpha-\eta\alpha)/(\alpha-1)}}\Big\}.
 \eeqlb
Observe that
 \beqnn
\frac{1}{2^n}\sum_{i=1}^{2^n}|\mathcal{U}_{n_i}(x_1,x_2)|
=\frac{1}{2^n}\sum_{i=1}^{2^n}\mathcal{U}_{n_i}(x_1,x_2)
+\frac{2}{2^n}\sum_{i=1}^{2^n}\mathcal{U}_{n_i}(x_1,x_2)^-,
 \eeqnn
which together with Lemma \ref{t3.5} and \eqref{3.25} implies
 \beqnn
 \ar\ar
\mbf{P}\Big\{\frac{1}{2^n}\sum_{i=1}^{2^n}|\mathcal{U}_{n_i}(x_1,x_2)|>
2r|x_1-x_2|^\eta,A_n^\varepsilon\Big\} \cr
 \ar\le\ar
\mbf{P}\Big\{\frac{1}{2^n}\sum_{i=1}^{2^n}\mathcal{U}_{n_i}(x_1,x_2)>
r|x_1-x_2|^\eta,A_n^\varepsilon\Big\}
+\mbf{P}\Big\{\frac{2}{2^n}\sum_{i=1}^{2^n}\mathcal{U}_{n_i}(x_1,x_2)^->
r|x_1-x_2|^\eta,A_n^\varepsilon\Big\} \cr
 \ar\le\ar
(C_\varepsilon r^{-1}|x_1-x_2|)^{C'_\varepsilon
r|x_1-x_2|^{(\eta-\eta_c)/2}} + C_\varepsilon\exp\Big\{
-\frac{C'_\varepsilon
r^{\alpha/(\alpha-1)}}{|x_1-x_2|^{(2-\alpha-\eta\alpha)/(\alpha-1)}}\Big\}.
 \eeqnn
As the same argument we can also get the same estimation for
$\mathcal{V}_s(x_1,x_2)$. It then
follows from \eqref{3.26} that
 \beqnn
 \ar\ar
\mbf{P}\Big\{\frac{1}{2^n}\sum_{i=1}^{2^n}|Z_s(x_1)-Z_s(x_2)|>8r|x_1-x_2|^\eta,
A_n^\varepsilon\Big\}
\cr
 \ar\ar\quad\le
(C_\varepsilon r^{-1}|x_1-x_2|)^{C'_\varepsilon
r|x_1-x_2|^{(\eta-\eta_c)/2}} + C_\varepsilon\exp\Big\{
-\frac{C'_\varepsilon
r^{\alpha/(\alpha-1)}}{|x_1-x_2|^{(2-\alpha-\eta\alpha)/(\alpha-1)}}\Big\}.
 \eeqnn

Define $Z_t^{n,\varepsilon}=Z_t1_{A_n^\varepsilon}$. By the proof of
Lemma \ref{t1.7}, $Z_t$ has a continuous modification $\tilde{Z}_t$
for fixed $t>0$. Then
$\tilde{Z}_t^{n,\varepsilon}:=\tilde{Z}_t1_{A_n^\varepsilon}$ is a
continuous modification of $Z_t^{n,\varepsilon}$ for fixed $t>0$. By
\cite[Lemma III.5.1]{GS74}, it is easy to see that
 \beqlb\label{3.13}
\mbf{P}\Big\{ \frac{1}{2^n}\sum_{k=1}^{2^n}\sup_{x,z\in
\mbb{K},|x-z|\le\delta}
|\tilde{Z}_{n_k}^{n,\varepsilon}(x)-\tilde{Z}_{n_k}^{n,\varepsilon}(z)|
\ge
r\mathcal{E}\Big(\Big[\log_2\frac{\mathbb{K}_0}{2\delta}\Big]\Big)
\Big\} \le Q\Big(\Big[\log_2\frac{\mathbb{K}_0}{2\delta}\Big],r\Big)
 \eeqlb
for all $n\ge1$ and $\delta>0$, where
 \beqnn
\mathcal{E}(m):=\sum_{l=m}^\infty 8 (2^{-l}\mbb{K}_0)^\eta
=\frac{8\mbb{K}_0^\eta }{1-2^{-\eta}}2^{-\eta m}
 \eeqnn
and
 \beqnn
Q(m,r):=\sum_{l=m}^\infty2^{l+1}\Big[(C_\varepsilon
r^{-1}2^{-l}\mbb{K}_0)^{C'_\varepsilon
r(2^{-l}\mbb{K}_0)^{(\eta-\eta_c)/2}}
+C_\varepsilon\exp\Big(-\frac{C'_\varepsilon
r^{\alpha/(\alpha-1)}}{(2^{-l}\mbb{K}_0)^{(2-\alpha-\eta\alpha)/(\alpha-1)}}
\Big) \Big].
 \eeqnn
It is easy to check that  for
 \beqlb\label{3.14}
Q(r):=\sum_{m=0}^\infty Q(m,r)<\infty,\quad Q(r)\to 0 \mbox{\quad
as\quad} r\to\infty.
 \eeqlb
Observe that for each $m,n\ge1$ and $1\le k\le 2^n$
 \beqnn
\sup_{x,z\in \mbb{K},\frac{\delta}{2^{m+1}}<|x-z|\le
\frac{\delta}{2^m}} \frac{|\tilde{Z}_{n_k}^{n,\varepsilon}(x)
-\tilde{Z}_{n_k}^{n,\varepsilon}(z)|}{|x-z|^\eta} \le \sup_{x,z\in
\mbb{K},|x-z|\le \frac{\delta}{2^m}}
|\tilde{Z}_{n_k}^{n,\varepsilon}(x)
-\tilde{Z}_{n_k}^{n,\varepsilon}(z)|\Big(\frac{\delta}{2^{m+1}}\Big)^{-\eta}.
 \eeqnn
This implies
 \beqnn
 \ar\ar
\Big\{\frac{1}{2^n}\sum_{k=1}^{2^n} \sup_{x,z\in
\mbb{K},\frac{\delta}{2^{m+1}}<|x-z|\le \frac{\delta}{2^m}}
\frac{|\tilde{Z}_{n_k}^{n,\varepsilon}(x)-\tilde{Z}_{n_k}^{n,\varepsilon}(z)|}
{|x-z|^\eta}\ge  r\Big\} \cr
 \ar\ar\quad\subset
\Big\{\frac{1}{2^n}\sum_{k=1}^{2^n}\sup_{x,z\in
\mbb{K},|x-z|\le\frac{\delta}{2^m}}
|\tilde{Z}_{n_k}^{n,\varepsilon}(x)-\tilde{Z}_{n_k}^{n,\varepsilon}(z)|
\ge  r\Big(\frac{\delta}{2^{m+1}}\Big)^{\eta}\Big\} \cr
 \ar\ar\quad\subset
\Big\{ \frac{1}{2^n}\sum_{k=1}^{2^n}\sup_{x,z\in
\mbb{K},|x-z|\le\frac{\delta}{2^m}}
|\tilde{Z}_{n_k}^{n,\varepsilon}(x)-\tilde{Z}_{n_k}^{n,\varepsilon}(z)|
\ge
rc_1\mathcal{E}\Big(\Big[\log_2\frac{\mbb{K}_0}{2^{-m}\delta}\Big]\Big)
\Big\},
 \eeqnn
where
$c_1:=\inf_{m\ge0}\mathcal{E}\Big(\Big[\log_2\frac{\mbb{K}_0}{2^{-m}\delta}\Big]\Big)^{-1}
\Big(\frac{\delta}{2^{m+1}}\Big)^{\eta}>0$. It follows from
\eqref{3.13} that for each $m\ge0$,
 \beqnn
 \ar\ar
\mbf{P}\Big\{\frac{1}{2^n}\sum_{k=1}^{2^n} \sup_{x,z\in
\mbb{K},\frac{\delta}{2^{m+1}}<|x-z|\le \frac{\delta}{2^m}}
\frac{|\tilde{Z}_{n_k}^{n,\varepsilon}(x)-\tilde{Z}_{n_k}^{n,\varepsilon}(z)|}
{|x-z|^\eta}\ge r\Big\} \cr
 \ar\ar\qquad\le
\mbf{P}\Big\{\frac{1}{2^n}\sum_{k=1}^{2^n} \sup_{x,z\in
\mbb{K},|x-z|\le \frac{\delta}{2^m}}
|\tilde{Z}_{n_k}^{n,\varepsilon}(x)-\tilde{Z}_{n_k}^{n,\varepsilon}(z)|
\ge
rc_1\mathcal{E}\Big(\Big[\log_2\frac{\mbb{K}_0}{2^{-m}\delta}\Big]\Big)
\Big\} \cr
 \ar\ar\qquad\le
Q\Big(\Big[\log_2\frac{\mbb{K}_0}{2^{-m}\delta}\Big],rc_1\Big),
 \eeqnn
which implies
 \beqnn
\mbf{P}\Big\{ \frac{1}{2^n}\sum_{k=1}^{2^n} \sup_{x,z\in
\mbb{K},0<|x-z|\le \delta}
\frac{|\tilde{Z}_{n_k}^{n,\varepsilon}(x)-\tilde{Z}_{n_k}^{n,\varepsilon}(z)|}
{|x-z|^\eta}\ge r\Big\} \le Q(rc_1).
 \eeqnn
Then by Fatou's lemma and \eqref{3.11} for each subsequence
$\{n':n'\ge1\}$ of $\{n:n\ge1\}$,
 \beqnn
 \ar\ar
\mbf{P}\Big\{\liminf_{n'\to\infty}
\frac{1}{2^{n'}}\sum_{k=1}^{2^{n'}} \sup_{x,z\in \mbb{K}\cap
\mbb{Q},0<|x-z|\le\delta} \frac{|Z_{{n'}_k}(x)-Z_{{n'}_k}(z)|}
{|x-z|^\eta}\ge r\Big\} \cr
 \ar\ar\quad\le
\liminf_{n'\to\infty}\mbf{P}\Big\{
\frac{1}{2^{n'}}\sum_{k=1}^{2^{n'}} \sup_{x,z\in
\mbb{K}\cap\mbb{Q},0<|x-z|\le\delta}
\frac{|Z_{{n'}_k}(x)-Z_{{n'}_k}(z)|} {|x-z|^\eta}\ge r\Big\} \cr
 \ar\ar\quad\le
\liminf_{{n'}\to\infty}\Big\{\mbf{P}\Big(
\frac{1}{2^{n'}}\sum_{k=1}^{2^{n'}} \sup_{x,z\in
\mbb{K},0<|x-z|\le\delta}
\frac{|\tilde{Z}_{{n'}_k}^{{n'},\varepsilon}(x)-\tilde{Z}_{{n'}_k}^{{n'},\varepsilon}(z)|}
{|x-z|^\eta}\ge r\Big) +\mbf{P}(A_{n'}^{\varepsilon,c})\Big\} \cr
 \ar\ar\quad\le
Q(rc_1)+4\varepsilon.
 \eeqnn

First letting $r\to\infty$ and then letting $\varepsilon\to0$ we
immediately have
 \beqnn
\liminf_{n'\to\infty} \frac{1}{2^{n'}}\sum_{k=1}^{2^{n'}}
\sup_{x,z\in \mbb{K}\cap \mbb{Q},0<|x-z|\le\delta}
\frac{|Z_{{n'}_k}(x)-Z_{{n'}_k}(z)|} {|x-z|^\eta}<\infty,
\quad\mbf{P}\mbox{-a.s.}
 \eeqnn
 by \eqref{3.14}. Letting $\delta\to\infty$ we have
 \beqlb\label{3.20}
\liminf_{n'\to\infty} \frac{1}{2^{n'}}\sum_{k=1}^{2^{n'}}
\sup_{x,z\in \mbb{K}\cap \mbb{Q}}
\frac{|Z_{{n'}_k}(x)-Z_{{n'}_k}(z)|} {|x-z|^\eta}<\infty,
\qquad\mbf{P}\mbox{-a.s.}
 \eeqlb
Thus \eqref{2.17} follows from \eqref{1.3} and
\eqref{1.25}--\eqref{1.26}. This completes the proof.   \qed

\section{Proof of Theorem \ref{t4.1}}

\setcounter{equation}{0}

\subsection{The proof}

In this subsection we prove the pathwise uniqueness of solution for
\eqref{1.1}.
Throughout this subsection we always assume that the assumptions of Theorem \ref{t4.1} hold. The proof of Theorem \ref{t4.1} adopts the arguments
from \cite{MPS06, MyP11}. By
considering a conditional probability, we may assume that the
initial states $X_0$ and $Y_0$ are both deterministic. For $n\ge1$ define
 \beqnn
a_n:=\exp\{-n(n+1)/2\}.
 \eeqnn
Then $a_{n+1}=a_n a_n^{2/n}$. Let $\psi_n\in C_c^\infty(\mbb{R})$
satisfy $supp(\psi_n)\subset(a_n,a_{n-1})$,
$\int_{a_n}^{a_{n-1}}\psi_n(x)dx=1$, and $0\le\psi_n(x)\le 2/(nx)$
for all $x>0$ and $n\ge1$. For $x\in\mbb{R}$ and $n\ge1$ let
 \beqnn
\phi_n(x):=\int_0^{|x|}dy\int_0^y\psi_n(z)dz.
 \eeqnn
Then $\|\phi'_n\|\le 1$,  $\phi_n(x)\to |x|$, and $\phi_n'(x)\to x/|x|$ for $x\neq 0$
as $n\to\infty$.
For $n\ge1$ and $y,z\in\mbb{R}$ put
 \beqnn
D_n(y,z):=\phi_n(y+z)-\phi_n(y)-z\phi_n'(y)
\mbox{~~and~~}
H_n(y,z):=\phi_n(y+z)-\phi_n(y).
 \eeqnn

Let $\Phi\in C_c^\infty(\mbb{R})$ satisfy $0\le \Phi\le 1$,
$supp(\Phi)\subset(-1,1)$ and $\int_{\mbb{R}}\Phi(x)dx=1$. Let
$\Phi_x^m(y)=\Phi^m(x,y):=m\Phi(m(x-y))$ for $x,y\in\mbb{R}$ and
$m\ge1$. For $t\ge0$ and $y\in\mbb{R}$ let $U_t(y):=X_t(y)-Y_t(y)$,
$V_t(y):=H(X_t(y))-H(Y_t(y))$ and $R_t(y):=G(X_t(y))-G(Y_t(y))$. By
the argument in Section 2, both $\{X_t:t\ge0\}$ and $\{Y_t:t\ge0\}$
satisfy equation \eqref{1.1b}. Using \eqref{1.1b} and It\^o's
formula we have
 \beqlb\label{4.1}
\phi_n(\<U_t,\Phi_x^m\>)
 \ar=\ar
\frac12\int_0^t\phi_n'(\<U_s,\Phi_x^m\>)\<U_s,\Delta\Phi_x^m\>ds
+\int_0^t\phi_n'(\<U_s,\Phi_x^m\>)\<R_s,\Phi_x^m\>ds  \cr
 \ar\ar
+\int_0^t\int_0^\infty\int_{\mbb{R}}
H_n(\<U_{s-},\Phi_x^m\>,zV_{s-}(y)\Phi_x^m(y))\tilde{N}(ds,dz,dx) \cr
 \ar\ar
+\int_0^tds\int_0^\infty m_0(dz)\int_{\mbb{R}}
D_n(\<U_s,\Phi_x^m\>,zV_s(y)\Phi_x^m(y))dy.
 \eeqlb
For $t>0$ let $\tilde{X}_t$ and $\tilde{Y}_t$ denote the continuous
modifications of $X_t$ and $Y_t$, respectively. Let
$\tilde{U}_t(y):=\tilde{X}_t(y)-\tilde{Y}_t(y)$,
$\tilde{V}_t(y):=H(\tilde{X}_t(y))-H(\tilde{Y}_t(y))$ and
$\tilde{R}_t(y):=G(\tilde{X}_t(y))-G(\tilde{Y}_t(y))$.

Suppose that $T,K>0$ and that $\Psi$ is a nonnegative and compactly
supported infinitely differentiable function on $[0,T]\times\mbb{R}$
satisfying
 \beqnn
\Psi_s(x)=0\mbox{ for all }(s,x)\in[0,T]\times[-K,K]^c.
 \eeqnn
By \eqref{4.1} and a stochastic Fubini's theorem, it is easy to see that
 \beqnn
 \ar\ar
\<\phi_n(\<U_t,\Phi_{\cdot}^m\>),\Psi_t\> \cr
 \ar=\ar
\sum_{i=1}^k[\<\phi_n(\<U_{t_i},\Phi_{\cdot}^m\>),\Psi_{t_i}\>-\<
\phi_n(\<U_{t_{i-1}},\Phi_{\cdot}^m\>),\Psi_{t_{i-1}}\>]\cr
 \ar=\ar
\sum_{i=1}^k[\<\phi_n(\<U_{t_i},\Phi_{\cdot}^m\>),\Psi_{t_{i-1}}\>
-\<\phi_n(\<U_{t_{i-1}},\Phi_{\cdot}^m\>),\Psi_{t_{i-1}}\>] +
\sum_{i=1}^k[\<\phi_n(\<U_{t_i},\Phi_{\cdot}^m\>),\Psi_{t_i}-\Psi_{t_{i-1}}\>]
\cr
 \ar=\ar
\frac12\int_0^t\sum_{i=1}^kI_i(s)
\<\phi_n'(\<U_s,\Phi_{\cdot}^m\>)\<U_s,\Delta\Phi_{\cdot}^m\>,\Psi_{t_{i-1}}\>ds
\cr
 \ar\ar
+\int_0^t\sum_{i=1}^kI_i(s)
\<\phi_n'(\<U_s,\Phi_{\cdot}^m\>)\<R_s,\Phi_{\cdot}^m\>,\Psi_{t_{i-1}}\>ds
+
\int_0^t\sum_{i=1}^kI_i(s)\<\phi_n(\<U_{t_i},\Phi_{\cdot}^m\>),\dot{\Psi}_s\>ds
\cr
 \ar\ar
+\int_0^t\int_0^\infty\int_{\mbb{R}}\sum_{i=1}^kI_i(s)
\<H_n(\<U_{s-},\Phi_{\cdot}^m\>,zV_{s-}(y)\Phi_{\cdot}^m(y)),
\Psi_{t_{i-1}}\>\tilde{N}(ds,dz,dx) \cr
 \ar\ar
+\int_0^tds\int_0^\infty m_0(dz)\int_{\mbb{R}}
\sum_{i=1}^kI_i(s)\<D_n(\<U_s,\Phi_{\cdot}^m\>,zV_s(y)\Phi_{\cdot}^m(y)),\Psi_{t_{i-1}}\>dy,
 \eeqnn
where $0=t_0<t_1<\cdots<t_k=t$ and $I_i(s):=1_{(t_{i-1},t_i]}(s)$.
Letting $\max_{1\le i\le k}(t_i-t_{i-1})$ converge to zero we have
$\mbf{P}$-a.s.
 \beqlb\label{4.2}
\<\phi_n(\<U_t,\Phi_{\cdot}^m\>),\Psi_t\>
 \ar=\ar
\frac12\int_0^t\<\phi_n'(\<U_s,\Phi_{\cdot}^m\>)\<U_s,\Delta\Phi_{\cdot}^m\>,\Psi_s\>ds \cr
 \ar\ar
+\int_0^t\<\phi_n'(\<U_s,\Phi_{\cdot}^m\>)\<R_s,\Phi_{\cdot}^m\>,\Psi_s\>ds
+\int_0^t\<\phi_n(\<U_s,\Phi_{\cdot}^m\>),\dot{\Psi}_s\>ds \cr
 \ar\ar
+\int_0^t\int_0^\infty\int_{\mbb{R}}
\<H_n(\<U_{s-},\Phi_{\cdot}^m\>,zV_{s-}(y)\Phi_{\cdot}^m(y)),\Psi_s\>\tilde{N}(ds,dz,dy)\cr
 \ar\ar
+\int_0^tds\int_0^\infty m_0(dz)\int_{\mbb{R}}
\<D_n(\<U_s,\Phi_{\cdot}^m\>,zV_s(y)\Phi_{\cdot}^m(y)),\Psi_s\>dy
\cr
 \ar=:\ar
I_1^{m,n}(t)+I_2^{m,n}(t)+I_3^{m,n}(t)+I_4^{m,n}(t)+I_5^{m,n}(t),
\qquad t\ge0.
 \eeqlb

For $k\ge1$ define a stopping time
$\gamma_k$ by
 \begin{equation}\label{def_gamma}
\gamma_k:=\inf\Big\{t\in(0,T]: \<X_t,1\>+\<Y_t,1\>>k \Big\}
\end{equation}
with the convention $\inf \emptyset=\infty$. By Definition
\ref{t1.10}, $\{1_{\{t=0\}}X_0(1)+1_{\{t>0\}}\<X_t,1\>:t\ge0\}$ and
$\{1_{\{t=0\}}Y_0(1)+1_{\{t>0\}}\<Y_t,1\>:t\ge0\}$ are c\'adl\'ag
processes. Thus
 \beqnn
\sup_{t\in(0,T]}[\<X_t,1\>+\<Y_t,1\>]<\infty,\quad
\mbf{P}\mbox{-a.s.},
 \eeqnn
which implies $\lim_{k\to\infty}\gamma_k=\infty$ almost surely. Let
$\{l':l'\ge1\}$ be the subsequence of $\{l:l\ge1\}$ that will be
determined later in Lemma \ref{t4.6}. For $0\le i\le 2^{l'}$ let
${l'}_i=i/2^{l'}$. For any nonnegative function $f$ define
 \beqnn
\bar{\int}_{(0,t]}f(s)ds :=\liminf_{l'\to\infty}\int_0^t
\sum_{i=1}^{2^{l'}}1_{({l'}_{i-1}T,{l'}_iT]}(s) f({l'}_iT)ds,\qquad
t>0
 \eeqnn
and
 \beqnn
\bar{\int}_{(0,t)}f(s)ds :=\lim_{t'\uparrow
t}\bar{\int}_{(0,t']}f(s)ds.
 \eeqnn
Then $\bar{\int}_{(0,t]}f(s)ds\le \bar{\int}_{(0,T]}f(s)ds$. For
fixed $K>0$ and $0<\eta<\eta_c=\frac{2}{\alpha}-1$ define $\sigma_k$
by
 \begin{equation}\label{def_sigma}
\sigma_k=\inf\Big\{t\in (0,T] :\bar{\int}_{(0,t]}\qquad
\sup_{\makebox[0pt]{$%
{{\scriptscriptstyle x\neq z,\atop{%
{x,z\in [-(K+1),K+1]\cap \mbb{Q}\atop }}}}%
$}%
} \frac{|X_s(x)-X_s(z)|\vee|Y_s(x)-Y_s(z)|}{|x-z|^\eta}ds> k \Big\}.
\end{equation}
By Theorem \ref{t1.1},
 \beqnn
\bar{\int}_{(0,T]}\qquad
\sup_{\makebox[0pt]{$%
{{\scriptscriptstyle x\neq z,\atop{%
{x,z\in [-(K+1),K+1]\cap \mbb{Q}\atop }}}}%
$}%
} \frac{|X_s(x)-X_s(z)|\vee|Y_s(x)-Y_s(z)|}{|x-z|^\eta}ds <\infty,
\quad \mbf{P}\mbox{-a.s.},
 \eeqnn
which implies $\lim_{k\to\infty}\sigma_k=\infty,\,
\mbf{P}\mbox{-a.s.}$

In the rest of this subsection we always write
 \beqnn
\tau_k:=\min\{\gamma_k,\sigma_k\}.
 \eeqnn
Before proving Theorem \ref{t4.1}, we state three important lemmas.
Similar to \cite[Lemma 2.2(b)]{MPS06} we have the following result.
\blemma\label{t4.2} For any stopping time $\tau$ and $t>0$, we have
 \beqlb\label{4.8}
\limsup_{m,n\to\infty}\mbf{E}\{
I_1^{m,n}(t\wedge \tau)\}
\le\frac12\,\mbf{E}\Big\{\int_0^{t\wedge \tau}ds\int_{\mbb{R}}
|U_s(x)|\Delta\Psi_s(x)dx\Big\},
 \eeqlb
 \beqnn
\lim_{m,n\to\infty}\mbf{E}\{ I_2^{m,n}(t\wedge \tau)\} =
\mbf{E}\Big\{\int_0^{t\wedge \tau}ds\int_{\mbb{R}} \mathop{\rm
sgn}(U_s(x))R_s(x)\Psi_s(x)dx\Big\}
 \eeqnn
and
 \beqnn
\lim_{m,n\to\infty}\mbf{E}\{
I_3^{m,n}(t\wedge \tau)\}
=\mbf{E}\Big\{\int_0^{t\wedge \tau}ds\int_{\mbb{R}}
|U_s(x)|\dot{\Psi}_s(x)dx\Big\}.
 \eeqnn
\elemma

\blemma\label{t4.3} For any stopping time $\tau$, any $t>0$ and
$m,n\ge1$, we have
 \beqlb\label{4.3}
\mbf{E}\{I_4^{m,n}(t\wedge \tau)\}=0.
 \eeqlb
\elemma

The first inequality of \eqref{4.17} equals to
$\beta>\frac{(\alpha-1)(\eta_c+1)}{(2-\alpha)\eta_c}$,
which also equals to $\eta_c^{-1}<\frac{(2-\alpha)\beta}{\alpha-1}-1$.
Thus there exist constants
$\varepsilon,\delta>0$ satisfying
$\eta_c^{-1}<\delta<\frac{(2-\alpha)\beta}{\alpha-1}-1$
and $\frac{\delta+1}{2-\alpha}<\varepsilon<\frac{\beta}{\alpha-1}$,
which implies
 \beqlb\label{4.25}
\frac{\delta+1}{2-\alpha}<\varepsilon<\frac{\delta\eta_c\beta}{\alpha-1}
\mbox{~~and~~}\varepsilon<\frac{\beta}{\alpha-1}.
 \eeqlb

\blemma\label{t4.4} If $m=a_{n-1}^{-\delta}$ for the $\delta$ in (\ref{4.25}), then for each $t>0$
and $k\ge1$,
 \beqnn
\lim_{n\to\infty}\mbf{E}\{I_5^{m,n}(t\wedge \tau_{k})\}=0.
 \eeqnn
\elemma

Deferring the proofs of Lemmas 4.3--4.5, we first present the main
proof.

\noindent{\it Proof of Theorem \ref{t4.1}.}
By the continuity of $x\mapsto \tilde{U}_t(x)$,
for each $x\in\mbb{R}$ and $t>0$,
 \beqlb\label{4.4}
\lim_{m\to\infty}\<\tilde{U}_t,\Phi_x^m\>
=\lim_{m\to\infty}\int_{\mbb{R}}\tilde{U}_t(x-\frac{y}{m})\Phi(y)dy
=\tilde{U}_t(x).
 \eeqlb
Note that $\|\phi'_n\|\le1$. Then for all $x_m\to x$ as $m\to\infty$, we have
that
 \beqnn
|\phi_n(x_m)-|x|| \le |\phi_n(x_m)-\phi_n(x)| +|\phi_n(x)-|x|| \le
|x_m-x|+|\phi_n(x)-|x||,
 \eeqnn
which converges to zero as $m,n\to\infty$. Now by \eqref{4.4} and
Fatou's lemma
 \beqnn
 \ar\ar
\mbf{E}\{\<|U_t|,\Psi_t\>1_{\{t\le\tau_{k}\}}\}
=
\mbf{E}\Big\{\<|\tilde{U}_t|,\Psi_t\>1_{\{t\le\tau_{k}\}}\Big\} \cr
 \ar=\ar
\mbf{E}\Big\{
\<\lim_{m,n\to\infty}\phi_n(\big\<\tilde{U}_t,
\Phi_{\cdot}^m\>),\Psi_t\big\>1_{\{t\le\tau_{k}\}}\Big\}
\le
\liminf_{m,n\to\infty}\mbf{E}\Big\{
\<\phi_n(\big\<\tilde{U}_t,
\Phi_{\cdot}^m\>),\Psi_t\big\>1_{\{t\le\tau_{k}\}}\Big\} \cr
 \ar=\ar
\liminf_{m,n\to\infty}\mbf{E}\Big\{
\<\phi_n(\big\<U_t,
\Phi_{\cdot}^m\>),\Psi_t\big\>1_{\{t\le\tau_{k}\}}\Big\}
\le
\liminf_{m,n\to\infty}\mbf{E}\Big\{
\<\phi_n(\big\<U_{t\wedge\tau_{k}},
\Phi_{\cdot}^m\>),\Psi_{t\wedge\tau_{k}}\big\>\Big\}.
 \eeqnn
Together with \eqref{4.2} and Lemmas \ref{t4.2}--\ref{t4.4} we have
 \beqnn
\mbf{E}\{\<|U_t|,\Psi_t\>1_{\{t\le\tau_{k}\}}\}
 \ar\le\ar
\mbf{E}\Big\{\int_0^{t\wedge\tau_{k}}ds
\int_{\mbb{R}}|U_s(x)|\big[\frac12\Delta\Psi_s(x)+\dot{\Psi}_s(x)\big]dx\Big\}
\cr
 \ar\ar
+\mbf{E}\Big\{\int_0^{t\wedge \tau_k}ds\int_{\mbb{R}} \mathop{\rm
sgn}(U_s(x))R_s(x)\Psi_s(x)dx\Big\}.
 \eeqnn
Letting $k\to\infty$ in the above inequality
we have
 \beqnn
\mbf{E}\{\<|U_t|,\Psi_t\>\}
 \ar\le\ar
\int_0^tds\int_{\mbb{R}}
\mbf{E}[|U_s(x)|]\big[\frac12\Delta\Psi_s(x)+\dot{\Psi}_s(x)\big]dx
\cr
 \ar\ar
+\int_0^tds\int_{\mbb{R}}\mbf{E}\big[\mathop{\rm
sgn}(U_s(x))R_s(x)\big]\Psi_s(x)dx.
 \eeqnn
This is similar to (34) in \cite{MPS06}. Then by the same argument
as in Theorem 1.6 of \cite{MPS06}, for any fixed $t>0$
and nonnegative $f\in C_c^\infty(\mbb{R})$, with $\Psi_s(x)$
replaced by $\bar{\Psi}_N(s,x):=(P_{t-s}f(x))g_N(x)$ for a proper
sequence of functions $(g_N)_{N\ge1}$ so that $g_N(x)\to 1$ for all
$x\in\mbb{R}$ and the first term on the right hand side of the above inequality tends to
zero as $N\to\infty$. Thus, we have
 \beqnn
\<\mbf{E}[|U_t|],f\> \le
\int_0^tds\int_{\mbb{R}}\mbf{E}\big[\mathop{\rm
sgn}(U_s(x))R_s(x)\big]P_{t-s}f(x)dx \le
\int_0^t\<\mbf{E}[r_0(|U_s|)],P_{t-s}f\>ds,
 \eeqnn
where condition (C2) was used in the last inequality. It is
elementary to check that the above inequality holds for each $f\in
\mathscr{B}(\mbb{R})^+$ satisfying $\lambda_0(f)<\infty$. This means
that for each $f\in \mathscr{B}(\mbb{R})^+$ satisfying
$\lambda_0(f)=1$,
 \beqnn
\<\mbf{E}[|U_t|],P_{T-t}f\> \le
\int_0^t\<\mbf{E}[r_0(|U_s|)],P_{T-s}f\>ds.
 \eeqnn
Then by the concaveness of $x\mapsto r_0(x)$ and Jensen's inequality,
 \beqlb\label{1.31}
 \ar\ar
\<\mbf{E}[|U_t|],P_{T-t}f\>
 \leq
\int_{\mbb{R}}\mbf{E}[|U_t(x)|]P_{T-t}f(x)dx \cr
 \ar\ar\quad=
\int_0^tds\int_{\mbb{R}}r_0\big(\mbf{E}[|U_s(x)|]\big)P_{T-s}f(x)dx
 \le
\int_0^tr_0\big(\<\mbf{E}[|U_s|],P_{T-s}f\>\big)ds.
 \eeqlb
Since $\int_{0+}r_0(z)^{-1}dz=\infty$, the above inequality implies that
$\<\mbf{E}[|U_t|],P_{T-t}f\>=0$ for all $t>0$. Thus
 \beqnn
\mbf{P}\{X_t(x)=Y_t(x)\mbox{~for~}\lambda_0\mbox{-a.e. }x\}=1
 \eeqnn
for all $t>0$. It follows that $\<X_t,f\> = \<Y_t,f\>$
$\mbf{P}$-a.s. for all $t> 0$ and $f\in \mcr{S}(\mbb{R})$. By the
right-continuities of $t\mapsto \<X_t,f\>$ and $t\mapsto \<Y_t,f\>$
we have $\mbf{P} \{\<X_t,f\> = \<Y_t,f\>$ for all $t>0\} = 1$ for
all $f\in \mcr{S}(\mbb{R})$. Considering a suitable sequence
$\{f_1,f_2, \cdots\}\subset \mcr{S}(\mbb{R})$ we can conclude
\eqref{4.16}. \qed

We now present the proofs of Lemmas \ref{t4.2}--\ref{t4.4}.

\noindent{\it Proof of Lemma \ref{t4.2}.}
By the same argument as Lemma 2.2(b) of \cite{MPS06},
 \beqlb\label{4.5}
\limsup_{m,n\to\infty}\mbf{E}\{
I_1^{m,n}(t\wedge \tau)\}
\le\limsup_{m,n\to\infty}\frac12
\,\mbf{E}\Big\{\int_0^{t\wedge \tau}ds\int_{\mbb{R}}
\phi'_n(\<U_s,\Phi_x^m\>)\<U_s,\Phi_x^m\>\Delta\Psi_s(x)dx\Big\}.
 \eeqlb
It follows from the dominated convergence and the continuity of
$x\mapsto \tilde{U}_s(x)$ that $\int_{-1}^1
\tilde{U}_s(x-\frac{y}{m})\Phi(y)dy\to \tilde{U}_s(x)$ as
$m\to\infty$. Now using the fact that $\phi'_n(x)\to \mathop{\rm
sgn}(x)$ as $n\to\infty$, we have
 \beqnn
\lim_{m,n\to\infty}\phi'_n(\<\tilde{U}_s,\Phi_x^m\>)\<\tilde{U}_s,\Phi_x^m\>
=\lim_{m,n\to\infty}\phi'_n(\<\tilde{U}_s,\Phi_x^m\>)
\int_{-1}^1\tilde{U}_s(x-\frac{y}{m})\Phi(y)dy
=|\tilde{U}_s(x)|.
 \eeqnn
Observe that $\|\phi'_n\|\le1$ for all $n\ge1$ and
 \beqnn
0\le\phi'_n(\<\tilde{U}_s,\Phi_x^m\>)\<\tilde{U}_s,\Phi_x^m\>
=\phi'_n(\<\tilde{U}_s,\Phi_x^m\>)\int_{\mbb{R}}
\tilde{U}_s(x-\frac{y}{m})\Phi(y)dy
\le \sup_{|y|\le K+1}|\tilde{X}_s(y)+\tilde{Y}_s(y)|.
 \eeqnn
Then by \eqref{4.5}, \eqref{1.6} and the dominated convergence
 \beqnn
\limsup_{m,n\to\infty}\mbf{E}\{
I_1^{m,n}(t\wedge \tau)\}
 \ar\le\ar
\limsup_{m,n\to\infty}\frac12\int_0^tds\int_{\mbb{R}}
\mbf{E}\Big\{\phi'_n(\<\tilde{U}_s,\Phi_x^m\>)
\<\tilde{U}_s,\Phi_x^m\>\Delta\Psi_s(x)1_{\{s\le\tau\}}\Big\}dx \cr
 \ar\le\ar
\frac12\int_0^tds\int_{\mbb{R}}
\mbf{E}\Big\{\limsup_{m,n\to\infty}
\phi'_n(\<\tilde{U}_s,\Phi_x^m\>)\<\tilde{U}_s,\Phi_x^m\>
\Delta\Psi_s(x)1_{\{s\le\tau\}}\Big\}dx \cr
 \ar=\ar
\frac12\int_0^tds\int_{\mbb{R}}
\mbf{E}\Big\{|\tilde{U}_s(x)|\Delta\Psi_s(x)1_{\{s\le\tau\}}\Big\}dx  \cr
 \ar=\ar
\frac12\,\mbf{E}\Big\{\int_0^{t\wedge \tau}ds\int_{\mbb{R}}
|U_s(x)|\Delta\Psi_s(x)dx\Big\} .
 \eeqnn

By the continuity of $x\mapsto \tilde{U}_s(x)$ and $x\mapsto
\tilde{R}_s(x)$, one also sees that
 \beqnn
 \ar\ar
\lim_{m,n\to\infty}\phi'_n(\<\tilde{U}_s,\Phi_x^m\>)\<\tilde{R}_s,\Phi_x^m\>
\cr
 \ar\ar\qquad=
\lim_{m,n\to\infty}\phi'_n(\<\tilde{U}_s,\Phi_x^m\>)
\int_{-1}^1\tilde{R}_s(x-\frac{y}{m})\Phi(y)dy =\mathop{\rm
sgn}(\tilde{U}_s(x))\tilde{R}_s(x).
 \eeqnn
By condition (C1) and the fact $\|\phi'_n\|\le 1$ one sees that
 \beqnn
 \ar\ar
|\phi'_n(\<\tilde{U}_s,\Phi_x^m\>)\<\tilde{R}_s,\Phi_x^m\>| \le
|\<\tilde{R}_s,\Phi_x^m\>|=\Big|\int_{\mbb{R}}\tilde{R}_s(x-\frac{y}{m})\Phi(y)dy\Big|
\cr
 \ar\ar\qquad\le
C\int_{\mbb{R}}[\tilde{X}_s(x-\frac{y}{m})+\tilde{Y}_s(x-\frac{y}{m})]\Phi(y)dy+C
\le C\sup_{|y|\le K+1}|\tilde{X}_s(y)+\tilde{Y}_s(y)|+C.
 \eeqnn
By the dominated convergence again that
 \beqnn
\lim_{m,n\to\infty}\mbf{E}\{
I_2^{m,n}(t\wedge \tau)\}
 \ar=\ar
\int_0^{t\wedge\tau}ds\int_{\mbb{R}}
\mbf{E}\Big\{\lim_{m,n\to\infty}\phi'_n(\<\tilde{U}_s,\Phi_x^m\>)
\<\tilde{R}_s,\Phi_x^m\>\Psi_s(x)\Big\}dx \cr
 \ar=\ar
\int_0^{t\wedge\tau}ds\int_{\mbb{R}} \mbf{E}\big\{\mathop{\rm
sgn}(U_s(x))R_s(x)\Psi_s(x)\big\}dx.
 \eeqnn

By the fact $\|\phi'_n\|\le 1$ again,
 \beqnn
 \ar\ar
\Big|\mbf{E}\{
I_3^{m,n}(t\wedge \tau)\}
-\mbf{E}\Big\{\int_0^{t\wedge \tau}ds\int_{\mbb{R}}
|U_s(x)|\dot{\Psi}_s(x)dx\Big\}\Big| \cr
 \ar\ar\qquad\le
\int_0^tds\int_{\mbb{R}}
\mbf{E}\{|\phi_n(\<U_s,\Phi_x^m\>)-|U_s(x)||\}|\dot{\Psi}_s(x)|dx \cr
 \ar\ar\qquad\le
\int_0^tds\int_{\mbb{R}}
\mbf{E}\Big\{
|\phi_n(\<U_s,\Phi_x^m\>)-\phi_n(U_s(x))|
+|\phi_n(U_s(x))-|U_s(x)||\Big\}|\dot{\Psi}_s(x)|dx \cr
 \ar\ar\qquad\le
\int_0^tds\int_{\mbb{R}}|\dot{\Psi}_s(x)|dx
\int_{\mbb{R}}\mbf{E}\Big\{|U_s(x-\frac{y}{m})-U_s(x)|\Big\}\Phi(y)dy \cr
 \ar\ar\qquad~~
+\int_0^tds\int_{\mbb{R}}
\mbf{E}\{|\phi_n(U_s(x))-|U_s(x)||\}|\dot{\Psi}_s(x)|dx.
 \eeqnn
Now by \eqref{1.5} and the dominated convergence one finishes the
proof. \qed

\noindent{\it Proof of Lemma \ref{t4.3}.}
For $t\ge0$ and $m,n\ge1$ let
 \beqnn
I_{4,1}^{m,n}(t) := \int_0^t\int_0^1\int_{\mbb{R}}
\<H_n(\<U_{s-},\Phi_{\cdot}^m\>,zV_{s-}(y)\Phi_{\cdot}^m(y)),\Psi_s\>\tilde{N}(ds,dz,dy)
\eeqnn and $I_{4,2}^{m,n}(t):=I_4^{m,n}(t)-I_{4,1}^{m,n}(t)$. By the
Burkholder-Davis-Gundy inequality (see p.195 of \cite{P05}), for
$\bar{\alpha}\in(\alpha,\frac\alpha\beta\wedge2)$ and $T>0$,
 \beqnn
 \ar\ar
\mbf{E}\Big\{\sup_{t\in[0,T]}|I_{4,1}^{m,n}(t\wedge \tau)|^{\bar{\alpha}}\Big\} \cr
 \ar\ar\qquad\le
C\,\mbf{E}\Big\{\Big[\int_0^T\int_0^1\int_{\mbb{R}}
|\<H_n(\<U_{s-},\Phi_{\cdot}^m\>,zV_{s-}(y)\Phi_{\cdot}^m(y)),\Psi_s\>|^2
N(ds,dz,dy)\Big]^{\frac{\bar{\alpha}}{2}}\Big\} \cr
 \ar\ar\qquad\le
C\,\mbf{E}\Big\{\int_0^T\int_0^1\int_{\mbb{R}}
|\<H_n(\<U_{s-},\Phi_{\cdot}^m\>,zV_{s-}(y)\Phi_{\cdot}^m(y)),\Psi_s\>|^{\bar{\alpha}}
N(ds,dz,dy)\Big\},
 \eeqnn
where for the last inequality we used the fact that
 \beqnn
\big|\sum_{i=1}^nx_i^2\big|^{\frac{\bar{\alpha}}{2}}\le
\sum_{i=1}^n|x_i|^{\bar{\alpha}}
 \eeqnn
for all $x_i\in\mbb{R}$ and $n\ge1$. Since $|H_n(y,z)|\le |z|$ for
all $y,z\in\mbb{R}$, and $\Psi$ is continuous and compactly
supported, then by the H\"older inequality, condition (C3) and Lemma
\ref{t1.4},
 \beqnn
 \ar\ar
\mbf{E}\Big\{\sup_{t\in[0,T]}|I_{4,1}^{m,n}(t\wedge \tau)|^{\bar{\alpha}}\Big\} \cr
 \ar\ar\qquad\le
C\,\mbf{E}\Big\{\int_0^T\int_0^1\int_{\mbb{R}}
\Big|z\int_{\mbb{R}}V_{s-}(y)\Phi_x^m(y)\Psi_s(x)dx\Big|^{\bar{\alpha}}
N(ds,dz,dy)\Big\}  \cr
 \ar\ar\qquad\le
C\,\mbf{E}\Big\{\int_0^T\int_0^1\int_{\mbb{R}}
z^{\bar{\alpha}}\Big[\int_{\mbb{R}}|V_{s-}(y)\Phi_x^m(y)\Psi_s(x)|^{\bar{\alpha}}dx\Big]
N(ds,dz,dy)\Big\} \cr
 \ar\ar\qquad=
C\int_0^Tds\int_0^1z^{\bar{\alpha}}m_0(dz)\int_{\mbb{R}}\Psi_s(x)^{\bar{\alpha}}dx
\int_{\mbb{R}}\mbf{E}\{|V_s(y)|^{\bar{\alpha}}\}\Phi_x^m(y)^{\bar{\alpha}}
dy \cr
 \ar\ar\qquad\le
C\int_0^Tds\int_{\mbb{R}}\Psi_s(x)^{\bar{\alpha}}dx
\int_{\mbb{R}}\mbf{E}\{|X_s(y)
-Y_s(y)|^{\bar{\alpha}\beta}\}\Phi_x^m(x)^{\bar{\alpha}}dy \cr
 \ar\ar\qquad\le
C\int_0^Tds\int_{\mbb{R}}\Psi_s(x)^{\bar{\alpha}}dx
\int_{\mbb{R}}\mbf{E}\{X_s(y)^{\bar{\alpha}\beta}
+Y_s(y)^{\bar{\alpha}\beta}\}\Phi_x^m(x)^{\bar{\alpha}}dy
<\infty.
 \eeqnn
Similarly,
 \beqnn
 \ar\ar
\mbf{E}\Big\{\sup_{t\in[0,T]}|I_{4,2}^{m,n}(t\wedge \tau)|\Big\} \cr
 \ar\ar\qquad\le
\mbf{E}\Big\{\int_0^T\int_1^\infty\int_{\mbb{R}}
|\<H_n(\<U_{s-},\Phi_{\cdot}^m\>,zV_{s-}(y)\Phi_{\cdot}^m(y)),\Psi_s\>|
N(ds,dz,dy)\Big\} \cr
 \ar\ar\qquad~~
+ \mbf{E}\Big\{\int_0^Tds\int_1^\infty m_0(dz)\int_{\mbb{R}}
|\<H_n(\<U_{s-},\Phi_{\cdot}^m\>,zV_{s-}(y)\Phi_{\cdot}^m(y)),\Psi_s\>|dy\Big\}
\cr
 \ar\ar\qquad=
2\int_0^Tds\int_1^\infty m_0(dz)\int_{\mbb{R}}
\mbf{E}\Big\{|\<H_n(\<U_{s-},\Phi_{\cdot}^m\>,zV_{s-}(y)\Phi_{\cdot}^m(y)),
\Psi_s\>|\Big\}dy \cr
 \ar\ar\qquad\le
2\int_0^Tds\int_1^\infty z m_0(dz)\int_{\mbb{R}}\Psi_s(x)dx
\int_{\mbb{R}}\mbf{E}\{|V_s(y)|\}|\Phi_x^m(y)|dy <\infty.
 \eeqnn
It follows that for $T>0$,
 \beqnn
\mbf{E}\Big\{\sup_{t\in[0,T]}|I_4^{m,n}(t\wedge \tau)|\Big\}<\infty.
 \eeqnn
Then by \cite[p.38]{P05}, $t\mapsto  I_4^{m,n}(t\wedge \tau)$
is a martingale, which implies \eqref{4.3}. \qed

To prove Lemma \ref{t4.4}, we only need to show the following two
lemmas. \blemma\label{t4.9} For $m,n,k,i\ge1$ and $t\in[0,T]$ let
 \beqnn
I_{5,1}^{m,n,k,i}(t):=\mbf{E}\Big\{\int_0^{t\wedge \gamma_k}ds
\int_0^{1/i} m_0(dz)\int_{\mbb{R}}
\<D_n(\<U_s,\Phi_{\cdot}^m\>,zV_s(y)\Phi_{\cdot}^m(y)),\Psi_s\>dy\Big\}.
 \eeqnn
Then
 \beqnn
I_{5,1}^{m,n,k,i}(t) \le C_Tkm(na_n)^{-1}i^{\alpha-2},\qquad
t\in[0,T].
 \eeqnn
\elemma

\blemma\label{t4.8} For $m,n,k,i\ge1$ and $t\in[0,T]$ let
 \beqnn
I_{5,2}^{m,n,k,i}(t):=\mbf{E}\Big\{\int_0^{t\wedge
\sigma_k}ds\int_{1/i}^\infty m_0(dz)\int_{\mbb{R}}
\<D_n(\<U_s,\Phi_{\cdot}^m\>,zV_s(y)\Phi_{\cdot}^m(y)),\Psi_s\>dy\Big\}.
 \eeqnn
Then
 \beqnn
I_{5,2}^{m,n,k,i}(t)\le C_T[k^\beta
m^{-\eta\beta}+a_{n-1}^\beta]i^{\alpha-1},\qquad t\in[0,T]. \eeqnn
\elemma

\noindent{\it Proof of Lemma \ref{t4.4}.} Recall that
$\delta,\varepsilon$ satisfy \eqref{4.25}.  We take
$i=a_{n-1}^{-\varepsilon}$ and $\eta<\eta_c$ satisfying \eqref{4.25}
with $\eta$ replaced by $\eta_c$. Then $\mbf{E}\{I_5^{m,n}(t\wedge
\tau_{k})\}\le I_{5,1}^{m,n,k,i}(t)+I_{5,2}^{m,n,k,i}(t)$ converges
to zero as $n\to\infty$ by Lemmas \ref{t4.9} and \ref{t4.8}.  \qed

We first present the proof for Lemma \ref{t4.9}.

\noindent{\it Proof of Lemma \ref{t4.9}.} Recall that $\psi_n(x)\le
2(na_n)^{-1}$. Then by (3.3) in \cite{LP12} and condition (C3),
 \beqnn
 \ar\ar
m^{-1}D_n(\<U_s,\Phi_x^m\>, zmV_s(x-y/m)\Phi(y)) \cr
 \ar=\ar
mz^2 V_s(x-y/m)^2\Phi(y)^2 \int_0^1\psi_n\big(|\<U_s,\Phi_x^m\>+zhm
V_s(x-y/m)\Phi(y)|\big) (1-h)dh \cr
 \ar\le\ar
Cm(na_n)^{-1}z^2|U_s(x-y/m)|^{2\beta}\Phi(y) \cr
 \ar\le\ar
Cm(na_n)^{-1}z^2|X_s(x-y/m)+Y_s(x-y/m)|^{((2\beta)\vee1)-1} \cr
 \ar\ar\quad
\times|1+X_s(x-y/m)+Y_s(x-y/m)|\Phi(y).
 \eeqnn
It follows that $\mbf{P}$-a.s.
 \beqlb\label{4.22}
 \ar\ar
\int_{-K}^K\Psi_s(x)dx\int_{-1}^1m^{-1}D_n(\<U_s,\Phi_x^m\>,
zmV_s(x-y/m)\Phi(y))dy \cr
 \ar\ar\quad\le
\frac{Cc_2mz^2\tilde{K}_s}{na_n}
\int_{-K}^Kdx\int_{-1}^1|1+X_s(x-y/m)+Y_s(x-y/m)|\Phi(y)dy \cr
 \ar\ar\quad\le
\frac{Cc_2mz^2\tilde{K}_s}{na_n} [1+\<X_s,1\>+\<Y_s,1\>]\le
\frac{Cc_2kmz^2\tilde{K}_s}{na_n}
 \eeqlb
on $\{s\le \gamma_k\}$, where
 \beqnn
c_2:=\sup_{(s,x)\in[0,T]\times[-K,K]}\Psi_s(x),\quad
\tilde{K}_s:=\sup_{|u|\le
K+1}|\tilde{X}_s(u)+\tilde{Y}_s(u)|^{((2\beta)\vee1)-1}.
 \eeqnn
Since $0<\beta<1$, $2\beta-1<\alpha_{-}$ for each
$\alpha_{-}\in(1,\alpha)$. Then by \eqref{1.6}, for each $0<s\le T$,
 \beqnn
\mbf{E}\{\tilde{K}_s\}\le 2\mbf{E}\Big\{1+ \sup_{|u|\le
K+1}|\tilde{X}_s(u)\vee\tilde{Y}_s(u)|^{\alpha_{-}}\Big\} \le
C_Ts^{-\frac{\alpha_{-}}{2}}.
 \eeqnn
The above inequality together with \eqref{4.22} leads to
 \beqnn
I_{5,1}^{m,n,k,i}(t) \le Ckm(na_n)^{-1}\int_0^{1/i}z^2m_0(dz) \le
Ckm(na_n)^{-1}i^{\alpha-2},
 \eeqnn
which finishes the proof. \qed

For $m,n,k,i\ge1$ and $t>0$ define
 \beqnn
J_{m,n,k,i}(t)
 \ar:=\ar
\int_{\mbb{R}}\Psi_t(x)dx\int_{1/i}^\infty zm_0(dz)
\int_{\mbb{R}}\Phi(y)dy \cr
 \ar\ar\quad\times
\int_0^1\tilde{D}_n\big(\<U_t,\Phi_x^m\>,mzhV_t(x-\frac{y}{m})\big)
V_t(x-\frac{y}{m})dh,
 \eeqnn
where $\tilde{D}_n(y,z)=\phi_n'(y+z)-\phi_n'(y)$ for all
$y,z\in\mbb{R}$ and $n\ge1$. Before proving Lemma \ref{t4.8} we need
to show two more lemmas. \blemma\label{t4.6} There is a subsequence
$\{l':l'\ge1\}$ of $\{l:l\ge1\}$ so that for each $m,n,k,i\ge1$,
$\mbf{P}$-a.s.
 \beqlb\label{4.13}
\lim_{l'\to\infty} \int_0^t
\sum_{j=1}^{2^{l'}}1_{({l'}_{j-1}T,{l'}_jT]}(s)
|J_{m,n,k,i}({l'}_jT)- J_{m,n,k,i}(s)|ds=0, \qquad t\in(0,T].
 \eeqlb
 \elemma \proof Observe that $\|\phi_n'\|\le1$.
Then by condition (C3)
 \beqnn
|M_{m,n,k}(x,y,z,h,t)|:=|\tilde{D}_n\big(\<U_t,\Phi_x^m\>,mzhV_t(x-\frac{y}{m})\big)
V_t(x-\frac{y}{m})| \le C|U_t(x-\frac{y}{m})|^\beta
 \eeqnn
and
 \beqnn
J_{m,n,k,i}(t)
 \le
c_2Ci^{\alpha-1}\int_{-K}^Kdx \int_{-1}^1|U_t(x-\frac{y}{m})|^\beta
dy,
 \eeqnn
where $c_2=\sup_{(s,x)\in[0,T]\times[-K,K]}\Psi_s(x)$.

Now by Lemma \ref{t1.4}, there is a constant $\delta\in(1,\alpha)$
so that for each $t\in(0,T]$, there is a set $K_t\subset\mbb{R}$ of
Lebesgue measure zero satisfying
 \beqlb\label{4.14}
\mbf{E}\big\{M_{m,n,k,i}(x,y,z,h,t)^\delta\big\} \le C_T
t^{-\beta\delta/2},\qquad x,y\in\mbb{R} \backslash K_t,
 \eeqlb
and
 \beqlb\label{4.20}
\mbf{E}\{J_{m,n,k,i}(t)^\delta\} \le C_Ti^{\delta(\alpha-1)}
t^{-\beta\delta/2}.
 \eeqlb
By Lemma \ref{t1.8} for each $t>0$ and $t_j\to t$ as $j\to\infty$,
there is a set $\bar{K}_t\subset\mbb{R}$ of Lebesgue measure zero so
that for each $ x\in\mbb{R} \backslash \bar{K}_t$, both
 \beqnn
X_{t_j}(x)\to X_t(x) \mbox{\quad and \quad} Y_{t_j}(x)\to Y_t(x)
 \eeqnn
in probability as $j\to\infty$. Then for each $x,y\in\mbb{R}
\backslash
 \bar{K}_t$,
 \beqnn
M_{m,n,k,i}(x,y,z,h,t_j) \rightarrow M_{m,n,k,i}(x,y,z,h,t)
 \eeqnn
in probability as $j\to\infty$. Together with \eqref{4.14} we have
 \beqnn
\lim_{j\to\infty}\mbf{E}\Big\{\big|M_{m,n,k,i}(x,y,z,h,t_j)
-M_{m,n,k,i}(x,y,z,h,t)\big|\Big\}=0, \qquad x,y\in\mbb{R}
\backslash (K_t\cup \bar{K}_t).
 \eeqnn
Using \eqref{4.20} and \cite[Theorem 4.5.2]{Chun} we then have
 \beqlb\label{4.21}
\lim_{t\downarrow
u}\mbf{E}\Big\{|J_{m,n,k,i}(t)-J_{m,n,k,i}(u)|\Big\}=0.
 \eeqlb
Using \eqref{4.20} again,
 \beqnn
 \ar\ar
\sup_{l\ge1}\int_0^T \Big|\sum_{j=1}^{2^l}1_{(l_{j-1}T,l_jT]}(s)
\mbf{E}\{|J_{m,n,k,i}(l_jT)-J_{m,n,k,i}(s)|\}\Big|^\delta ds \cr
 \ar\ar\qquad\le
\sup_{l\ge1}\int_0^T \sum_{j=1}^{2^l}1_{(l_{j-1}T,l_jT]}(s)
\mbf{E}\{J_{m,n,k,i}(l_jT)^\delta+J_{m,n,k,i}(s)^\delta\} ds<\infty.
 \eeqnn
It then follows from \eqref{4.21} and \cite[Theorem 4.5.2]{Chun} again that
 \beqnn
 \ar\ar
\lim_{l\to\infty}\mbf{E}\Big\{\int_0^T
\sum_{j=1}^{2^l}1_{(l_{j-1}T,l_jT]}(s)
|J_{m,n,k,i}(l_jT)-J_{m,n,k,i}(s)|ds\Big\} \cr
 \ar\ar\qquad=
\lim_{l\to\infty}\int_0^T \sum_{j=1}^{2^l}1_{(l_{j-1}T,l_jT]}(s)
\mbf{E}\{|J_{m,n,k,i}(l_jT)-J_{m,n,k,i}(s)|\}ds=0,
 \eeqnn
which implies that
 \beqnn
\int_0^T \sum_{i=1}^{2^l}1_{(l_{j-1}T,l_jT]}(s)
|J_{m,n,k,i}(l_jT)-J_{m,n,k,i}(s)|ds
 \eeqnn
goes to $0$ in probability as $l\to\infty$. Then for each
$m,n,k,i\ge1$, there is a subsequence $\{l':=l'(m,n,k,i):l'\ge1\}$
of $\{l:l\ge1\}$ so that \eqref{4.13} holds. Then one can choose a
proper subsequence of $\{l:l\ge1\}$ which is independent of
$m,n,k,i$ so that \eqref{4.13} holds. \qed

\blemma\label{t4.7} Let $\{l':l'\ge1\}$ be the subsequence of
$\{l:l\ge1\}$  in Lemma \ref{t4.6}. Then $\mbf{P}$-a.s.
 \beqnn
\lim_{l'\to\infty}\int_0^t
\sum_{j=1}^{2^{l'}}1_{({l'}_{j-1}T,{l'}_jT]}(s)
J_{m,n,k,i}({l'}_jT)1_{\{{l'}_jT<\sigma_k\}} ds =\int_0^t
J_{m,n,k,i}(s)1_{\{s<\sigma_k\}} ds, \quad t\in(0,T].
 \eeqnn
  \elemma

\proof
It follows from \eqref{4.20} that
 \beqnn
\mbf{E}\Big\{\int_0^TJ_{m,n,k,i}(s)ds\Big\}<\infty,
 \eeqnn
which implies that $\mbf{P}$-a.s.
 \beqnn
\int_0^TJ_{m,n,k,i}(s)ds<\infty.
 \eeqnn
Then by Lemma \ref{t4.6} and the dominated convergence one obtains
$\mbf{P}$-a.s.
 \beqnn
 \ar\ar
\lim_{l'\to\infty}\int_0^t
\sum_{j=1}^{2^{l'}}1_{({l'}_{j-1}T,{l'}_jT]}(s)
|J_{m,n,k,i}({l'}_jT)1_{\{{l'}_jT< \sigma_k\}}
-J_{m,n,k,i}(s)1_{\{s<\sigma_k\}}|ds \cr
 \ar\ar\qquad\le
\lim_{l'\to\infty}\int_0^t
\sum_{j=1}^{2^{l'}}1_{({l'}_{j-1}T,{l'}_jT]}(s)
|J_{m,n,k,i}({l'}_jT)-J_{m,n,k,i}(s)|1_{\{{l'}_jT< \sigma_k\}}ds \cr
 \ar\ar\qquad~~~
+ \lim_{l'\to\infty}\int_0^t
\sum_{j=1}^{2^{l'}}1_{({l'}_{j-1}T,{l'}_jT]}(s)
J_{m,n,k,i}(s)|1_{\{{l'}_jT< \sigma_k\}} -1_{\{s<\sigma_k\}}|ds=0,
 \eeqnn
which completes the proof. \qed

We are now ready to show Lemma \ref{t4.8}.

\noindent{\it Proof of Lemma \ref{t4.8}.} In the following let $t>0$
and $m,n,k\ge1$ be fixed. By Taylor's formula, dominated convergence
and Lemma \ref{t4.7} we have
 \beqlb\label{4.12}
I_{5,2}^{m,n,k,i}(t)
 \ar=\ar
\mbf{E}\Big\{\int_0^{t\wedge\sigma_{k}}ds
\int_{-K}^K\Psi_s(x)dx\int_{1/i}^\infty zm_0(dz)
\int_{-1}^1\Phi(y)dy \cr \ar\ar\qquad
\times\int_0^1\tilde{D}_n\big(\<U_s,\Phi_x^m\>,mzhV_s(x-\frac{y}{m})\big)
V_s(x-\frac{y}{m})dh \Big\}\cr
 \ar=\ar
\mbf{E}\Big\{\int_0^tJ_{m,n,k,i}(s) 1_{\{s<\sigma_k\}}ds \Big\} =
\mbf{E}\Big\{\bar{\int}_{(0,t)}
J_{m,n,k,i}(s)1_{\{s<\sigma_k\}}ds\Big\},
 \eeqlb
where recall that $\tilde{D}_n(y,z)=\phi_n'(y+z)-\phi_n'(y)$. Let
 \beqnn
\tilde{J}_{m,n,k,i}(t)
 \ar:=\ar
\int_{\mbb{R}}\Psi_t(x)dx\int_{1/i}^\infty zm_0(dz)
\int_{\mbb{R}}\Phi(y)dy \cr
 \ar\ar\quad\times
\int_0^1\tilde{D}_n\big(\<\tilde{U}_t,\Phi_x^m\>,mzh\tilde{V}_t(x-\frac{y}{m})\big)
\tilde{V}_t(x-\frac{y}{m})dh.
 \eeqnn
 Observe that for each
fixed $t>0$, $\tilde{U}_t$ and $\tilde{V}_t$ are the continuous
modifications of $U_t$ and $V_t$, respectively. Then it is
elementary to check that for each $t>0$,
$\<|\tilde{U}_t-U_t|,1\>=\<|\tilde{V}_t-V_t|,1\>=0$, $\mbf{P}$-a.s.
This implies $\tilde{J}_{m,n,k,i}(t)=J_{m,n,k,i}(t)$ for all
$t\in(0,\infty)\cap\mbb{Q}$, $\mbf{P}$-a.s. Together with
\eqref{4.12} we have $\mbf{P}$-a.s.
 \beqlb\label{4.6}
I_{5,2}^{m,n,k,i}(t) \le \mbf{E}\Big\{\bar{\int}_{(0,t)}
\tilde{J}_{m,n,k,i}(s)1_{\{s<\sigma_k\}}ds\Big\}.
 \eeqlb
For fixed $s$ and $x$ let $x_{s,m}\in[-1,1]$ be a value satisfying
 \beqnn
|\tilde{V}_s(x-\frac{x_{s,m}}{m})|
=\inf_{y\in[-1,1]}\{|\tilde{V}_s(x-\frac{y}{m})|\}.
 \eeqnn
It follows from \eqref{4.6} that
 \beqlb\label{4.7}
 \ar\ar
I_{5,2}^{m,n,k,i}(t) \cr
 \ar\le\ar
\mbf{E}\Big\{\int_0^{t\wedge\sigma_{k}}ds
\int_{-K}^K\Psi_s(x)dx\int_{1/i}^\infty zm_0(dz)
\int_{-1}^1\Phi(y)dy \cr \ar\ar\quad
\times\int_0^1\Big|\tilde{D}_n\big(\<\tilde{U}_s,\Phi_x^m\>,mzh\tilde{V}_s(x-\frac{y}{m})\big)
[\tilde{V}_s(x-\frac{y}{m})-\tilde{V}_s(x-\frac{x_{s,m}}{m})]\Big|dh
\Big\}\cr
 \ar\ar
+ \mbf{E}\Big\{\int_0^tds \int_{-K}^K\Psi_s(x)dx\int_{1/i}^\infty
zm_0(dz) \int_{-1}^1\Phi(y)dy\cr
 \ar\ar\quad
\times\int_0^1\Big|\tilde{D}_n\big(\<\tilde{U}_s,\Phi_x^m\>,mzh\tilde{V}_s(x-\frac{y}{m})\big)
\tilde{V}_s(x-\frac{x_{s,m}}{m})\Big|
1_{\{\tilde{V}_s(x-\frac{x_{s,m}}{m})\neq0\}}dh \Big\} \cr
 \ar=:\ar
I_{5,2,1}^{m,n,k,i}(t)+I_{5,2,2}^{m,n,i}(t).
 \eeqlb
We can finish the proof in  two steps.

{\bf Step 1.} We first estimate $I_{5,2,1}^{m,n,k,i}(t)$. Since for
fixed $s>0$, $\tilde{X}_s$ and $\tilde{Y}_s$ are the continuous
modifications of $X_s$ and $Y_s$, respectively, then we have
$\mbf{P}$-a.s.
 \beqnn
X_s(x)=\tilde{X}_s(x),\,\,Y_s(x)=\tilde{Y}_s(x), \qquad
x\in\mbb{R}\cap \mbb{Q},\,\,s\in(0,\infty)\cap \mbb{Q}.
 \eeqnn
Combining this with the definition of $\sigma_k$ and
$\bar{\int}_{(0,t]}$, we have $\mbf{P}$-a.s.
 \beqlb\label{4.9}
\sigma_k \ar=\ar \tilde{\sigma}_k :=\inf\Big\{t\in (0,T]
:\bar{\int}_{(0,t]} \sup_{x,z\in [-(K+1),K+1]\cap \mbb{Q},x\neq z}
\cr
 \ar\ar\qqquad\qqquad\qqquad
\frac{|\tilde{X}_s(x)-\tilde{X}_s(z)|\vee
|\tilde{Y}_s(x)-\tilde{Y}_s(z)|}{|x-z|^\eta}ds> k\Big\} \cr
 \ar=\ar
\inf\Big\{t\in (0,T]:\bar{\int}_{(0,t]}\sup_{x,z\in [-(K+1),K+1],x\neq z} \cr
 \ar\ar\qqquad\qqquad\qqquad
\frac{|\tilde{X}_s(x)-\tilde{X}_s(z)|\vee
|\tilde{Y}_s(x)-\tilde{Y}_s(z)|}{|x-z|^\eta}ds> k\Big\}.
 \eeqlb
By the H\"older inequality,
 \beqnn
 \ar\ar
\bar{\int}_{(0,t\wedge\tilde{\sigma}_k)} \sup_{|x|\le
K,|y|\vee|v|\le1}|\tilde{X}_s(x-\frac{y}{m})
-\tilde{X}_s(x-\frac{v}{m})|^{\beta}ds \cr
 \ar\ar\quad\le
t^{1-\beta}\Big[\bar{\int}_{(0,t\wedge\tilde{\sigma}_k)}
\sup_{|x|\le K,|y|\vee|v|\le1}|\tilde{X}_s(x-\frac{y}{m})
-\tilde{X}_s(x-\frac{v}{m})|ds\Big]^\beta \cr
 \ar\ar\quad\le
t^{1-\beta}(2/m)^{\eta\beta}\Big[\bar{\int}_{(0,t\wedge\tilde{\sigma}_k)}
\sup_{|x|\le K,|y|\vee|v|\le1,y\neq v}
\frac{|\tilde{X}_s(x-\frac{y}{m}) -\tilde{X}_s(x-\frac{v}{m})|}
{|y/m-v/m|^\eta}ds\Big]^\beta \cr
 \ar\ar\quad\le
2^{\eta\beta}t^{1-\beta}m^{-\eta\beta}k^{\beta},
 \eeqnn
and the same estimation holds for $\tilde{Y}$. Then by
\eqref{4.9} we have  $\mbf{P}$-a.s.
 \beqnn
 \ar\ar
\bar{\int}_{(0,t\wedge\sigma_k)} \sup_{|x|\le K,|y|\vee|v|\le1}
|\tilde{V}_s(x-\frac{y}{m})-\tilde{V}_s(x-\frac{v}{m})|ds \cr
 \ar\le\ar
C\bar{\int}_{(0,t\wedge\tilde{\sigma}_k)}\sup_{|x|\le
K,|y|\vee|v|\le1} \Big[|\tilde{X}_s(x-\frac{y}{m})
-\tilde{X}_s(x-\frac{v}{m})|^{\beta} +|\tilde{Y}_s(x-\frac{y}{m})
-\tilde{Y}_s(x-\frac{v}{m})|^{\beta}\Big]ds \cr
 \ar\le\ar
2^{\eta\beta+1}Ct^{1-\beta}m^{-\eta\beta}k^{\beta}.
 \eeqnn
Observe that $|\tilde{D}_n(y,z)|\le2$ for all $n\ge1$ and
$y,z\in\mbb{R}$. It then follows that
 \beqlb\label{4.10}
I_{5,2,1}^{m,n,k,i}(t)
 \ar\le\ar
\frac{2c_2i^{\alpha-1}}{\alpha-1}\mbf{E}\Big\{\bar{\int}_{(0,t\wedge\sigma_k)}ds
\int_{-K}^Kdx \int_{-1}^1 |\tilde{V}_s(x-\frac{y}{m})-
\tilde{V}_s(x-\frac{x_{s,m}}{m})|\Big\} \Phi(y)dy \cr
 \ar\le\ar
\frac{2c_2i^{\alpha-1}}{\alpha-1}\int_{-K}^Kdx \int_{-1}^1
\mbf{E}\Big\{\bar{\int}_{(0,t\wedge\sigma_k)}  \cr
 \ar\ar\qquad\qquad
\sup_{|x|\le K,|y|\vee|v|\le1}
|\tilde{V}_s(x-\frac{y}{m})-\tilde{V}_s(x-\frac{v}{m})| ds\Big\}
\Phi(y)dy \cr
 \ar\le\ar
2^{\eta\beta+2}k^{\beta}c_2C Kt^{1-\beta}
m^{-\eta\beta}i^{\alpha-1}(\alpha-1)^{-1},
 \eeqlb
where $c_2=\sup_{(s,x)\in[0,T]\times[-K,K]}\Psi_s(x)$.

{\bf Step 2.} We then estimate $I_{5,2,2}^{m,n,i}(t)$. Since
$supp(\phi_n'')\subset (a_n,a_{n-1})$, then
$\tilde{D}_n\big(y,z\big)=0$ for $y\ge a_{n-1}$ and $z\ge0$. It then
follows that for each $y,z\ge0$,
 \beqlb\label{4.23}
\tilde{D}_n(y,z)=\tilde{D}_n(y,z)1_{\{|y|<a_{n-1}\}}.
 \eeqlb
One can also get \eqref{4.23} for the case $y,z\le0$.

By the H\"older continuity of $H$, there is a constant $c_3>0$ so
that $|\tilde{U}_s(x-\frac{u}{m})|\ge
c_3|\tilde{V}_s(x-\frac{u}{m})|^{1/\beta}$ for all $u\in[-1,1]$.
Then
 \beqlb\label{4.24} \{|\<\tilde{U}_s,\Phi_x^m\>|<a_{n-1}\}
\subset \Big\{|\tilde{V}_s(x-\frac{x_{s,m}}{m})|
<(c_3^{-1}a_{n-1})^\beta\Big\}.
 \eeqlb
To verify \eqref{4.24}, if
 \beqnn
|\tilde{V}_s(x-\frac{x_{s,m}}{m})| \ge(c_3^{-1}a_{n-1})^\beta,
 \eeqnn
then $\tilde{V}_s(x-\frac{x_{s,m}(x)}{m})\neq0$. This implies that
$\tilde{V}_s(x-\frac{u}{m})\neq0$ for all $u\in[-1,1]$. Then by the
continuity of $u\mapsto \tilde{V}_s(u)$ and the mean value theorem,
$\tilde{V}_s(x-\frac{u}{m})>0$ for all $u\in[-1,1]$, or
$\tilde{V}_s(x-\frac{u}{m})<0$ for all $u\in[-1,1]$. On the other
hand, $H$ is a nondecreasing function (condition (C4) in Section 1).
Then $\tilde{U}_s(x-\frac{u}{m})>0$ as
$\tilde{V}_s(x-\frac{u}{m})>0$ and $\tilde{U}_s(x-\frac{u}{m})<0$ as
$\tilde{V}_s(x-\frac{u}{m})<0$. Therefore,
 \beqnn
 \ar\ar
\Big|\int_{\mbb{R}}\tilde{U}_s(x-\frac{u}{m})\Phi(u)du\Big|
=\int_{-1}^1|\tilde{U}_s(x-\frac{u}{m})|\Phi(u)du \cr
 \ar\ar\qquad\ge
c_3\int_{-1}^1|\tilde{V}_s(x-\frac{u}{m})|^{\frac{1}{\beta}}\Phi(u)du
\ge c_3|\tilde{V}_s(x-\frac{x_{s,m}(x)}{m})|^{\frac{1}{\beta}} \ge
a_{n-1},
 \eeqnn
which implies \eqref{4.24}.

By \eqref{4.23} one can also see that
 \beqnn
 \ar\ar
\tilde{D}_n\big(\<\tilde{U}_s,\Phi_x^m\>,mzh\tilde{V}_s(x-\frac{x_{s,m}}{m})\Phi(y)\big)
1_{\{\tilde{V}_s(x-\frac{x_{s,m}}{m})\neq0\}} \cr
 \ar\ar\quad=
\tilde{D}_n\big(\<\tilde{U}_s,\Phi_x^m\>,mzh\tilde{V}_s(x-\frac{x_{s,m}}{m})\Phi(y)\big)
1_{\{|\tilde{V}_s(x-\frac{x_{s,m}}{m})\neq0,\,\<\tilde{U}_s,\Phi_x^m\>|<a_{n-1}\}}.
 \eeqnn
Putting together \eqref{4.24} with $|\tilde{D}_n(y,z)|\le2$ for all
$y,z\in\mbb{R}$ and $m\ge1$ we have that
 \beqnn
 \ar\ar
\Big|\tilde{D}_n\big(\<\tilde{U}_s,\Phi_x^m\>,mzh\tilde{V}_s(x-\frac{x_{s,m}}{m})
\Phi(y)\big)
\tilde{V}_s(x-\frac{x_{s,m}}{m})\Big| \cr
 \ar\ar\qquad\le
2|\tilde{V}_s(x-\frac{x_{s,m}}{m})|
1_{\{\<\tilde{U}_s,\Phi_x^m\>|<a_{n-1}\}}
\le2(c_3^{-1}a_{n-1})^\beta.
 \eeqnn
Then
 \beqnn
I_{5,2,2}^{m,n,i}(t) \le 2(c_3^{-1}a_{n-1})^\beta \int_0^tds
\int_{-K}^K\Psi_s(x)dx\int_{1/i}^\infty zm_0(dz) \le
C_Ta_{n-1}^\beta i^{\alpha-1},\qquad t\in[0,T].
 \eeqnn
Combining  with \eqref{4.7} and \eqref{4.10}, we finish the proof.
 \qed

\subsection{A remark on the proof}
We remark that we do not consider the increased H\"older regularity
near its zero (used in Mytnik and Perkins (2011) for proving the
pathwise uniqueness of SPDE driven by Gaussian white noise), i.e.
the difference of two solutions is jointly H\"older continuous with
the H\"older exponent in space in (0,1) and with H\"older exponent
in time in (0,1/2) when the difference is close to zero. But for the
SPDE \eqref{1.1}, it is hard to establish the similar result of
H\"older regularity for the difference of two solutions because the
regularities of the solutions in time for fixed spacial point may be
bad. For example, for super-Brownian motion (i.e. $G\equiv0$,
$H(x)=x^\beta$ and $p=1$ in \eqref{1.1}) it was proved in
\cite[Theorem 1.2]{MP03} that for any $t,\delta>0$ and almost every
spatial point $x\in\mbb{R}$ fixed, the essential supremum of the
solution over time  interval $(t,t+\delta)$ is infinity.

In this paper the proof of pathwise uniqueness for the solution to
\eqref{1.01} relies on the H\"older continuity of the solution at a
fixed time. If one uses the H\"older continuity at any given spatial
point where the H\"older exponent $\bar{\eta}_c=(3/\alpha-1)\wedge1$
is bigger than $\eta_c$, it appears that the criterion for pathwise
uniqueness could be improved with
$\beta>\frac{(\alpha-1)(\bar{\eta}_c+1)}{(2-\alpha)\bar{\eta}_c}$
in Theorem \ref{t4.1}
(with $\eta_c$ replaced by $\bar{\eta}_c$ in \eqref{4.25} and
$1<\alpha<(\sqrt{17}+1)/4$ for case $p=1$). But there is a problem
with this approach, which we explain below.
For the stoping time $\sigma_k $ one gets an equation similar to \eqref{4.1} with $t$ replaced by $t\wedge
\sigma_k(x)$. Then by the same argument
as in \eqref{4.2} we have
 \beqnn
 \ar\ar
\<\phi_n(\<U_{t\wedge\tau_k(\cdot)},\Phi_{\cdot}^m\>),\Psi_t\> \cr
 \ar=\ar
\frac12\int_0^t\<\phi_n'(\<U_s,\Phi_{\cdot}^m\>)
\<U_s,\Delta\Phi_{\cdot}^m\>1_{\{{s\le\tau_k(\cdot)}\}},\Psi_s\>ds \cr
 \ar\ar
+\int_0^t\<\phi_n'(\<U_s,\Phi_{\cdot}^m\>)\<R_s,\Phi_{\cdot}^m\>
1_{\{{s\le\tau_k(\cdot)}\}},\Psi_s\>ds
+\int_0^t\<\phi_n(\<U_{s\wedge\tau_k(\cdot)},\Phi_{\cdot}^m\>),\dot{\Psi}_s\>ds
\cr
 \ar\ar
+\int_0^t\int_0^\infty\int_{\mbb{R}}
\<H_n(\<U_{s-},\Phi_{\cdot}^m\>,zV_{s-}(y)\Phi_{\cdot}^m(y))
1_{\{{s\le\tau_k(\cdot)}\}},\Psi_s\>\tilde{N}(ds,dz,dy)\cr
 \ar\ar
+\int_0^tds\int_0^\infty m_0(dz)\int_{\mbb{R}}
\<D_n(\<U_s,\Phi_{\cdot}^m\>,zV_s(y)\Phi_{\cdot}^m(y))
1_{\{{s\le\tau_k(\cdot)}\}},\Psi_s\>dy \cr
 \ar=:\ar
\hat{I}_1^{m,n}(t,k)+\hat{I}_2^{m,n}(t,k)+\hat{I}_3^{m,n}(t,k)
+\hat{I}_4^{m,n}(t,k)+\hat{I}_5^{m,n}(t,k),
 \eeqnn
where $U_t,V_t,\Phi_x^m,D_n,\phi_n,\Psi_t,\gamma_k$ defined below
\eqref{1.06}, $\tau_k(x):=\gamma_k\wedge\sigma_k(x)$,
$H_n(y,z):=\phi_n(y+z)-\phi_n(y)$ and $R_s$ denotes the difference
of compositions of the two solutions into function $G$. As Lemmas
\ref{t4.2}--\ref{t4.3}, for each $k\ge1$ one can get
 \beqnn
\lim_{m,n\to\infty}\mbf{E}\{
\hat{I}_2^{m,n}(t,k)\}
=
\mbf{E}\Big\{\int_0^tds\int_{\mbb{R}}
\mathop{\rm sgn}(U_s(x))Q_s(x)\Psi_s(x)1_{\{{s\le\tau_k(x)}\}}dx\Big\}
 \eeqnn
and
 \beqnn
\lim_{m,n\to\infty}\mbf{E}\{
\hat{I}_3^{m,n}(t,k)\}
=\mbf{E}\Big\{\int_0^tds\int_{\mbb{R}}
|U_{s\wedge\tau_k(x)}(x)|\dot{\Psi}_s(x)dx\Big\},
\quad \mbf{E}\{\hat{I}_4^{m,n}(t,k)\}=0.
 \eeqnn
Similar to Lemma \ref{t4.4}, we also have that if
$m=a_{n-1}^{-\delta}$ for $\delta>0$, then for each $t>0$ and
$k\ge1$,
 \beqnn
\lim_{m,n\to\infty}\mbf{E}\{\hat{I}_5^{m,n}(t,k)\}=0 .
 \eeqnn
But it is hard to deal with
 \beqnn
\mbf{E}\{
\hat{I}_1^{m,n}(t,k)\}.
 \eeqnn
The difficulty comes from the fact  that
$x\mapsto1_{\{{s\le\tau_k(x)}\}}$ is not continuous. So we cannot
use the same argument as in Lemma 2.2(b) of \cite{MPS06} to obtain
an inequality like \eqref{4.5}.

\section{Appendix: proofs of Proposition \ref{t1.5} and Lemma \ref{t1.4}}

\setcounter{equation}{0}

Before proving Proposition \ref{t1.5}, we state a lemma.

\blemma\label{t1.2b} Let $t\in[0,T]$ be fixed. For any $k\ge1$, $\lambda>0$ and
$f\in C(\mbb{R})$ satisfying $\lambda_0(|f|)<\infty$ we have
$\mbf{P}$-a.s.
 \beqlb\label{2.10}
 \ar\ar
\<X_{t\wedge
\tilde{\tau}_k},P_{t-(t\wedge\tilde{\tau}_k)+\lambda}f\>
=X_0(P_{t+\lambda}f)+
\int_0^t\<G(X_s),P_{t-s+\lambda}f\>1_{\{s\le\tilde{\tau}_k\}}ds  \cr
 \ar\ar\qquad
 + \int_0^t\int_0^\infty\int_{\mbb{R}}\int_0^{H(X_{s-}(u))^\alpha}
z  P_{t-s+\lambda}f(u)1_{\{s\le
\tilde{\tau}_k\}}\tilde{N}_0(ds,dz,du,dv),
 \eeqlb
where $\tilde{\tau}_k$ is the stopping time defined in \eqref{1.21}.
 \elemma
 \proof
We consider a partition
$\Delta_n:=\{0=t_0<t_1<\cdots<t_n=t\}$ of $[0,t]$. Let
$|\Delta_n|:=\max_{1\le i\le n}|t_i-t_{i-1}|$. Let $f_\lambda:=P_\lambda f$.
It is clear that
$\frac{dP_sf_\lambda(x)}{ds}=\frac12 P_sf_\lambda''(x)$ for $s\ge0$.
For $k\ge1$ and $s\in[0,T]$, let $Z_k(s)=X_{s\wedge \tilde{\tau}_k}$.
By Proposition \ref{t1.9}
 \beqnn
\<Z_k(t),f_\lambda\>
 \ar=\ar
X_0(f_\lambda)+
\frac12\int_0^t\<X_s,f_\lambda''\>1_{\{s\le\tilde{\tau}_k\}}ds +
\int_0^t\<G(X_s),f_\lambda\>1_{\{s\le\tilde{\tau}_k\}}ds \cr
 \ar\ar
+ \int_0^t\int_0^\infty\int_{\mbb{R}}\int_0^{H(X_{s-}(u))^\alpha}
z f_\lambda(u)1_{\{s\le
\tilde{\tau}_k\}}\tilde{N}_0(ds,dz,du,dv).
 \eeqnn
It follows that
 \beqlb\label{2.8}
 \ar\ar
\<Z_k(t),P_{t-(t\wedge\tilde{\tau}_k)}f_\lambda\> \cr
 \ar=\ar
X_0(P_{t+\lambda}f) + \sum\limits_{i=1}^n
\<Z_k(t_i),P_{t-(t_i\wedge\tilde{\tau}_k)}f_\lambda
-P_{t-(t_{i-1}\wedge\tilde{\tau}_k)}f_\lambda\>
\cr
 \ar\ar
+ \sum\limits_{i=1}^n
\big[\<Z_k(t_i),P_{t-(t_{i-1}\wedge\tilde{\tau}_k)}f_\lambda\> -
\<Z_k(t_{i-1}),P_{t-(t_{i-1}\wedge\tilde{\tau}_k)}f_\lambda\>\big]
\cr
 \ar=\ar
X_0(P_{t+\lambda}f) +
\frac12\sum\limits_{i=1}^n
\int_{t-(t_{i-1}\wedge\tilde{\tau}_k)}^{t-(t_i\wedge\tilde{\tau}_k)}
\<Z_k(t_i),P_sf_\lambda''\>ds \cr
 \ar\ar
+ \frac12\sum\limits_{i=1}^n\int_{t_{i-1}}^{t_i}
\<X_s,P_{t-(t_{i-1}\wedge\tilde{\tau}_k)}f_\lambda''\>1_{\{s\le\tilde{\tau}_k\}}ds
+\sum\limits_{i=1}^n\int_{t_{i-1}}^{t_i}
\<G(X_s),P_{t-(t_{i-1}\wedge\tilde{\tau}_k)+\lambda}f\>1_{\{s\le\tilde{\tau}_k\}}ds
\cr
 \ar\ar
+ \sum\limits_{i=1}^n\int_{t_{i-1}}^{t_i}
\int_0^\infty\int_{\mbb{R}}\int_0^{H(X_{s-}(u))^\alpha }
z
P_{t-(t_{i-1}\wedge\tilde{\tau}_k)+\lambda}f(u)1_{\{s\le
\tilde{\tau}_k\}}\tilde{N}_0(ds,dz,du,dv) \cr
 \ar=\ar
X_0(P_{t+\lambda}f) + \frac12\int_0^t  \sum\limits_{i=1}^nI_i(s)
\big[\<X_s,P_{t-(t_{i-1}\wedge\tilde{\tau}_k)}f_\lambda''\>
-\<Z_k(t_i),P_{t-s}f_\lambda''\>\big]1_{\{s\le \tilde{\tau}_k\}}ds
\cr \ar\ar +\int_0^t  \sum\limits_{i=1}^nI_i(s)
\<G(X_s),P_{t-(t_{i-1}\wedge\tilde{\tau}_k)+\lambda}f\>1_{\{s\le\tilde{\tau}_k\}}ds
\cr
 \ar\ar
+ \int_0^t\int_0^\infty\int_{\mbb{R}}\int_0^{H(X_{s-}(u))^\alpha}
z \sum\limits_{i=1}^nI_i(s)
P_{t-(t_{i-1}\wedge\tilde{\tau}_k)+\lambda}f(u)1_{\{s\le
\tilde{\tau}_k\}}\tilde{N}_0(ds,dz,du,dv),
 \eeqlb
where $I_i(s):=1_{(t_{i-1},t_i]}(s)$.

Since $\|f_\lambda''\|<\infty$ and $\<X_s,1\>\le k$ on
$\{s\le \tilde{\tau}_k\}$,
then by the dominated convergence, $\mbf{P}$-a.s.
 \beqlb\label{2.9}
 \ar\ar
\lim\limits_{|\Delta_n|\to0}\int_0^t\sum\limits_{i=1}^nI_i(s)
|\<X_s,P_{t-(t_{i-1}\wedge\tilde{\tau}_k)}f_\lambda''\> -
\<Z_k(t_i),P_{t-s}f_\lambda''\>|1_{\{s\le \tilde{\tau}_k\}}ds \cr
 \ar\ar\le
\lim\limits_{|\Delta_n|\to0}\int_0^t\sum\limits_{i=1}^nI_i(s)
|\<X_s,P_{t-(t_{i-1}\wedge\tilde{\tau}_k)}f_\lambda''-P_{t-s}f_\lambda''\>|1_{\{s\le
\tilde{\tau}_k\}}ds \cr
 \ar\ar\quad
+ \lim\limits_{|\Delta_n|\to0}\int_0^t \sum\limits_{i=1}^nI_i(s)
|\<X_s,P_{t-s}f_\lambda''\>-\<X_{t_i\wedge\tilde{\tau}_k},P_{t-s}f_\lambda''\>|1_{\{s\le
\tilde{\tau}_k\}} ds \cr
 \ar\ar\le
\int_0^t\sum\limits_{i=1}^nI_i(s)
\lim\limits_{|\Delta_n|\to0}|
\<X_s,P_{t-(t_{i-1}\wedge\tilde{\tau}_k)}f_\lambda''-P_{t-s}f_\lambda''\>|1_{\{s\le
\tilde{\tau}_k\}}ds \cr
 \ar\ar\quad
+\int_0^t \sum\limits_{i=1}^n\lim\limits_{|\Delta_n|\to0}I_i(s)
|\<X_s,P_{t-s}f_\lambda''\>-\<X_{t_i\wedge\tilde{\tau}_k},P_{t-s}f_\lambda''\>|
1_{\{s\le
\tilde{\tau}_k\}} ds =0,
 \eeqlb
where the right continuities of  $t'\mapsto P_{t'}f_\lambda''$ and
$t'\mapsto \<X_{t'},P_{t-s}f_\lambda''\>$ were used in the last equation.

By the Lipschitz condition of $G$ and the dominated convergence we
can also have $\mbf{P}$-a.s.
 \beqlb\label{2.11}
 \ar\ar
\lim\limits_{|\Delta_n|\to0}\int_0^t  \sum\limits_{i=1}^nI_i(s)
|\<G(X_s),P_{t-(t_{i-1}\wedge\tilde{\tau}_k)+\lambda}f\>
-\<G(X_s),P_{t-s+\lambda}f\>|1_{\{s\le\tilde{\tau}_k\}}ds \cr
 \ar\ar\le
\lim\limits_{|\Delta_n|\to0}\int_0^t  \sum\limits_{i=1}^nI_i(s)
\<G(X_s),|P_{t-(t_{i-1}\wedge\tilde{\tau}_k)+\lambda}f
-P_{t-s+\lambda}f|\>|1_{\{s\le\tilde{\tau}_k\}}ds =0.
 \eeqlb

Observe that for $s\in[0,T]$
 \beqnn
f(s,u,\lambda,n,k)
 :=
\sum\limits_{i=1}^nI_i(s)
|P_{t-(t_{i-1}\wedge\tilde{\tau}_k)+\lambda}f(u) -
P_{t-s+\lambda}f(u)| \le 2\|f\|
 \eeqnn
and $f(s,u,\lambda,n,k)$ converges to zero by the right continuity of $t'\mapsto P_{t'}f$ for $s\le \tilde{\tau}_k$ as $\Delta_n\to0$. By the same argument as in
\eqref{1.32} and \eqref{1.17},
 \beqnn
\mbf{E}\Big\{\int_0^{t\wedge\tilde{\tau}_k}ds\int_{\mbb{R}}
[1+X_s(u)^p] P_{\lambda+T}f(u)du\Big\}<\infty
 \eeqnn
for each $k\ge1$. Then by the dominated convergence and
Burkholder-Davis-Gundy inequality it is easy to
see that as $\Delta_n\to0$,
 \beqlb\label{2.7}
 \ar\ar
\mbf{E}\Big\{\Big|\int_0^t\int_1^\infty\int_{\mbb{R}}\int_0^{H(X_{s-}(u))^\alpha }
 z \Big[\sum\limits_{i=1}^nI_i(s)
P_{t-(t_{i-1}\wedge\tilde{\tau}_k)+\lambda}f(u) \cr
 \ar\ar\qqquad\qqquad\qquad\qqquad\quad
- P_{t-s+\lambda}f(u)\Big] 1_{\{s\le \tilde{\tau}_k\}}
\tilde{N}_0(ds,dz,du,dv)\Big|\Big\} \cr
 \ar\ar\qquad\le
2\int_1^\infty zm_0(dz)\mbf{E}\Big\{\int_0^tds\int_{\mbb{R}}
H(X_s(u))^\alpha f(s,u,\lambda,n,k) 1_{\{s\le \tilde{\tau}_k\}}du\Big\}
\cr
 \ar\ar\qquad\le
2C\int_1^\infty zm_0(dz)\,\mbf{E}\Big\{\int_0^tds\int_{\mbb{R}}
[1+X_s(u)^p] f(s,u,\lambda,n,k) 1_{\{s\le \tilde{\tau}_k\}}du\Big\} \to0
 \eeqlb
and
 \beqlb\label{2.13}
\ar\ar
\mbf{E}\Big\{\Big|\int_0^t\int_0^1\int_{\mbb{R}}\int_0^{H(X_{s-}(u))^\alpha }
 z \Big[\sum\limits_{i=1}^nI_i(s)
P_{t-(t_{i-1}\wedge\tilde{\tau}_k)+\lambda}f(u) \cr
 \ar\ar\qqquad\qquad\qquad\qquad\qquad\quad
- P_{t-s+\lambda}f(u)\Big] 1_{\{s\le \tilde{\tau}_k\}}
\tilde{N}_0(ds,dz,du,dv)\Big|^2\Big\} \cr
 \ar\ar\quad\le
C\|f\|\int_0^1z^2m_0(dz)
\mbf{E}\Big\{\int_0^tds\int_{\mbb{R}}[1+X_s(u)^p]f(s,u,\lambda,n,k)
1_{\{s\le \tilde{\tau}_k\}}du\Big\}
 \to0.
 \eeqlb
Now it is obvious that \eqref{2.10} follows from
\eqref{2.8}--\eqref{2.11} and  \eqref{2.7}--\eqref{2.13}. \qed

Now we are ready to present proof of Proposition \ref{t1.5}.

\noindent{\it Proof of \eqref{1.2}.} Recall that the stopping time
$\tilde{\tau}_k$ is defined in \eqref{1.21}. Let $f\in B(\mbb{R})$
with $\lambda_0(|f|)<\infty$ in this step. For each $n\ge1$ and
$x\in\mbb{R}$ define $f_n(x)=n\int_{x-1/n}^xf(y)dy$. Then $f_n\in
C(\mbb{R})$ and $\lambda_0(|f_n|)\le\lambda_0(|f|)<\infty$ by
integration by parts. Then \eqref{2.10} holds with $f$ replaced by
$f_n$ by Lemma \ref{t1.2b}. By the right continuity of $t'\mapsto
P_{t'}f_n$ and the same argument in \eqref{2.7} and \eqref{2.13},
 \beqnn
\mbf{E}\Big\{\Big|\int_0^t\int_0^\infty\int_{\mbb{R}}\int_0^{H(X_{s-}(u))^\alpha}
z  [P_{t-s+\lambda}f_n(u)-P_{t-s}f_n(u)]1_{\{s\le
\tilde{\tau}_k\}}\tilde{N}_0(ds,dz,du,dv)\Big|\Big\}
 \to0
 \eeqnn
as $\lambda\to0$.
Since \eqref{2.10} holds with $f$ replaced by $f_n$,
taking $\lambda\to0$ we get
 \beqlb\label{2.10b}
 \ar\ar
\<X_{t\wedge
\tilde{\tau}_k},P_{t-(t\wedge\tilde{\tau}_k)}f_n\>
=X_0(P_tf_n)+
\int_0^t\<G(X_s),P_{t-s}f_n\>1_{\{s\le\tilde{\tau}_k\}}ds  \cr
 \ar\ar\qquad
 + \int_0^t\int_0^\infty\int_{\mbb{R}}\int_0^{H(X_{s-}(u))^\alpha}
z  P_{t-s}f_n(u)1_{\{s\le
\tilde{\tau}_k\}}\tilde{N}_0(ds,dz,du,dv).
 \eeqlb
Observe that $\|f_n\|\le \|f\|<\infty$
and $\lim_{n\to\infty}f_n(x)=f(x)$, $\lambda_0$-a.e. $x$.
Then letting $n\to\infty$ in \eqref{2.10b},
by the dominated convergence and the same argument in \eqref{2.7}--\eqref{2.13} again,
we obtain
 \beqlb\label{2.10c}
 \ar\ar
\<X_{t\wedge
\tilde{\tau}_k},P_{t-(t\wedge\tilde{\tau}_k)}f\>
=X_0(P_{t}f)+
\int_0^t\<G(X_s),P_{t-s}f\>1_{\{s\le\tilde{\tau}_k\}}ds  \cr
 \ar\ar\qquad
 + \int_0^t\int_0^\infty\int_{\mbb{R}}\int_0^{H(X_{s-}(u))^\alpha}
z  P_{t-s}f(u)1_{\{s\le
\tilde{\tau}_k\}}\tilde{N}_0(ds,dz,du,dv),
 \eeqlb
which implies \eqref{1.2} by taking $k\to\infty$.
\qed

\noindent{\it Proof of \eqref{1.3}}.
Let $t>0$ and $f\in B(\mbb{R})$ with $\lambda_0(|f|)<\infty$ be fixed.
By Fubini's theorem,
 \beqlb\label{1.18}
 \ar\ar
\int_0^{t\wedge \tilde{\tau}_k}\int_1^\infty
\int_{\mbb{R}}\int_0^{H(X_{s-}(u))^\alpha } z  P_{t-s}f(u) \tilde{N}_0(ds,dz,du,dv)  \cr
 \ar=\ar
\int_0^{t\wedge \tilde{\tau}_k}\int_1^\infty
\int_{\mbb{R}}\int_0^{H(X_{s-}(u))^\alpha } z  P_{t-s}f(u) N_0(ds,dz,du,dv) \cr
 \ar\ar
-\int_0^{t\wedge \tilde{\tau}_k}ds\int_1^\infty m_0(dz)
\int_{\mbb{R}}du \int_0^{H(X_{s-}(u))^\alpha } z  P_{t-s}f(u) dv \cr
 \ar=\ar
\int_{\mbb{R}}f(x)\Big[\int_0^{t\wedge
\tilde{\tau}_k}\int_1^\infty\int_{\mbb{R}}\int_0^{H(X_{s-}(u))^\alpha } z  p_{t-s}(x-u)
N_0(ds,dz,du,dv)\Big]dx \cr
 \ar\ar
-\int_{\mbb{R}}f(x)\Big[\int_0^{t\wedge
\tilde{\tau}_k}ds\int_1^\infty
m_0(dz)\int_{\mbb{R}}du\int_0^{H(X_{s-}(u))^\alpha } z  p_{t-s}(x-u) dv\Big]dx \cr
 \ar=\ar
\int_{\mbb{R}}f(x)\Big[\int_0^{t\wedge
\tilde{\tau}_k}\int_1^\infty\int_{\mbb{R}}\int_0^{H(X_{s-}(u))^\alpha } z
p_{t-s}(x-u) \tilde{N}_0(ds,dz,du,dv)\Big]dx.
 \eeqlb
By stochastic Fubini's theorem (see e.g. \cite[Theorem 7.24]{Li11}),
to prove
 \beqlb\label{1.19}
 \ar\ar
\int_0^{t\wedge \tilde{\tau}_k}\int_0^1
\int_{\mbb{R}}\int_0^{H(X_{s-}(u))^\alpha } z  P_{t-s}f(u) \tilde{N}_0(ds,dz,du,dv)  \cr
 \ar=\ar
\int_{\mbb{R}}f(x)\Big[\int_0^{t\wedge
\tilde{\tau}_k}\int_0^1\int_{\mbb{R}}\int_0^{H(X_{s-}(u))^\alpha } z
p_{t-s}(x-u) \tilde{N}_0(ds,dz,du,dv)\Big]dx,~~  \mbf{P}\mbox{-a.s.},
 \eeqlb
we only need to verify
 \beqlb\label{1.20}
 \ar\ar
\mbf{E}\Big\{\int_{\mbb{R}}f(x)dx\int_0^{t\wedge \tilde{\tau}_k}ds
\int_0^1m_0(dz)\int_{\mbb{R}}du\int_0^{H(X_{s-}(u))^\alpha }
z^2  p_{t-s}(x-u)^2dv\Big\} \cr
 \ar\le\ar
C\int_0^1 z^2 m_0(dz)\mbf{E}\Big\{\int_0^{t\wedge \tilde{\tau}_k}ds\int_{\mbb{R}}|f(x)|dx\int_{\mbb{R}}
[1+X_s(u)^p] p_{t-s}(x-u)^2du\Big\}<\infty.
 \eeqlb

Indeed, for the case $0<p<1$, by an argument similar to \eqref{1.32},
 \beqnn
\int_0^{t\wedge \tilde{\tau}_k}ds\int_{\mbb{R}}|f(x)|dx\int_{\mbb{R}}
X_s(u)^p p_{t-s}(x-u)^2du
\le
\lambda_0(|f|)\int_0^t[2\pi s]^{-\frac{\bar{\alpha}-1}{2}}
[1+[2\pi s]^{-1/2}k]ds<\infty.
 \eeqnn
For the case $p=1$,
 \beqnn
\int_0^{t\wedge \tilde{\tau}_k}ds\int_{\mbb{R}}|f(x)|dx\int_{\mbb{R}}
X_s(u) p_{t-s}(x-u)^2du
\le
k\|f\|\int_0^{t}[2\pi(t-s)^{-1/2}]ds<\infty.
 \eeqnn
Observe that
 \beqnn
X_s(u)^p p_{t-s}(x-u)^2=
[X_s(u)^{p-\delta}p_{t-s}(x-u)^{1-\delta}]
\times[X_s(u)^\delta p_{t-s}(x-u)^{1+\delta}].
 \eeqnn
For the case $1<p<2$ choosing $\delta=(q-p)/(q-1)\in(0,1)$ for $q$ given in Assumption 1.4, it is elementary to see that
$(p-\delta)/(1-\delta)\wedge \delta>1/2$ by the fact $q>3p/(3-\alpha)$.
Then by the H\"older inequality,
 \beqnn
 \ar\ar
\int_0^{t\wedge \tilde{\tau}_k}ds\int_{\mbb{R}}|f(x)|dx\int_{\mbb{R}}
X_s(u)^p p_{t-s}(x-u)^2du \cr
 \ar\ar\quad\le
\|f\|\Big[\int_0^{t\wedge \tilde{\tau}_k}ds\int_{\mbb{R}}dx\int_{\mbb{R}}
X_s(u)^q p_{t-s}(x-u)du\Big]^{1-\delta} \cr
 \ar\ar\qquad
\times\Big[\int_0^{t\wedge \tilde{\tau}_k}ds\int_{\mbb{R}}dx\int_{\mbb{R}}
X_s(u) p_{t-s}(x-u)^{1+1/\delta}du\Big]^{\delta} \cr
 \ar\ar\quad\le
\|f\|\Big[\int_0^{t\wedge \tilde{\tau}_k}ds\int_{\mbb{R}}
X_s(u)^q du\Big]^{1-\delta}\times\Big|\int_0^{t\wedge \tilde{\tau}_k}[2\pi(t-s)^{-1/2}]^{1/\delta}\<X_s,1\>ds\Big|^{\delta} \cr
 \ar\ar\quad\le
k\|f\|[2\pi]^{-1}\Big[\int_0^ts^{-1/(2\delta)}ds\Big]^{\delta}<\infty,
 \eeqnn
which implies \eqref{1.20}.

Combining \eqref{1.18}, \eqref{1.19} and \eqref{2.10c}, we have
$\mbf{P}$-a.s.
 \beqnn
 \ar\ar
\int_{\mbb{R}}X_{t\wedge \tilde{\tau}_k}(x)P_{t-{(t\wedge
\tilde{\tau}_k)}}f(x)dx \cr
 \ar=\ar
\int_{\mbb{R}}\Big[\int_{\mbb{R}}p_t(x-z)X_0(dz) +\int_0^{t\wedge
\tilde{\tau}_k}ds\int_{\mbb{R}}p_{t-s}(x-z)G(X_s(z))dz \cr
 \ar\ar
+ \int_0^{t\wedge \tilde{\tau}_k}
\int_0^\infty\int_{\mbb{R}}\int_0^{H(X_{s-}(u))^\alpha }
z  p_{t-s}(x-u)
\tilde{N}_0(ds,dz,du,dv) \Big]f(x)dx.
 \eeqnn
Letting $k\to\infty$ one completes the proof. \qed

\noindent{\it Proof of Lemma \ref{t1.4}.}
If \eqref{1.22} holds for
$1\vee p<\bar{p}<\alpha$,
the rest can be given by the Jensen inequality.
So in the following we always assume that $1\vee p<\bar{p}<\alpha$.

{\bf Step 1.}
Note that
 \beqlb\label{1.23}
 \ar\ar
\Big|\int_0^tds\int_{\mbb{R}}p_{t-s}(x-z)G(X_s(z))dz\Big|^{\bar{p}} \cr
 \ar\ar\quad\le
Ct^{\bar{p}}\int_0^tds\int_{\mbb{R}}p_{t-s}(x-z)G(X_s(z))^{\bar{p}}dz \cr
 \ar\ar\quad\le
Ct^{\bar{p}}\int_0^tds\int_{\mbb{R}}p_{t-s}(x-z)[1+X_s(z)^{\bar{p}}]dz.
 \eeqlb
Recalling the stopping time $\tilde{\tau}_k$ defined in
\eqref{1.21}, one can see that
 \beqnn
\bar{Z}_k(t,x)
 \ar:=\ar
\int_0^{t\wedge\tilde{\tau}_k} \int_0^\infty\int_{\mbb{R}}
\int_0^{H(X_{s-}(u))^\alpha } z
p_{t-s}(x-u)\tilde{N}_0(ds,dz,du,dv)   \cr
 \ar=\ar
\int_0^t \int_0^1\int_{\mbb{R}} \int_0^{H(X_{s-}(u))^\alpha }
z
p_{t-s}(x-u)1_{\{s\le\tilde{\tau}_k\}} \tilde{N}_0(ds,dz,du,dv)\cr
 \ar\ar
+ \int_0^t \int_1^\infty\int_{\mbb{R}}
\int_0^{H(X_{s-}(u))^\alpha } z
p_{t-s}(x-u)1_{\{s\le\tilde{\tau}_k\}} \tilde{N}_0(ds,dz,du,dv) \cr
 \ar=:\ar
\bar{Z}_{k,1}(t,x)+\bar{Z}_{k,2}(t,x).
 \eeqnn
By (1.6) in \cite{saint} and the fact $u^p\le u^{\bar{p}}+1$ for $u\ge0$,
for $\alpha<\hat{p}<2$ we have
 \beqnn
\mbf{E}\Big\{ |\bar{Z}_{k,1}(t,x)|^{\hat{p}}  \Big\}
 \ar\le\ar
C\int_0^1 z^{\hat{p}} m_0(dz)\mbf{E}\Big\{\int_0^tds
\int_{\mbb{R}}H(X_{s-}(u))^\alpha p_{t-s}(x-u)^{\hat{p}}
1_{\{s\le\tilde{\tau}_k\}} du \Big\} \cr
 \ar\le\ar
C\,\mbf{E}\Big\{\int_0^t(t-s)^{-\frac{\hat{p}-1}{2}}ds
\int_{\mbb{R}}[1+X_s(u)^p] p_{t-s}(x-u) 1_{\{s\le\tilde{\tau}_k\}}
du \Big\}  \cr
 \ar\le\ar
C\,\mbf{E}\Big\{\int_0^t(t-s)^{-\frac{\hat{p}-1}{2}}ds
\int_{\mbb{R}}[1+X_s(u)^{\bar{p}}] p_{t-s}(x-u) 1_{\{s\le\tilde{\tau}_k\}}
du \Big\}
 \eeqnn
and
 \beqnn
\mbf{E}\Big\{ |\bar{Z}_{k,2}(t,x)|^{\bar{p}}  \Big\}
 \ar\le\ar
C\int_1^\infty z^{\bar{p}} m_0(dz)\mbf{E}\Big\{\int_0^tds \int_{\mbb{R}}
H(X_{s-}(u))^\alpha p_{t-s}(x-u)^{\bar{p}} 1_{\{s\le\tilde{\tau}_k\}} du
\Big\} \cr
 \ar\le\ar
C\,\mbf{E}\Big\{\int_0^t(t-s)^{-\frac{\bar{p}-1}{2}}ds
\int_{\mbb{R}}[1+X_s(u)^p] p_{t-s}(x-u) 1_{\{s\le\tilde{\tau}_k\}}
du \Big\} \cr
 \ar\le\ar
C\,\mbf{E}\Big\{\int_0^t(t-s)^{-\frac{\bar{p}-1}{2}}ds
\int_{\mbb{R}}[1+X_s(u)^{\bar{p}}] p_{t-s}(x-u) 1_{\{s\le\tilde{\tau}_k\}}
du \Big\}.
 \eeqnn
Then one obtains that
 \beqnn
 \ar\ar
\mbf{E}\{|\bar{Z}_k(t,x)|^{\bar{p}}\}
 \le
2\mbf{E}\{|\bar{Z}_{k,1}(t,x)|^{\bar{p}}\}
+2\mbf{E}\{|\bar{Z}_{k,2}(t,x)|^{\bar{p}}\} \cr
 \ar\ar\qquad\le
C\Big\{\mbf{E}[|\bar{Z}_{k,1}(t,x)|^{\hat{p}} ]
+ \mbf{E}[|\bar{Z}_{k,2}(t,x)|^{\bar{p}}]+1\Big\} \cr
 \ar\ar\qquad\le
C\,\mbf{E}\Big\{\int_0^t[(t-s)^{-\frac{\hat{p}-1}{2}}+(t-s)^{-\frac{\bar{p}-1}{2}}]ds
\int_{\mbb{R}}[1+X_s(u)^{\bar{p}}] p_{t-s}(x-u)
1_{\{s\le\tilde{\tau}_k\}} du \Big\}+C.
 \eeqnn

Combining this with \eqref{1.3} and \eqref{1.23}, we have
 \beqlb\label{1.40}
 \ar\ar
\mbf{E}\Big\{\int_0^Tdt\int_{\mbb{R}}X_t(y)^{\bar{p}} p_{T-t}(x-y)1_{\{t\le
\tilde{\tau}_k\}} dy\Big\} \cr
 \ar\le\ar
3\int_0^Tdt\int_{\mbb{R}}X_0(p_t(y-\cdot))^{\bar{p}} p_{T-t}(x-y)dy
+3 \int_0^Tdt\int_{\mbb{R}}\mbf{E}\{|\bar{Z}_k(t,y)|^{\bar{p}}\} p_{T-t}(x-y) dy  \cr
 \ar\ar
+3\,\mbf{E}\Big\{\int_0^Tdt\int_{\mbb{R}}
\Big|\int_0^tds\int_{\mbb{R}}p_{t-s}(y-u)G(X_s(u))du\Big|^{\bar{p}}
p_{T-t}(x-y)1_{\{t\le \tilde{\tau}_k\}} dy\Big\} \cr
 \ar\le\ar
C\int_0^Tdt\int_0^t[t^{\bar{p}}+(t-s)^{-\frac{\hat{p}-1}{2}}+(t-s)^{-\frac{\bar{p}-1}{2}}]
\mbf{E}\Big\{\int_{\mbb{R}}[1+X_s(u)^{\bar{p}}] p_{T-s}(x-u)
1_{\{s\le\tilde{\tau}_k\}} du \Big\}ds \cr
 \ar\ar
+CX_0(1)^{\bar{p}}T^{\frac{2-\bar{p}}{2}}+CT \cr
 \ar=\ar
C\,\mbf{E}\Big\{\int_0^Tds\int_{\mbb{R}}[1+X_s(u)^{\bar{p}}] p_{T-s}(x-u)
1_{\{s\le\tilde{\tau}_k\}} du
\int_s^T[t^{\bar{p}}+(t-s)^{-\frac{\hat{p}-1}{2}}+(t-s)^{-\frac{\bar{p}-1}{2}}]
dt\Big\} \cr
 \ar\ar
+CX_0(1)^{\bar{p}}T^{\frac{2-\bar{p}}{2}}+CT \cr
 \ar\le\ar
C(T^{\bar{p}+1}+T^{\frac{3-\hat{p}}{2}}+T^{\frac{3-\bar{p}}{2}})
\mbf{E}\Big\{\int_0^Tds\int_{\mbb{R}}[1+X_s(u)^{\bar{p}}] p_{T-s}(x-u)
1_{\{s\le\tilde{\tau}_k\}} du \Big\} \cr
 \ar\ar
+CX_0(1)^{\bar{p}}T^{\frac{2-\bar{p}}{2}}+CT.
 \eeqlb
In view of \eqref{1.32} and \eqref{1.17},
 \beqnn
\mbf{E}\Big\{\int_0^Tdt\int_{\mbb{R}}X_t(y)^{\bar{p}} p_{T-t}(x-y)1_{\{t\le
\tilde{\tau}_k\}} dy\Big\} <\infty.
 \eeqnn
Taking $\tilde{T}_0>0$ satisfying $K':=C(\tilde{T}_0^{\bar{p}+1}+\tilde{T}_0^{\frac{3-\hat{p}}{2}}+\tilde{T}_0^{\frac{3-\bar{p}}{2}})<1$,
for all $T\in[0,\tilde{T}_0]$ and $k\ge1$ we have
 \beqnn
\mbf{E}\Big\{\int_0^Tdt\int_{\mbb{R}}X_t(y)^{\bar{p}} p_{T-t}(x-y)1_{\{t\le
\tilde{\tau}_k\}} dy\Big\} \le
(1-K')^{-1}\Big[CX_0(1)^{\bar{p}}T^{\frac{2-\bar{p}}{2}}+CT\Big].
 \eeqnn
Then by the monotone convergence theorem
 \beqlb\label{3.1}
\sup\limits_{T\in[0,\tilde{T}_0]}\mbf{E}\Big\{\int_0^Tdt
\int_{\mbb{R}}X_t(x)^{\bar{p}} p_{T-t}(x-y) dx\Big\}<\infty.
 \eeqlb

{\bf Step 2.}
In this step we prove that \eqref{1.22} holds with $T$ replaced by the $\tilde{T}_0$
specified in Step 1.
Observe that for $0<r<1$,
 \beqlb\label{3.2}
 \ar\ar
\int_0^T(T-t)^{-\frac{r}{2}}dt\int_{\mbb{R}}
X_0(p_t(y-\cdot))^{\bar{p}} p_{T-t}(x-y)dy \cr
 \ar\ar\quad\le
[(2\pi)^{-1}X_0(1)]^{\bar{p}} \int_0^T(T-t)^{-\frac{r}{2}} t^{-\frac{{\bar{p}}}{2}}dt
\le
[(2\pi)^{-1}X_0(1)]^{\bar{p}} T^{1-\frac{r}{2}-\frac{{\bar{p}}}{2}}dt.
 \eeqlb
For $r\in(0,1)$ and $\delta\in[1,2)$,
 \beqlb\label{3.3}
 \ar\ar
\int_0^T(T-t)^{-\frac{r}{2}}dt\int_{\mbb{R}} p_{T-t}(x-y)dy\int_0^tds
\int_{\mbb{R}} X_s(u)^{\bar{p}} p_{t-s}(y-u)^\delta du \cr
 \ar\ar\qquad\le
C\int_0^T(T-t)^{-\frac{r}{2}}dt\int_0^t(t-s)^{-\frac{\delta-1}{2}}ds
\int_{\mbb{R}} X_s(u)^{\bar{p}} p_{T-s}(x-u) du \cr
 \ar\ar\qquad=
C\int_0^Tds\int_{\mbb{R}} X_s(u)^{\bar{p}} p_{T-s}(x-u) du
\int_s^T(T-t)^{-\frac{r}{2}}(t-s)^{-\frac{\delta-1}{2}}dt \cr
 \ar\ar\qquad\le
CT^{\frac{3-r-\delta}{2}}\int_0^Tds\int_{\mbb{R}}
X_s(u)^{\bar{p}} p_{T-s}(x-u) du.
 \eeqlb
Similar to the argument in \eqref{1.40}, combining \eqref{1.3} and \eqref{3.1}--\eqref{3.3}, it is easy
to see that for $0<r<1$,
 \beqlb\label{3.6}
\sup\limits_{T\in[0,\tilde{T}_0]}T^{\frac{r}{2}}\mbf{E}\Big\{\int_0^T (T-t)^{-\frac{r}{2}}dt
\int_{\mbb{R}}X_t(y)^{\bar{p}} p_{T-t}(x-y)dy\Big\}<\infty.
 \eeqlb

By (1.6) of \cite{saint} again we have
 \beqnn
 \ar\ar
\mbf{E}\Big\{\Big|\int_0^t \int_0^1\int_{\mbb{R}}
\int_0^{H(X_{s-}(u))^\alpha }
z
p_{t-s}(x-u)\tilde{N}_0(ds,dz,du,dv)  \Big|^{\hat{p}}\Big\} \cr
 \ar\ar\qquad\le
C\int_0^1z^{\hat{p}}m_0(dz)\,\mbf{E}\Big\{\int_0^tds
\int_{\mbb{R}}H(X_{s-}(u))^\alpha p_{t-s}(x-u)^{\hat{p}}du  \Big\} \cr
 \ar\ar\qquad\le
C\int_0^1z^{\hat{p}}m_0(dz)\,\mbf{E}\Big\{\int_0^tds
\int_{\mbb{R}}[1+X_s(u)^p] p_{t-s}(x-u)^{\hat{p}}du  \Big\} \cr
 \ar\ar\qquad\le
C\,\mbf{E}\Big\{\int_0^t(t-s)^{-\frac{\hat{p}-1}{2}}ds
\int_{\mbb{R}}[1+X_s(u)^{\bar{p}}] p_{t-s}(x-u)du  \Big\}
 \eeqnn
for $\alpha<\hat{p}<2$
and
 \beqnn
 \ar\ar
\mbf{E}\Big\{\Big|\int_0^t \int_1^\infty\int_{\mbb{R}}
\int_0^{H(X_{s-}(u))^\alpha }
z
p_{t-s}(x-u)\tilde{N}_0(ds,dz,du,dv)  \Big|^{\bar{p}} \Big\} \cr
 \ar\ar\qquad\le
C\int_1^\infty z^{\bar{p}}m_0(dz)\,\mbf{E}\Big\{\int_0^tds\int_{\mbb{R}}
H(X_{s-}(u))^\alpha p_{t-s}(x-u)^{\bar{p}}du \Big\}\cr
 \ar\ar\qquad\le
C\int_1^\infty z^{\bar{p}}m_0(dz)\,\mbf{E}\Big\{\int_0^tds\int_{\mbb{R}}
[1+X_s(u)^{\bar{p}}] p_{t-s}(x-u)^{\bar{p}}du \Big\} \cr
 \ar\ar\qquad\le
C\,\mbf{E}\Big\{\int_0^t(t-s)^{-\frac{\bar{p}-1}{2}}ds\int_{\mbb{R}}
[1+X_s(u)^{\bar{p}}] p_{t-s}(x-u)du \Big\},
 \eeqnn
which implies
 \beqnn
 \ar\ar
\mbf{E}\Big\{\Big|\int_0^t \int_0^\infty\int_{\mbb{R}}
\int_0^{H(X_{s-}(u))^\alpha }
z p_{t-s}(x-u)\tilde{N}_0(ds,dz,du,dv)  \Big|^{\bar{p}} \Big\} \cr
 \ar\ar\quad\le
C\,\mbf{E}\Big\{\int_0^t
\big[(t-s)^{-\frac{\bar{p}-1}{2}}+(t-s)^{-\frac{\hat{p}-1}{2}}\big]ds
\int_{\mbb{R}}p_{t-s}(x-y)[1+X_s(y)^{\bar{p}}]dy\Big\}+ C.
 \eeqnn
By \eqref{1.3} and \eqref{1.23} again, we have
 \beqlb\label{3.7}
\mbf{E}\{X_t(x)^{\bar{p}}\}
 \ar\le\ar
Ct^{-\frac{\bar{p}}{2}}+
C\,\mbf{E}\Big\{\Big|\int_0^tds\int_{\mbb{R}}
p_{t-s}(x-z)G(X_s(z))dz\Big|^{\bar{p}}\Big\} \cr
 \ar\ar
+C\,\mbf{E}\Big\{\Big|\int_0^t \int_0^\infty\int_{\mbb{R}}
\int_0^{H(X_{s-}(u))^\alpha }
z
p_{t-s}(x-u)\tilde{N}_0(ds,dz,du,dv)  \Big|^{\bar{p}} \Big\} \cr
 \ar\le\ar
C\,\mbf{E}\Big\{\int_0^t
\big[t^{\bar{p}}+(t-s)^{-\frac{\bar{p}-1}{2}}+(t-s)^{-\frac{\hat{p}-1}{2}}\big]ds
\int_{\mbb{R}}p_{t-s}(x-y)[1+X_s(y)^{\bar{p}}]dy\Big\} \cr
 \ar\ar
+Ct^{-\frac{\bar{p}}{2}}+C.
 \eeqlb
Then by \eqref{3.6} one sees that \eqref{1.22} holds with $T$
replaced by $\tilde{T}_0$.

{\bf Step 3.}
Similar to Step 1, for $\tilde{\gamma}\in(0,1)$
and $0\le\tilde{T}_1\le \tilde{T}_0\wedge \tilde{T}_2$ with
$\tilde{\gamma}^{\bar{p}+1} \tilde{T}_2^{\bar{p}}(\tilde{T}_2-\tilde{T}_1)\le\tilde{T}_0^{\bar{p}+1}$,
 \beqlb\label{1.41}
 \ar\ar
\mbf{E}\Big\{\int_{\tilde{\gamma}\tilde{T}_1}^{\tilde{\gamma}\tilde{T}_2}dt\int_{\mbb{R}}X_t(y)^{\bar{p}}
p_{\tilde{\gamma}\tilde{T}_2-t}(x-y)1_{\{t\le \tilde{\tau}_k\}} dy\Big\} \cr
 \ar\le\ar
C\int_{\tilde{\gamma}\tilde{T}_1}^{\tilde{\gamma}\tilde{T}_2}dt
\int_0^t[t^{\bar{p}}+(t-s)^{-\frac{\hat{p}-1}{2}}
+(t-s)^{-\frac{\bar{p}-1}{2}}]
\mbf{E}\Big\{\int_{\mbb{R}}[1+X_s(u)^{\bar{p}}] \cr
 \ar\ar
\times p_{\tilde{\gamma}\tilde{T}_2-s}(x-u)
1_{\{s\le\tilde{\tau}_k\}} du \Big\}ds
+CX_0(1)^{\bar{p}}\tilde{T}_2^{\frac{2-\bar{p}}{2}}+C\tilde{T}_2 \cr
 \ar\le\ar
C\int_{\tilde{\gamma}\tilde{T}_1}^{\tilde{\gamma}\tilde{T}_2}
\mbf{E}\big\{\int_{\mbb{R}}[1+X_s(u)^{\bar{p}}] p_{\tilde{\gamma}\tilde{T}_2-s}(x-u)
1_{\{s\le\tilde{\tau}_k\}} du \Big\}ds \cr
 \ar\ar\qqquad\qqquad
\times\int_s^{\tilde{\gamma}\tilde{T}_2}[t^{\bar{p}}+(t-s)^{-\frac{\hat{p}-1}{2}}
+(t-s)^{-\frac{\bar{p}-1}{2}}]dt \cr
 \ar\ar
+C\int_0^{\tilde{\gamma}\tilde{T}_1} ds
\int_{\mbb{R}}[1+\mbf{E}\{X_s(u)^{\bar{p}}\}] p_{\tilde{\gamma}\tilde{T}_2-s}(x-u)du \int_{\tilde{\gamma}\tilde{T}_1}^{\tilde{\gamma}\tilde{T}_2}
[t^{\bar{p}}+(t-s)^{-\frac{\hat{p}-1}{2}}
+(t-s)^{-\frac{\bar{p}-1}{2}}]dt \cr
 \ar\ar
+CX_0(1)^{\bar{p}}\tilde{T}_0^{\frac{2-\bar{p}}{2}}+C\tilde{T}_0 \cr
 \ar\le\ar
C[(\tilde{\gamma}\tilde{T}_2)^{\bar{p}} \tilde{\gamma}(\tilde{T}_2-\tilde{T}_1)+\tilde{T}_0^{\frac{3-\hat{p}}{2}}
+\tilde{T}_0^{\frac{3-\bar{p}}{2}}]
\mbf{E}\Big\{\int_{\tilde{\gamma}\tilde{T}_1}^{\tilde{\gamma}\tilde{T}_2}ds
\int_{\mbb{R}}[1+X_s(u)^{\bar{p}}] p_{\tilde{\gamma}\tilde{T}_2-s}(x-u)
1_{\{s\le\tilde{\tau}_k\}} du\Big\} \cr
 \ar\ar
+
C[\tilde{\gamma}(1+\tilde{\gamma})^{\bar{p}} \tilde{T}_0^{\bar{p}+1}
+\tilde{T}_0^{\frac{3-\hat{p}}{2}}+\tilde{T}_0^{\frac{3-\bar{p}}{2}}]
\int_0^{\tilde{T}_1}ds\int_{\mbb{R}}[1+s^{-\frac{\bar{p}}{2}}] p_{\tilde{\gamma}\tilde{T}_2-s}(x-u)du
+C_{\tilde{T}_0} \cr
 \ar\le\ar
C[\tilde{\gamma}^{\bar{p}+1} \tilde{T}_2^{\bar{p}}(\tilde{T}_2-\tilde{T}_1)+\tilde{T}_0^{\frac{3-\hat{p}}{2}}
+\tilde{T}_0^{\frac{3-\bar{p}}{2}}]
\mbf{E}\Big\{\int_{\tilde{\gamma}\tilde{T}_1}^{\tilde{\gamma}\tilde{T}_2}ds
\int_{\mbb{R}}[1+X_s(u)^{\bar{p}}] \cr
 \ar\ar\qqquad
\times p_{\tilde{\gamma}\tilde{T}_2-s}(x-u)
1_{\{s\le\tilde{\tau}_k\}} du \Big\} +C_{\tilde{T}_0}\cr
 \ar\le\ar
K'\,
\mbf{E}\Big\{\int_{\tilde{\gamma}\tilde{T}_1}^{\tilde{\gamma}\tilde{T}_2}ds
\int_{\mbb{R}}[1+X_s(u)^{\bar{p}}] p_{\tilde{\gamma}\tilde{T}_2-s}(x-u)
1_{\{s\le\tilde{\tau}_k\}} du \Big\} +C_{\tilde{T}_0},
 \eeqlb
where the assertion in Step 2 was used in the third inequality.
This implies
 \beqlb\label{1.42}
\sup\limits_{\tilde{T}_1\in[0,\tilde{T}_0],\tilde{\gamma}^{\bar{p}+1} \tilde{T}_2^{\bar{p}}(\tilde{T}_2-\tilde{T}_1)\le\tilde{T}_0^{\bar{p}+1}}
\mbf{E}\Big\{\int_{\tilde{\gamma}\tilde{T}_1}^{\tilde{\gamma}\tilde{T}_2}dt
\int_{\mbb{R}}X_t(y)^{\bar{p}} p_{\tilde{\gamma}\tilde{T}_2-t}(x-y) dy\Big\} <\infty.
 \eeqlb
Then by the assertion in Step 2 again,
 \beqnn
\sup\limits_{\tilde{T}_1\in[0,\tilde{T}_0],\tilde{\gamma}^{\bar{p}+1} \tilde{T}_2^{\bar{p}}(\tilde{T}_2-\tilde{T}_1)\le\tilde{T}_0^{\bar{p}+1}}
\mbf{E}\Big\{\int_0^{\tilde{\gamma}\tilde{T}_2}dt
\int_{\mbb{R}}X_t(y)^{\bar{p}} p_{\tilde{\gamma}\tilde{T}_2-t}(x-y) dy\Big\} <\infty.
 \eeqnn
Similar to \eqref{3.6} we have
 \beqlb\label{3.8}
\sup\limits_{\tilde{T}_1\in[0,\tilde{T}_0],\tilde{\gamma}^{\bar{p}+1} \tilde{T}_2^{\bar{p}}(\tilde{T}_2-\tilde{T}_1)\le\tilde{T}_0^{\bar{p}+1}}
(\tilde{\gamma}\tilde{T}_2)^{\frac{r}{2}}
\mbf{E}\Big\{\int_0^{\tilde{\gamma}\tilde{T}_2} (\tilde{\gamma}\tilde{T}_2-t)^{-\frac{r}{2}}dt
\int_{\mbb{R}}p_{\tilde{\gamma}\tilde{T}_2-t}(x-y)X_t(y)^{\bar{p}}dy\Big\}<\infty
 \eeqlb
for $r\in(0,1)$.
This together this with \eqref{3.7} shows that \eqref{1.22} holds
with $T$ replaced by $\tilde{\gamma}  \tilde{T}_2$,
where $\tilde{\gamma}^{\bar{p}+1} \tilde{T}_2^{\bar{p}}(\tilde{T}_2-\tilde{T}_0)\le\tilde{T}_0^{\bar{p}+1}$.

{\bf Step 4.}
Since $1+\frac12+\frac13+\cdots+\frac1n\le 1+\ln n$, one can chose
$\tilde{\gamma}\in(0,1)$ so that
$\sup_{n\ge1}\tilde{\gamma}^{\bar{p}+1}(1+\ln n)^{\bar{p}}/n\le1$, which implies
 \beqnn
\sup_{n\ge1}\frac{\tilde{\gamma}^{\bar{p}+1}}{n}
\big[1+\frac12+\frac13+\cdots+\frac1n\big]^{\bar{p}}\le1.
 \eeqnn
Observe that Step 2 proves that \eqref{1.22} holds
with $T$ replaced by $\tilde{\gamma} \tilde{T}_0$.
With $\tilde{T}_1$ and $\tilde{T}_2$ replaced by $\tilde{T}_1', \, 0\le\tilde{T}_1' \le \tilde{T}_0 $ and
$(1+\frac12)\tilde{T}_1'$, respectively, in \eqref{1.41}--\eqref{1.42},
we get \eqref{1.22} with $T$ replaced by
$\tilde{\gamma}(1+\frac12)\tilde{T}_0$.
Repeating the above argument, for each $n\ge1$,
with $\tilde{T}_1$ and $\tilde{T}_2$ replaced by
$(1+\frac12+\frac13+\cdots+\frac{1}{n-1})\tilde{T}_1'$ and
$(1+\frac12+\frac13+\cdots+\frac{1}{n})\tilde{T}_1'$, respectively, in \eqref{1.41}--\eqref{1.42}, we can get \eqref{1.22} with $T$ replaced by $\tilde{\gamma}(1+\frac12+\frac13+\cdots+\frac1n)\tilde{T}_0$, which completes the proof.
 \qed

\section*{Acknowledgments}
The first author is deeply grateful to Professors Zenghu Li and Hui
He for their encouragement and helpful discussions. The first author
also thanks Concordia University where the early version of this
paper was finished during his visit. The two authors would like
thank the anonymous referees and editor for the careful reading of
the manuscript, and for a number of useful comments and suggestions
that have improved the presentation of the paper.


\end{document}